\documentclass[11 pt,a4paper,leqno]{article}
\usepackage[all]{xy}
\usepackage {makeidx}
\usepackage[french,english]{babel}
\usepackage{stmaryrd}

\usepackage{xcolor}
\usepackage{tikz}
\usetikzlibrary{patterns}

\makeindex

\RequirePackage{amsmath}
\RequirePackage{amssymb}
\RequirePackage{ifthen}
\RequirePackage{theorem}
\RequirePackage{pslatex}

\usepackage{amsmath}
\usepackage{amssymb}
\usepackage{ifthen}
\usepackage{theorem}
\usepackage{pslatex}

\title{\sc Geometry of the space of monodromy data}

\author{Jean-Pierre Ramis
\footnote{Institut de France (Acad\'emie des Sciences) and
  Institut de Math\'ematiques de Toulouse,
  CNRS UMR 5219, Universit\'e Paul Sabatier (Toulouse 3),
  118 route de Narbonne, 31062 Toulouse CEDEX 9, France;
  E-mail: ramis.jean-pierre@wanadoo.fr},
Jacques Sauloy
\footnote{Toulouse; E-mail:jacquessauloy@gmail.com}}

\theorembodyfont{\sl} 
\newtheorem{thm}{Theorem}[section]
\newtheorem{lem}[thm]{Lemma}
\newtheorem{prop}[thm]{Proposition}
\newtheorem{cor}[thm]{Corollary}

\theorembodyfont{\rm}
\newtheorem{rmk}[thm]{Remark}
\newtheorem{rems}[thm]{Remarks}
\newtheorem{exa}[thm]{Example}
\newtheorem{exas}[thm]{Examples}

\newtheorem{defn}[thm]{Definition}

\numberwithin{equation}{thm}

\def\Proof{\textsl{Proof. - }}
\def\Finpr{\hfill $\Box$ \\}
\def\Finprcourt{\hfill $\Box$}

\def\N{{\mathbf N}}                      % ensemble des entiers naturels
\def\Z{{\mathbf Z}}                      % ensemble des entiers relatifs
\def\C{{\mathbf C}}                      % ensemble des complexes
\def\Cs{{\mathbf C^*}}                   % ensemble des complexes non nuls

\def\Id{{\text{Id}}}                     % Identite
\def\Ker{{\text{Ker}}}                   % Noyau
\def\Sp{{\text{Sp}}}                     % Spectre (algèbre linéaire)
\def\ii{{\text{i}}}                      % imaginaire
\def\Diag{{\text{Diag}}}                 % matrices diagonales

\def\M{{\mathcal{M}}}                    % corps des fonctions meromorphes
\def\Linc{{\text{Lin}_\C}}
\def\GL{{\text{GL}}}                     % groupe lineaire
\def\GLn{{\text{GL}_{n}}}                 % groupe lineaire a l'ordre n
\def\SLn{{\text{SL}_{n}}}                 % groupe special lineaire a l'ordre n
\def\Mat{{\text{Mat}}}                   % ensemble des matrices
\def\Matn{{\text{Mat}_{n}}}               % ensemble des matrices a l'ordre n
\def\GLnc{\GLn(\C)}                      % groupe lineaire sur C
\def\SLnc{\SLn(\C)}                      % groupe spécial lineaire sur C
\def\Matnc{\Matn(\C)}                    % ensemble des matrices sur C
\def\Dnc{{\text{D}_n(\C)}}               % matrices diagonales inversibles d'ordre n
\def\Ddc{{\text{D}_2(\C)}}               % matrices diagonales inversibles d'ordre 2
\def\ie{{\emph{i.e.}}}
\def\eg{{\emph{e.g.}}}
\def\tq{\mid}
\def\Tq{\bigm|}

\def\cf{{\emph{cf.}}}
\def\Cf{{\emph{Cf.}}}

\def\Kg{{\C(x)}}
\def\Rwg{{\mathcal{O}(\C^{*})}}

\def\thq{{\theta_{q}}}                            % Notre fonction theta preferee
\def\sq{\sigma_q}                                 % dilatation f(z) \mapsto f(qz)
\def\Eq{{\mathbf{E}_{q}}}                         % courbe elliptique \C^*/q^\Z
\def\x{{\underline{x}}}

\def\E{{\mathcal{E}}}                   
\def\F{{\mathcal{F}}}               
\def\G{{\mathcal{G}}}

\def\iff{{if and only if}}

\def\Xf{{\mathbf{X}}}
\def\Pr{{\mathbf{P}}}
\def\e{{\underline{e}}}
\def\f{{\underline{f}}}
\def\u{{\underline{u}}}
\def\v{{\underline{v}}}
\def\w{{\underline{w}}}

\def\Affi{{\mathbf{A}}}
\def\Aff{{\C}}
\def\Af{{\C}}
\def\T{{\mathbf{T}}}
\def\Tn{{\mathbf{T}_n}}
\def\Spec{{\text{Spec}}}    % spectre (géométrie algébrique affine)

\def\O{{\mathcal{O}}}

\def\Bs{{\mathcal{B^*}}}
\def\cl{{\text{cl}}}
\def\zero{{\underline{0}}}
\def\Uu{{\mathcal{U}}}

\begin{document}

\maketitle

%%%%%%%%%%%%%%%%%%%%%%%%%%%%%%%%%%%%%%%%%%%%%%%%%%%%%%%%%%%%%%%%%%%%%%%%%%%%%

\begin{abstract}
  In \cite{ORS}, with Yousuke Ohyama, we defined and studied a ``space of monodromy data''
  underlying the well known derivation of $q$-Painlevé VI equation from ``$q$-isomonodromy''
  conditions by Jimbo and Sakai \cite{JimboSakai}. In \cite{JR}, Nalini Joshi and Pieter
  Roffelsen pursued our work. However, both \cite{ORS} and \cite{JR} are ambiguous on some
  foundational algebro-geometric matters. We proceed here to provide sound bases.
\end{abstract}

\selectlanguage{french}

\begin{abstract}
  Dans \cite{ORS}, avec Yousuke Ohyama, nous avons défini et étudié un ``espace des données
  de monodromie'' sous-jacent à la célèbre déduction de l'équation de $q$-Painlevé VI à partir
  de conditions de ``$q$-isomonodromie'' par Jimbo et Sakai \cite{JimboSakai}. Dans \cite{JR},
  Nalini Joshi and Pieter Roffelsen ont prolongé notre travail. Cependant, aussi bien \cite{ORS}
  que \cite{JR} présentent des ambiguïtés sur certains aspects de leurs fondements
  algébro-géométriques. Nous en proposons ici des bases rigoureuses. 
\end{abstract}

\selectlanguage{english}

\tableofcontents

%%%%%%%%%%%%%%%%%%%%%%%%%%%%%%%%%%%%%%%%%%%%%%%%%%%%%%%%%%%%%%%%%%%%%%%%%%%%%

% 0

\setcounter{section}{-1}

\section{Introduction}

In \cite{ORS} we attempted a description of the right-hand side of the Riemann-Hilbert correspondence
for a particular family of $q$-difference equations related to $q$-Painlev\'e VI. We did that using
a variant of the usual Birkhoff connection matrix, in which the \emph{local contributions} at $0$
and $\infty$ are ripped off; the interest of such a procedure had been demonstrated (with somewhat
different motivations) in \cite{JSENS}, where the idea was first introduced. In the first part of
\cite{ORS}, we studied the space of such matrices (up to relevant equivalence relations) from the
point of view of algebraic geometry, using mainly tools of (bi)linear algebra\footnote{In the second
part of \emph{loc. cit.} we introduced another approach through so-called \emph{Mano decomposition}.
This aspect will appear here only in a particular case.}. However, we did not reach a complete
description of the space of interest; in particular, we neither properly defined the quotient
structures involved, nor proved or disproved the smoothness properties that we were led to conjecture.
We were not either able to recognize among the standard list of classical algebraic surfaces (Segre,
del Pezzo, Kummer, \dots) a candidate model\footnote{We were not even able to decide if a good model
is a rational surface or not.}. \\

In \cite{JR}, Nalini Joshi and Pieter Roffelsen take up our angle of attack but develop further the
necessary calculations, to the point of solving some of the questions we had left unanswered. However,
we feel that some points stay unclear at the level of rigorous definitions: quotients are introduced
without any comment, let alone justification, about their nature (set theoretic ? analytic ? algebro-
geometric ?); the authors call (improperly) \emph{embedding} an injective map from the \emph{set} of
monodromy data to a projective space whose image is a locally closed algebraic sub-variety and do not
compare their different ``embeddings" from the point of view of algebraic geometry; etc. \\

In this work, we try to exploit the clever and efficient calculations of \cite{JR} while giving them
a sound basis in the domain of complex algebraic geometry. In particular, we properly define the
quotient spaces at stake and justify their existence and properties, and this allows us to give
rigorously a definition of the morphisms and a proof of their relevant properties. We define on
the \emph{set} of monodromy data a structure of geometric quotient $\F$ (this is precisely stated
in \ref{subsubsection:FtocalFquotientgeom}) and prove that
it is an affine, smooth, rational algebraic variety. We compare with $\F$ the different
structures on the \emph{set} of monodromy data \emph{introduced} in \cite{ORS} and \cite{JR}.
We get in particular an \emph{embedding}, \emph{in the sense of algebraic geometry}\footnote{That
is a morphism whose image is a locally closed subvariety isomorphic to the source.}, of $\F$
in $\C^4$, whose image is the \emph{affine Segre surface} introduced in \cite{JR}. Therefore this
affine Segre surface is a \emph{smooth} algebraic surface. \\

In \cite{ORS}, using the Mano decomposition and in relation with some properties of partial
reducibility, we described $16$ ``lines" on the space of monodromy data. We prove that they
correspond to the $16$ lines on the Segre surface. \\

Unfortunately, we cannot prove that $\F$ is isomorphic\footnote{We have only a bijective
morphism.} to the surface $\G$ defined in \cite{ORS} as an algebraic structure on the space
of monodromy data. Therefore our conjecture about the smoothness of $\G$, under some generic
hypothesis (\cf\ \cite[Conjecture 7.10]{ORS}), remains open\footnote{The assertion in
\cite[Remark 2.19]{JR}, is, in that sense, optimistic. What is proved in \cite{JR} is the existence
of a structure of a (possibly non separated) analytic manifold on the set of monodromy data: \cf\
Theorem 2.18. This does not imply our Conjecture 7.10.}. We return to the ``model $\G$'' in section
\ref{section:GModel}. \\

Although the heart of the article revolves around rank $2$ linear systems considered for generic
values of the parameters, we propose (as we did in \cite{ORS}) some more general results that could
allow for a future study of other cases.

\subsection*{Acknowledgements.}

The first author thanks Nalini Joshi and Pieter Roffelsen for interesting discussions about
the $16$ lines on a Segre surface. Both authors thank Yousuke Ohyama for sharing his knowledge
of Painlev\'e equations.

%%%%%%%%%%%%%%%%%%%%%%%%%%%%%%%%%%%%%%%%%%%%%%%%%%%%%%%%%%%%%%%%%%%%%%%%%%%%%

% 1

\section{Preparatory material}

% 1.1

\subsection{General preliminaries}

% 1.1.1

\subsubsection{General notations and conventions}

In all the text, if $E$ is a complex linear space, $E^\vee := \Linc(E,\C)$ denotes the dual of $E$.
We write $E^* := E \setminus \{0\}$. If $E$ is given as a product $\prod E_\alpha$, then we define
$E^{(*)} := \prod E_\alpha^*$ (hopefully, no ambiguity will arise). In a slight abuse, letting $\Pr(E)$
the projective space $E^*/\Cs$, we shall set $\Pr(E^{(*)}) := \prod \Pr(E_\alpha)$. \\

When a group $G$ operates (on the left) on $E$, we write $G x \subset E$ the orbit of $x \in E$ and
$G_x \subset G$ its stabilizer (or isotropy subgroup). \\

We write $P \sim Q$ for the conjugacy relation in $\GLnc$ (actually, the symbol $\sim$ will be used
for some more equivalence relations) and $\Sp\ P$ the spectrum of $P \in \Matnc$. We denote $\Dnc$
the subgroup of diagonal matrices in the linear group $\GLnc$. Diagonal matrices are abreviated as:
$$
\Diag(\lambda_1,\ldots,\lambda_n) := \begin{pmatrix}
  \lambda_1 & 0 & \ldots & 0 \\
  0 & \lambda_2 & \ldots & 0 \\
  \vdots & \vdots & \ddots & \vdots \\
  0 & 0 & \ldots & \lambda_n
  \end{pmatrix}.
$$
For even shorter abreviations, we sometimes identify $\Dnc$ with $\Cs^n$ and
$\Diag(\lambda_1,\ldots,\lambda_n) \in \Dnc$ with
$\underline{\lambda} := (\lambda_1,\ldots,\lambda_n) \in \Cs^n$. \\

For any two $c, d \in \Cs$, we write $c \equiv d$ the congruence relation modulo the subgroup $q^\Z$:
$$
\forall c,d \in \Cs \;,\; c \equiv d \underset{def}{\Longleftrightarrow}
c/d \in q^\Z := \{q^k \tq k \in \Z\} \subset \Cs,
$$
and $c \not\equiv d$ its negation. The words ``congruent'', ``congruence'', etc, applied to elements
of $\Cs$, will refer to that relation. \\

Whenever $A(x)$, $F(x)$ \dots, are invertible matrices of functions, we write $\sq F(x) := F(qx)$ and
$F[A] := (\sq F) A F^{-1}$ (``gauge transformation'')

% 1.1.2

\subsubsection{The context}
\label{subsubsection:thecontext}

In \cite{ORS} we studied the following situation; let $R := \Diag(\rho_1,\ldots,\rho_n) \in \Dnc$
and $S := \Diag(\sigma_1,\ldots,\sigma_n) \in \Dnc$ be fixed with the following \emph{strong non
resonancy} assumption:
$$
\forall i,j \in \{1,\ldots,n\} \;,\;
i \neq j \Longrightarrow
\left( \rho_i \not\equiv \rho_j \text{~and~} \sigma_i \not\equiv \sigma_j \right).
$$
Let also $\mu \in \N^*$, $N := \mu n$ and $x_1,\ldots,x_N \in \Cs$ be pairwise noncongruent; we write
$\x := \{x_1,\ldots,x_N\}$. Then we introduced the following sets:
$$
E_{R,S,\x} := \left\{A_0 + \cdots + A_\mu x^\mu \Tq
\begin{cases}
  \text{all~} A_i \in \Matnc, A_0,A_\mu \in \GLnc, \\
  A_0 \sim R, A_\mu \sim S, \\
  \det A(x_1) = \cdots = \det A(x_N) = 0,
  \end{cases}
\right\}
$$
and the quotient set $\E_{R,S,\x}$ of $E_{R,S,\x}$ by \emph{rational gauge equivalence relation} $\sim$:
$$
\E_{R,S,\x} := \dfrac{E_{R,S,\x}}{\sim}, \text{~where~}
\forall A,B \in E_{R,S,\x} \;,\; A \sim B \underset{def}{\Longleftrightarrow}
\exists F \in \GLn(\Kg) \;:\; B = F[A].
$$

Let $V := \{M \in \Matn(\Rwg) \tq \sq M = R M (S x^\mu)^{-1}\}$, a complex linear space of dimension
$\mu n^2$ (as we shall soon see); we also introduced its subset:
$$
F_{R,S,\x} := \{M \in V \tq \det M \neq 0, \det M \text{~vanishes on~} \x\}.
$$
The group $\Dnc \times \Dnc$ acts linearly on $V$ through the formula:
$$
(\Gamma,\Delta).M := \Gamma M \Delta^{-1}.
$$
The subset $F_{R,S,\x} \subset V$ is stable under that action; letting $\sim$ the induced equivalence
relation on it, we define the quotient set (later abreviated as $\F$, since the local data $R,S,\x$
will be fixed):
$$
\F_{R,S,\x} :=\dfrac{F_{R,S,\x}}{\sim} \cdot
$$

We defined (and proved) a bijection (``Riemann-Hilbert-Birkhoff correspondence'') from $\E_{R,S,\x}$
to $\F_{R,S,\x}$. Our main goal here is to complete the geometric description of $\F_{R,S,\x}$, mainly
in the case $n = \mu = 2$ (the so-called ``Jimbo-Sakai case'' or ``JS case'' for short). \\

It is an essential feature of the target (``right hand side'') of our correspondence that the
condition $\sq M = R M (S x^\mu)^{-1}$ on $M = (m_{i,j})_{1 \leq i,j \leq n} \in V \subset \Matn(\Rwg)$
splits into $n^2$ independent conditions $\sq m_{i,j} = (\rho_i/\sigma_j) x^{-\mu} m_{i,j}$, \ie\
$m_{i,j} \in V_{i,j}$, where $V_{i,j} := V_{\mu,\rho_i/\sigma_j}$, such spaces being defined as:
$$
\forall c \in \Cs \;,\; \forall k \in \N^* \;,\; V_{k,c} := \{m \in \Rwg \tq \sq m = c x^{-k} m\},
$$
each such linear space having dimension $\dim_\C V_{k,c} = k$. So we have a natural identification
$V := \prod\limits_{1 \leq i,j \leq n} V_{i,j}$ (whence our contention that $\dim_\C V = \mu n^2$). \\

Another essential feature is that the action of $\Dnc \times \Dnc$ on $F_{R,S,\x}$, which can naturally
extended to $V \supset F_{R,S,\x}$ as a linear action, splits correspondingly:
$$
(\Gamma,\Delta).(m_{i,j})_{1 \leq i,j \leq n} = \left((\gamma_i/\delta_j) m_{i,j}\right)_{1 \leq i,j \leq n}, 
$$
where $\Gamma =: \Diag(\gamma_1,\ldots,\gamma_n)$ and $\Delta =: \Diag(\delta_1,\ldots,\delta_n)$.
These two features will somehow be ``axiomatized'' in \ref{subsection:invariantsquotientsforgeneraln}.

% 1.2

\subsection{Some facts about quotients}

Since we intend to clarify the algebro-geometric structure of our spaces of interest, it seems fit
to make explicit our framework. We deal with the elementary theory of algebraic varieties over $\C$
(not necessarily affine, nor irreducible, nor separated but always reduced) and linear algebraic
groups (automatically smooth but possibly reducible). \\

A word about terminology: since we do \emph{not} require our algebraic varieties to be separated,
whenever they are, we call them \emph{separated varieties}. Of course all affine, quasi-affine,
projective and quasi-projective varieties are separated, as will be all varieties $X$ for which
we seek to construct a quotient. However, we shall in the end obtain some non separated quotients
$X \rightarrow Y$ (\ie\ $X$ is separated but $Y$ is not). At any rate, in the whole of
\ref{subsubsection:groupactions} and \ref{subsubsection:quotients}, \emph{all varieties are assumed
to be separated}. This convention will be relaxed from \ref{subsubsection:nonseparatedquotients} on
(except in explicitly stated special cases).

% 1.2.1

\subsubsection{Group actions}
\label{subsubsection:groupactions}

\paragraph{Basic formalism.}

Unless otherwise stated, the words ``closed'', ``open'', ``dense'', etc, refer to Zariski topology.
An affine algebraic variety $X$ is totally determined (and conversely) by its affine algebra
$\Aff[X]$ (a finite type reduced $\C$-algebra) and morphisms of affine varieties $X \rightarrow Y$
correspond functorially to morphisms of $\C$-algebras $\Aff[Y] \rightarrow \Aff[X]$. A general
morphism of varieties $\phi: X \rightarrow Y$ is locally given in the above form, \ie\ by morphisms
of $\C$-algebras $\Aff[V] \rightarrow \Aff[\phi^{-1}(V)]$ where the affine sets $V$ cover $Y$.
Important examples of \emph{non affine} varieties, to keep in mind, are $\C^n \setminus \{0\}$
and $\Pr^{n-1}(\C)$ (whenever $n \geq 2$). \\
All our algebraic groups are affine (as varieties), hence linear (\ie\ realized as closed subgroups
of some $\GLnc$); main references: \cite{Borel,Humphreys,SpringerLAG}, also see \cite{SpringerGAIV}.
For any morphism of algebraic groups $f: G \rightarrow H$ (\ie\ $f$ is a morphism of varieties and
a group morphism), the image $f(G) \subset H$ is a closed subgroup (the kernel $\Ker\ f \subset G$
obviously is !); and, if $f$ is injective, it is a closed immersion. The neutral component $G^\circ$
(\ie\ the connected component of the identity $1$) is a closed normal subgroup, the irreducible
components of $G$ are also its connected components $g G^\circ = G^\circ g$ and the quotient group
$G/G^\circ$ is finite. \\
From the structure theory of linear algebraic groups, we retain the following definitions and facts:
\begin{itemize}
\item The group $G$ is \emph{semi-simple} if its only connected normal solvable subgroup is trivial.
  The special linear groups $\SLnc$ and the tori $\Cs^n$ are semi-simple.
\item The group $G$ is \emph{reductive} if its only connected normal unipotent subgroup is trivial.
  The group $\GLnc$, and all semi-simple groups are reductive.
\end{itemize}   

\paragraph{Rational actions.}

Our main general references on algebraic group actions and quotients are \cite{Dolgachev,NewsteadTata}
and, for some special facts, \cite{PopovVinberg}; also note the survey \cite{Brion} and the short
summary \cite{NewsteadHAL} of \cite{NewsteadTata}. \\
Let $G$ an algebraic group and $X$ an algebraic variety. A \emph{rational action} of $G$ on $X$
is a morphism of varieties $G \times X \rightarrow X$, $(g,x) \mapsto gx$ which is a group action
(\ie\ $g'(gx) = (g'g)x$ and $1 x = x$). We shall just speak of an action and say that $X$ is a
\emph{$G$-variety}. Stabilizers $G_x$ of elements $x \in X$ are automatically closed subgroups
of $G$. Typical examples, to keep in mind, are:
\begin{enumerate}
\item $\Cs$ acting on $\C^n$, resp. on $\C^n \setminus \{0\}$, by homotheties;
\item $\Cs$ acting on $\C^2$, resp. on $\C^2 \setminus \{0\}$, by $t.(x,y) := (tx,t^{-1}y)$.
\end{enumerate}
An \emph{equivariant} morphism $\phi: X \rightarrow Y$ of $G$-varieties, also called a $G$-morphism,
is one such $\phi(gx) = g \phi(x)$ for all $g \in G$, $x \in X$. If $Y$ is endowed with the trivial
action (such that $g y = y$ for all $g \in G$, $y \in Y$) we say that $\phi$ is \emph{invariant}
(so $\phi(gx) = \phi(x)$ for all $g \in G$, $x \in X$). If $\phi$ is invariant, the induced morphism
of algebras $\Aff[U] \rightarrow \Aff[\phi^{-1}(U)]$ (generally defined for every morphism of varieties
$\phi: X \rightarrow Y$ and every open subset $U \subset Y$) has image in the subalgebra
$\Aff[\phi^{-1}(U)]^G$ of $G$-invariant functions. \\
A (rational) \emph{$G$-module} is a finite dimensional $\C$-linear space $V$ with a linear
$G$-action (\ie\ $x \mapsto gx$ is linear for all $g \in G$); equivalently, it is a linear
representation and the corresponding group morphism $G \rightarrow \GL(V)$ is a morphism of
algebraic groups. Then every algebraic subset of $V$ which is $G$-invariant is canonically a
$G$-variety. Moreover, every affine $G$-variety can be obtained in this way. (The definition
of $G$-module is usually extended to infinite dimensional $\C$-linear spaces $V$, provided
$V = \bigcup V_i$ where all $V_i$ are $G$-invariant and finite dimensional.)

\paragraph{Orbits.}

Let $X$ a $G$-variety. Then the \emph{orbit} $G x$ of any $x \in X$ is locally closed, \ie\
it is open in its Zariski closure $\overline{Gx}$. Hence the complementary subset
$\overline{Gx} \setminus Gx$ is a disjoint union of orbits, each having dimension $< \dim Gx$. \\
We have $\dim Gx = \dim G - \dim G_x$ and the function $x \mapsto \dim G x$ is lower semicontinuous,
\ie:
$$
\{x \in X \tq \dim G x \leq n\} \text{~is closed for every~} n \geq 0.
$$
Thus orbits of minimal dimension are closed and every orbit contains (at least) a closed orbit. \\
Note that if $X$ is a $G$-variety and $\phi: X \rightarrow Y$ an invariant morphism, every fiber
$\phi^{-1}(y)$ is closed and $G$-invariant, whence a union of the orbits $Gx$ such that $\phi(x) = y$.
So, the quotient map from the set of orbits\footnote{The notation is provisional and can be forgotten.
For quotients, we shall rather use $X/G$ and the like.} $\dfrac{X}{G}$ to $Y$ is well defined.
Looking for a \emph{separated} quotient (which a variety has to be), we would hope that map to be
bijective; this would at least require that all orbits be closed, which is seldom the case.

% 1.2.2

\subsubsection{Quotients}
\label{subsubsection:quotients}

\paragraph{Categorical quotients and orbit spaces.}

There are two natural ways to define the quotient of a $G$-variety $X$ by the action of $G$:
\begin{itemize}
\item Consider the quotient set $Y$, endow it with the quotient topology and the sheaf of
  $G$-invariant functions; this does not generally yield an algebraic variety.
\item Categorically: this will be detailed next; by definition, if it exists, such a quotient is
  an algebraic variety and unique up to isomorphism. But it may have bad geometrical properties.
\end{itemize}
The morphism of varieties $\phi: X \rightarrow Y$ is a \emph{categorical quotient} if it is
$G$-invariant (\ie\ constant on $G$-orbits) and initial for that property (\ie\ any $G$-invariant
morphism $\psi: X \rightarrow Z$ factors uniquely as $\psi = f \circ \phi$, where
$f: Y \rightarrow Z$ is a morphism). \\
In the case of $\Cs$ acting on $\C^n$ by homotheties, $0$ belongs to the closure of all orbits, so
all $G$-invariant $\psi: X \rightarrow Z$ are constant, so the categorical quotient is trivial
(a constant map to a point). \\
If the categorical quotient $\phi: X \rightarrow Y$ \emph{separates orbits}, \ie\ the fiber
$\phi^{-1}(y)$ is an orbit for every $y \in Y$, it is called an \emph{orbit space}. \\
In the case of $\Cs$ acting on $\C^n \setminus \{0\}$ by homotheties, the natural projection
$\C^n \rightarrow \Pr^{n-1}(\C)$ is an orbit space.

\paragraph{Good quotients and geometric quotients.}

The following definition comes from \cite{NewsteadTata}:
\begin{defn}
  A \emph{good quotient} is an affine\footnote{Recall that $\phi$ is said to be \emph{affine}
  if for each affine open subset $U \subset Y$, its preimage $\phi^{-1}(U) \subset X$ is affine.}
  morphism $\phi: X \rightarrow Y$ which satisfies the following properties: \\
  (i) $\phi$ is surjective. \\
  (ii) $\phi$ is $G$-invariant. \\  
  (iii) For every open subset $U \subset Y$, the induced morphism
  $\Aff[U] \rightarrow \Aff[\phi^{-1}(U)]^G$ is an isomorphism; in particular, $\Aff[Y]$ can be
  identified with $\Aff[X]^G$, which completely defines $Y$. \\
  (iv) For every closed $G$-invariant subset $W \subset X$, the set $\phi(W) \subset Y$ is closed. \\
  (v) For every pair of closed $G$-invariant subsets $W_1,W_2 \subset X$, if $W_1 \cap W_2 = \emptyset$
  then $\phi(W_1) \cap \phi(W_2) = \emptyset$. \\
  We write such a good quotient as $X//G$ (thus omitting the structural morphism $\phi$).
\end{defn}
The following properties are consequences:
\begin{cor}
  For every open subset $U \subset Y$, the morphism $\phi^{-1}(U) \rightarrow U$ is a categorical
  quotient. (In particular, $\phi:X \rightarrow Y$ is a categorical quotient.) If moreover the
  action of $G$ on $\phi^{-1}(U)$ is closed (\ie\ the orbits are closed), then it is an orbit space.
\end{cor}
Moreover,``$\phi$ separates orbits as much as topologically feasible'':
\begin{cor}
  Let $x_1,x_2 \in X$. Then
  $\phi(x_1) = \phi(x_2) \Leftrightarrow \overline{G x_1} \cap \overline{G x_2} \neq \emptyset$.
  (The implication $\Leftarrow$ is easy and holds for all invariant morphisms.) 
\end{cor}
\begin{cor}
  The fibers of $\phi: X \rightarrow X//G$ are closed if and only if they all have the same
  dimension, which in turn is equivalent to: the $\phi^{-1}(y)$ are the orbits, \ie\ $X//G$ is
  an orbit space.
\end{cor}
The following definition also comes from \cite{NewsteadTata}:
\begin{defn}
  \label{defn:geometricquotient}
  A \emph{geometric quotient} is a good quotient which is moreover an orbit space. \\
  We shall write a geometric quotients as $X/G$.
\end{defn}
The following characterization of a geometric quotient is taken as definition in
\cite{Brion,PopovVinberg}:
\begin{prop}
  \label{prop:geometricquotient}
  The morphism $\phi: X \rightarrow Y$ is a geometric quotient if and only if: \\
  (i) it is surjective and its fibers $\phi^{-1}(y)$, $y \in Y$, are exactly the $G$-orbits; \\
  (ii) $Y$ has the quotient topology, \ie\ $U \subset Y$ is open if and only if 
  $\phi^{-1}(U) \subset X$ is open; \\
  (iii) for every open subset $U \subset Y$, the induced morphism
  $\Aff[U] \rightarrow \Aff[\phi^{-1}(U)]^G$ is an isomorphism. \\
  Condition (ii) may be replaced by the equivalent one: (ii') $\phi$ is an open map.
\end{prop}
For example, $\C^n \setminus \{0\} \rightarrow \Pr^{n-1}(\C)$ (for $n \geq 2$) is a geometric
quotient, and so is $G \rightarrow G/H$ for any closed subgroup of an affine algebraic group $G$. \\
At any rate, we see that the geometric quotient $X/G$ is actually the topological quotient
$\pi: X \rightarrow \dfrac{X}{G}$ endowed with the sheaf $(\pi_* \mathcal{O}_X)^G$, where
$\mathcal{O}_X$ denotes the structural sheaf on the variety $X$ and $(\pi_* \mathcal{O}_X)^G$
the invariant subsheaf under the obvious action of $G$ on the direct image $\pi_* \mathcal{O}_X$.
Last, we quote from \cite[theorem 4.2]{PopovVinberg}:
\begin{thm}
  Let $X$ a $G$-variety and $\phi: X \rightarrow Y$ a surjective morphism such that its fibers
  $\phi^{-1}(y)$ are the orbits. Assume $X$ irreducible and $Y$ normal. Then $\phi$ is a geometric
  quotient.
\end{thm}

\paragraph{Action of reductive groups on affine varieties.}

Let $G$ a reductive group and $X$ an affine $G$-variety. The following are proved in
\cite{NewsteadTata}.
\begin{thm}
  There exist an affine variety $Y$ and a morphism $\phi:X \rightarrow Y$ which is a good quotient.
  It is the unique affine variety with affine algebra $\Aff[Y] := \Aff[X]^G$.
\end{thm}
The following precisions are proved in \cite{Brion}:
\begin{cor}
  The algebra of invariants functions $\Aff[X]^G$ is an affine algebra; let $f_1,\ldots,f_n$
  a set of generators. Then $X//G$ can be described as the image of the map $G \rightarrow \C^n$,
  $x \mapsto (f_1(x),\ldots,f_n(x))$, which is a closed subset of $\C^n$.
\end{cor}
  Now let $X' \subset X$ a closed $G$-invariant subset and let $\phi': X' \rightarrow X'//G$ as in
  the previous theorem. Since $\phi'$ is a categorical quotient and since the restriction
  $\phi_{|X'}: X' \rightarrow X//G$ is $G$-invariant, it factors through a morphism
  $X'//G \rightarrow X//G$.
\begin{cor}
  $X'//G \rightarrow X//G$ is a closed immersion, so $X'//G$ is (canonically identified to) a
  closed subset of $X//G$.
\end{cor}
\begin{cor}
  Let $X',X'' \subset X$ a pair of closed $G$-invariant subsets; then
  $\phi(X' \cap X'') = \phi(X') \cap \phi(X'')$.
\end{cor}
\begin{cor}
  Every fiber $\phi^{-1}(y)$ contains a unique closed orbit.
\end{cor}
So $X//G$ can be seen as ``the space of closed orbits''.
\begin{cor}
  If $X$ is irreducible, resp. normal, so is $X//G$.
\end{cor}

\paragraph{Stable points and geometric quotients.}

Let again $G$ a reductive group and $X$ an affine $G$-variety.
\begin{defn}
  \label{defn:stablepoints1}
  The point $x \in X$ is \emph{stable} if its orbit $Gx$ is closed and its stabilizer $G_x$ is finite.
\end{defn}
\begin{thm}
  \label{thm:stablepoints1}
  (i) The set $X^s$ of stable points is a (possibly empty) $G$-invariant open subset of $X$ and
  $\phi(X^s) \subset X//G$ is open too, so (provided $X^s$ is non empty) $\phi(X^s) = X^s//G$. \\
  (ii) The restriction $X^s \rightarrow \phi(X^s) = X^s//G$ is actually a geometric quotient $X^s/G$
  (again provided $X^s$ is non empty).
\end{thm}

% 1.2.3

\subsubsection{Non separated quotients}
\label{subsubsection:nonseparatedquotients}

We shall be confronted to the following situation: a separated $G$-variety $X$ can be covered by
$G$-invariant open subsets $X_i$ for which we are able to define geometric (resp. good, resp.
categorical) quotients $X_i \rightarrow Y_i$. Can we patch together the $Y_i$ to obtain a geometric
(resp. good, resp. categorical) quotient $X \rightarrow Y$ ? Under some circumstances, the answer
is yes, but $Y$ might be non separated. Since our basic way to construct geometric quotients involves
affine varieties, we shall only consider affine quotients $X_i \rightarrow Y_i$. 

\paragraph{Glueing of affine varieties.}

So we consider a family $(Y_i)_{i \in I}$ of affine varieties. Following \cite[\S 8.2]{Dolgachev},
we define \emph{glueing data} for that family as:
\begin{enumerate}
\item for each $(i,j) \in I \times I$, an affine open subset $U_{i,j} \subset Y_i$;
\item for each $(i,j) \in I \times I$, an isomorphism $f_{j,i}: U_{i,j} \rightarrow U_{j,i}$.
\end{enumerate}
Those are subject to the following compatibility conditions:
\begin{enumerate}
\item for each $i \in I$, $U_{i,i} = Y_i$ and $f_{i,i} = \Id_{Y_i}$;
\item for each $(i,j,k) \in I \times I \times I$, $f_{j,i}(U_{i,j} \cap U_{i,k}) \subset U_{j,k}$;
\item for each $(i,j,k) \in I \times I \times I$, $f_{k,j} \circ f_{j,i} = f_{k,i}$ on
  $U_{i,j} \cap U_{i,k}$.
\end{enumerate}
Note that restriction to the domain $U_{i,j} \cap U_{i,k}$ is a logical necessity in the third
condition, and that the second condition is then required for the third to have a meaning. \\
We can then define on the disjoint union $\bigsqcup Y_i$ an equivalence relation $\sim$ which is
the smallest such that $y_i \sim f_{j,i}(y_i)$ for each $(i,j) \in I \times I$ and $y_i \in U_{i,j}$.
Then the quotient set $\dfrac{\bigsqcup Y_i}{\sim}$ admits a natural structure of algebraic variety
(possibly non separated) such that:
\begin{enumerate}
\item each $Y_i$ embeds into $Y$, so we shall identify $Y_i$ as an affine open subset of $Y$;
\item a regular function $g$ on $Y$ is the same thing as a family of compatible regular functions
  $g_i$ on $Y_i$ (\ie\ such that $g_j \circ f_{j,i} = g_i$ on $U_{i,j}$).
\end{enumerate}
\begin{cor}
  A morphism of algebraic varieties $\psi: Y \rightarrow Z$ is the same thing as a family of
  compatible morphisms $\psi_i: Y_i \rightarrow Z$ (\ie\ such that $\psi_j \circ f_{j,i} = \psi_i$
  on $U_{i,j}$).
\end{cor}
\begin{cor}
  Let $X$ an algebraic variety, $(X_i)_{i \in I}$ an open covering of $X$ and let
  $\phi_i: X_i \rightarrow Y_i$ a family of morphisms such that $\phi_i(X_i \cap X_j) \subset U_{i,j}$
  for all $i,j$, and moreover supposed compatible (\ie\ such that   $\phi_j \circ f_{j,i} = \phi_i$
  on $X_i \cap X_j$). Then they can be uniquely patched (in an obvious sense) into a morphism
  $\phi: X \rightarrow Y$.
\end{cor}
  
\paragraph{Glueing of affine good and geometric quotients.} 
    
Let $G$ a reductive group, $X$ a $G$-variety and $(X_i)_{i \in I}$ an affine open covering of $X$.
We assume that the $X_i$ are $G$-invariant and that\footnote{Note that if we assume $X$ to be
separated (which will be the case in our applications), it is automatically true that the
intersection of two affine open subsets is affine.} $X_i \cap X_j$ is affine for all $i,j$. Then,
by \ref{subsubsection:quotients}, there are affine good quotients $\phi_i: X_i \rightarrow Y_i$
and their universal properties provide glueing data on the $U_{i,j} := \phi_i(X_i \cap X_j)$. From
the above, we obtain a unique morphism $\phi: X \rightarrow Y$ such that each restriction
$\phi_{| X_i}$ is identified with the affine good quotient $\phi_i: X_i \rightarrow Y_i \subset Y$;
and the $Y_i$ are an affine open covering of $Y$. \\
To obtain geometric quotients, we must extend the definition of stable points to non affine varieties:
\begin{defn}
  \label{defn:stablepoints2}
  The point $x \in X$ is \emph{stable} if belongs to a $G$-invariant affine subset of $X$, 
  its orbit $Gx$ is closed and its stabilizer $G_x$ is finite.
\end{defn}
Now the following fact is stated in \cite{NewsteadTata} for separated varieties, but it obviously
extends to the case the target $Y$ is not separated, whether we use as definition 
\ref{defn:geometricquotient} or \ref{prop:geometricquotient}; here we assume $X$ to be separated
(but $Y$ is arbitrary):
\begin{prop}
  Being a good or geometric quotient is a local property in the following sense: \\
  (i) If $\phi: X \rightarrow Y$ is a good, resp. a geometric quotient, then so is
  $\phi^{-1}(U) \rightarrow U$ for every open $U \subset Y$. \\
  (ii) If there is an open covering $(U_i)$ of $Y$ such that all the $\phi^{-1}(U_i) \rightarrow U_i$
  are good, resp. geometric quotients, then so is $\phi: X \rightarrow Y$. It follows that every
  map $\phi^{-1}(U) \rightarrow U$ is a categorical quotient (and an orbit space if all orbits
  in $\phi^{-1}(U)$ are closed).
\end{prop}
Combining with theorem \ref{thm:stablepoints1}, we get:
\begin{thm}
  \label{thm:stablepoints2}
  (i) The set $X^s$ of stable points is a (possibly empty) $G$-invariant open subset of $X$ and
  $\phi(X^s) \subset X//G$ is open too, so (provided $X^s$ is non empty) $\phi(X^s) = X^s//G$. \\
  (ii) The restriction $X^s \rightarrow \phi(X^s) = X^s//G$ is actually a geometric quotient $X^s/G$
  (again provided $X^s$ is non empty).
\end{thm}

% 1.2.4

\subsubsection{A useful criterion of isomorphy}
\label{subsubsection:IsomorphyCriterion}

Since we intend to consider non separated quotients, we shall have to use the following classical
criterion in that extended case, so we state and prove it.

\begin{prop}
  \label{prop:IsomorphyCriterion}
  Let $\varphi: X \rightarrow Y$ be a bijective morphism of algebraic varieties. We suppose that
  $X$ is separated and that $Y$ is normal, non necessarily separated. Then $\phi$ is an
  \emph{isomorphism}.
\end{prop}
\Proof
The problem is local on $Y$, therefore we can suppose that $Y$ is separated. Then, we can apply
\cite[Proposition (3.17),page $46$]{MumfordCPV}. It implies that $\varphi$ is \emph{birational}.
It follows from Zariski's Main Theorem \cite[(3.20),page $46$]{MumfordCPV} that a finite birational
morphism from a variety to a normal variety is an isomorphism, therefore $\varphi$ is an isomorphism.
\Finpr

% 1.2.5

\subsubsection{An example of non separated geometric quotient}
\label{subsubsection:NSGQ}

The following example will be needed later: it involves the diagonal action of $\Cs$ on
$\left(\Pr^1(\C)\right)^n$. \\

First, some basic notations and a baby version ($n = 1$) of the example. We identify $\Pr^1(\C)$
with $\hat{\C} := \C \cup \{\infty\}$ by $[a:b] \mapsto a/b$. We write $\rho$, resp. $\tilde{\rho}$
for the ``coordinate'' $a/b$, resp. for its inverse $1/\rho = b/a$. There is a natural action of $\Cs$
on $\Pr^1(\C)$ given (for $\lambda \in \Cs$) by $[a:b] \mapsto [\lambda a:b]$, \ie\ by
$\rho \mapsto \lambda \rho$ and $\tilde{\rho} \mapsto \lambda^{-1} \tilde{\rho}$. That action has
two fixed points $0$ and $\infty$. Let $\Theta_1 := \{0,\infty\}$ the set of fixed points. All
points of $\Pr^1(\C) \setminus \Theta_1 = \Cs$ have trivial stabilizer and a non closed orbit in
$\Pr^1(\C)$ but a closed orbit in the invariant open subset $\Pr^1(\C) \setminus \Theta_1$. The action
on that subset has a geometric quotient:
$$
\left(\Pr^1(\C) \setminus \Theta_1\right)/\Cs = \Cs/\Cs = \Pr^0(\C) = \{\bullet\}.
$$
Now we look at $\left(\Pr^1(\C)\right)^n$ for an arbitrary $n \geq 2$. The $\rho$ and $\tilde{\rho}$
coordinates of the components give rise to maps $\rho_1,\ldots,\rho_n$ and
$\tilde{\rho}_1,\ldots,\tilde{\rho}_n$ from $\left(\Pr^1(\C)\right)^n$ to $\hat{\C}$. The diagonal
action of $\Cs$ on $\left(\Pr^1(\C)\right)^n$ is characterized (for $\lambda \in \Cs$) by 
$\rho_i \mapsto \lambda \rho_i$ and $\tilde{\rho}_i \mapsto \lambda^{-1} \tilde{\rho}_i$ for
$i = 1,\ldots,n$. The set of fixed points is $\Theta_n := \{0,\infty\}^n$. The only closed orbits
are those of fixed points $p \in \Theta_n$: they are singletons $\Cs p = \{p\}$. Every other point
$p \in \left(\Pr^1(\C)\right)^n \setminus \Theta_n$ has a trivial stabilizer. Its orbit has boundary
$\overline{\Cs p} \setminus \Cs p \subset \Theta_n$, so it is closed in the invariant open subset
$\left(\Pr^1(\C)\right)^n \setminus \Theta_n$. \\

In the affine open subset $\C^n \subset \left(\Pr^1(\C)\right)^n$, the only fixed point is
$\zero:= (0,\ldots,0)$. All other points have trivial stabilizer and an orbit which is closed
in the invariant open subset $\C^n \setminus \{\zero\}$. The latter has a geometric quotient
under the diagonal $\Cs$-action:
$$
\C^n \setminus \{\zero\} \rightarrow \left(\C^n \setminus \{\zero\}\right)/\Cs = \Pr^{n-1}(\C).
$$
More generally, let $p \in \Theta_n$. We shall use the following notations for the complementary
subsets of indices:
$$
I := \{i \in \{1,\ldots,n\} \tq \rho_i(p) = 0\} \text{~and~}
J := \{j \in \{1,\ldots,n\} \tq \rho_j(p) = \infty\} =
\{j \in \{1,\ldots,n\} \tq \tilde{\rho}_j(p) = 0\}.
$$
Then we introduce an invariant open subset of $\left(\Pr^1(\C)\right)^n$:
\begin{align*}
  U_p
  &:= \left\{p \in \left(\Pr^1(\C)\right)^n \tq
  \begin{cases}
    i \in I \Rightarrow \rho_i(p) \neq \infty, \\
    j \in J \Rightarrow \tilde{\rho}_j(p) \neq \infty,
  \end{cases} \right\} \\
    &:= \left\{p \in \left(\Pr^1(\C)\right)^n \tq
  \begin{cases}
    i \in I \Rightarrow \rho_i(p) \neq \infty, \\
    j \in J \Rightarrow \rho_j(p) \neq 0,
  \end{cases} \right\} 
\end{align*}
so that in particular $U_\zero = \C^n$. Note that the map $s_I$ which transforms $\rho_j$ into
$\tilde{\rho}_j$ for $j \in J$ (and leaves $\rho_i$ invariant for $i \in I$ is an automorphism
of $\left(\Pr^1(\C)\right)^n$ which sends $U_p$ to $\C^n$ and $p$ to $\zero$ (and conversely since
it is an involution). Therefore $U_p \cap \Theta_n = \{p\}$ and all points of $U_p \setminus \{p\}$
have trivial stabilizer and an orbit open in $U_p \setminus \{p\}$. So here again we have a geometric
quotient $\Uu_p$ of $U_p \setminus \{p\}$ and a commutative diagram:
$$
\xymatrix{
  U_p \setminus \{p\} \ar@<0ex>[r]^{s_I} \ar@<0ex>[d] & \C^n \setminus \{\zero\} \ar@<0ex>[d] \\
  \Uu_p \ar@<0ex>[r]^\simeq  & \Pr^{n-1}(\C)
}
$$
Now each $U_p \cap U_{p'}$, $p \neq p' \in \Theta_n$, is a Zariski-dense affine open subset (actually
a product of copies of $\C$ and of $\Cs$) in both $U_p \setminus \{p\}$ and $U_{p'} \setminus \{p'\}$
(indeed $p' \not\in U_p$ and conversely). It therefore projects to a Zariski-dense open subset
$\Uu_{p,p'}$ of $\Uu_p$ and also to a Zariski-dense open subset $\Uu_{p',p}$ of $\Uu_{p'}$, with a
canonical isomorphism $\Uu_{p,p'} \rightarrow \Uu_{p',p}$. We use those isomorphisms to glue the $\Uu_p$
along the $\Uu_{p,p'}$, yielding an irreducible algebraic variety $\mathcal{K}_n$ which is a geometric
quotient:
$$
\bigcup_{p \in \Theta_n} \left(U_p \setminus \{p\}\right) = \left(\Pr^1(\C)\right)^n \setminus \Theta_n
\rightarrow \mathcal{K}_n := \left(\left(\Pr^1(\C)\right)^n \setminus \Theta_n\right)/\Cs.
$$
Note that the $\Uu_p$ are compact (they are projective spaces) but they are Zariski-dense open subsets
of $\mathcal{K}_n$, which is therefore not separated.

% 1.3

\subsection{Some facts about toric varieties}
\label{subsection:toricvarieties}

Some of our most interesting applications will fall under this heading. Exceptionally, in the whole
of \ref{subsection:toricvarieties}, we denote $\Affi[X]$ the affine algebra of an affine variety $X$
so as to avoid confusion with the group algebra $\C[\Lambda]$ of a group $\Lambda$ (or the monoid
algebra $\C[M]$ of a monoid $M$). Conversely, an affine algebra being given, we denote $\Spec\ A$
the corresponding affine variety; it can be canonically realized with underlying set the set of
all algebra morphisms $A \rightarrow \C$.

% 1.3.1

\subsubsection{Reminders on tori}

General references for the structure of tori are \cite{Borel,Humphreys,SpringerLAG,SpringerGAIV}.

\paragraph{Tori and their character groups.}

Let $\Tn$ the linear algebraic group $\Cs^n$. Its affine algebra is
$\Affi[\Tn] = \C[X_1,X_1^{-1},\ldots,X_n,X_n^{-1}]$. In more intrinsic terms, a ($n$-dimensional)
torus $T$ is a linear algebraic group group isomorphic to $\Tn$. Let:
$$
\Xf(T) := \{\lambda \in \Affi[T] \tq \lambda: T \rightarrow \Cs \text{~is a group morphism}\}
$$
the group of \emph{characters} of $T$. Then $\Lambda := \Xf(T)$ is a finitely generated free abelian
group and the affine algebra of $T$ is the group algebra of $\Lambda$:
$$
\Affi[T] = \C[\Lambda] = \bigoplus_{\lambda \in \Lambda} \C \lambda.
$$
The group $\Lambda^\vee := \{\text{group morphisms~} \Lambda \rightarrow \Cs\}$ has a natural
structure of torus (if $\Lambda = \Z^n$, it is canonically identified with $\Cs^n$) and the map
$T \rightarrow \Lambda^\vee$, $g \mapsto (\lambda \mapsto \lambda(g))$ is an isomorphism of
algebraic groups. Therefore $T$ is completely determined (as a group and as a variety) by $\Xf(T)$.

\paragraph{Quotients of tori.}

Let $S \subset T$ a closed subgroup of the torus $T$. We keep the notation $\Lambda := \Xf(T)$.
Let $S$ a closed subgroup of $T$ and let $\Lambda'$ is the subgroup of $\Lambda$ made of those
characters which vanish on $S$: thus $\Lambda'$ is a finitely generated free abelian group.
Then $T' := T/S$ has a natural structure of torus with affine algebra $\Affi[T'] = \C[\Lambda']$. \\
Moreover, any generating set of the subgroup $\Lambda'$ is a family of equations defining $S$
(\ie\ a generating set of the ideal of $S$ in $T$); and $S$ is a torus if, and only if, the subgroup
$\Lambda'$ is saturated, \ie\ $\Lambda/\Lambda'$ has no torsion (whence is free abelian). In that
case, $\Affi[S] = \C[\Lambda/\Lambda']$. \\
Note that in that case we have two dual exact sequences:
$$
1 \rightarrow S \rightarrow T \rightarrow T' \rightarrow 1 \text{~and~}
0 \rightarrow \Lambda' \rightarrow \Lambda \rightarrow \Lambda/\Lambda' \rightarrow 0.
$$
\begin{exa}
  The subgroup $S := \{(a,b) \in \Cs^2 \tq a^2 = b^2\}$ of $T := \mathbf{T}_2$ is not a subtorus
  (it is no even connected), but the subgroup $S := \{(a,b,c,d) \in \Cs^4 \tq ab = cd\}$ of
  $T := \mathbf{T}_4$ is a subtorus isomorphic to $\mathbf{T}_3$.
\end{exa}

\paragraph{Linear representations of tori.}

Let $V$ is a finite dimensional space endowed with a rational\footnote{This means that the
corresponding map $T \rightarrow \GL(V)$ is a morphism of algebraic groups. We shall omit
the word ``rational'' in the sequel.} linear representation of $T$. Then:
$$
V = \bigoplus_{\lambda \in \Lambda} V_\lambda, \text{~where~}
\forall \lambda \in \Lambda \;,\; V_\lambda := \{v \in V \tq \forall g \in T \;,\; g.v = \lambda(g) v\}.
$$
If $V$ is infinite dimensional but a union of finite dimensional $T$-invariant subspaces such that
the corresponding linear representations of $T$ are rational, the same conclusion holds.

\paragraph{Actions of tori.}

If a torus $T = \Lambda^\vee = \Spec\ \C[\Lambda]$ acts on an affine algebraic variety $X$, it acts
dually on $\Affi[X]$ (for $g \in T$ and $f \in \Affi[X]$, we set $gf: x \mapsto f(g^{-1}x)$) and this
is an infinite dimensional linear representation, whence a decomposition
$\Affi[X] = \bigoplus\limits_{\lambda \in \Lambda} \Affi[X]_\lambda$. \\
In the particular case of the left multiplication action of $T$ on itself, the corresponding
decomposition is just $\Affi[T] = \C[\Lambda] = \bigoplus\limits_{\lambda \in \Lambda} \C \lambda$. \\
We quote Sumihiro's theorem from \cite[p. 10]{OdaConvex}: if the torus $T$ acts (rationally) on
an irreducible normal variety, then the later is the union of $T$-invariant affine open subsets.

% 1.3.2

\subsubsection{Toroidal embeddings and toric varieties}

General references here are \cite{Brion,Dolgachev} (for toric varieties) and
\cite{KKMSD,OdaTorus,OdaConvex} (for toroidal embeddings). We mostly follow those sources, except
that we distinguish toroidal embeddings from toric varieties; and we do not systematically assume
normality.

\paragraph{Toroidal embeddings.}

A \emph{toroidal embedding}, of the torus $T$ is a variety $X$ admitting $T$ as an open dense subset
(thus $X$ is irreducible) and an action of $T$ extending the left multiplication action of $T$ on
itself, whence a commutative diagram (vertical arrows are inclusions, horizontal arrows are actions):
$$
\xymatrix{
T \times T \ar@<0ex>[r] \ar@<0ex>[d] & T \ar@<0ex>[d] \\
T \times X \ar@<0ex>[r]  & X
}
$$
If $X$ is affine, the dominant inclusion $T \subset X$ induces an injection $\Affi[X] \subset \Affi[T]$.
The splitting of $\Affi[T]$ in eigenspaces $\Affi[T]_\lambda$ induces a splitting of $\Affi[X]$ (under
the action of $T$):
$$
\Affi[X] = \bigoplus_{\lambda \in M} \Affi[X]_\lambda = \bigoplus_{\lambda \in M} \C \lambda = \C[M],
$$
where $M$ is a semigroup, the submonoid $\Lambda \cap \Affi[X]$ of those characters of $T$ which
can be extended as regular maps on the whole of $X$. The semigroup $M$ is finitely generated (as
a semigroup) and it generates the group $\Lambda$. The variety $X$ is normal if, and only if,
$M$ is \emph{saturated} in $\Lambda$, \ie\ for $r \in \N^*$, $\lambda \in \Lambda$, one has the
implication $r \lambda \in M \Rightarrow \lambda \in M$. \\
The way $M$ embeds into $\Lambda$ admits a combinatorial description (``fans'') which encodes the
way $X'$ is completed into $X$ by the addition of orbits of smaller dimension. \\
A morphism of toroidal embeddings of $T$ is a $T$-equivariant morphism. Then there is a bijective
correspondence $M \mapsto \Spec\ \C[M]$ from finitely generated sub semigroups of $\Lambda$ to
isomorphism classes of affine toroidal embeddings of $T$; saturated semigroups correspond to normal
affine toroidal embeddings. More generally, morphisms of affine toroidal embeddings correspond (in
a contravariant way) to inclusions of semigroups (see \cite[pp 4,6]{KKMSD}).

\paragraph{Toric varieties.}

Here we follow\footnote{The groups $G$, $B$, $U$ of \emph{loc.cit.} are here $T$, $T$, $1$.}
\cite[\S 2.2]{Brion}. \\
A \emph{toric variety} is a normal irreducible $T$-variety $X$ (for some torus $T$) containing
an open orbit (which is then Zariski dense since $X$ is irreducible). \\
Let $X' = T x$, $x \in X$, such an orbit. Then the morphism $T \rightarrow X$, $g \mapsto g x$ is
dominant, so $\Affi[X]$ is (identified to) a subalgebra of $\Affi[T] = \C[\Lambda]$. More precisely,
there is a finitely generated monoid $M \subset \Lambda$ such that $\Affi[X] = \C[M]$. \\
The surjective morphism $T \rightarrow X'$, $g \mapsto g x$, realizes $X'$ as a homogeneous space
under $T$. By \cite[\S 1.1]{Brion}, it is isomorphic (as a variety) to the torus $T' := T/T_x$
(remember $T_x$ is the stabilizer of $x$). Let $\Lambda'$ the group of characters of $T'$, so that
$\Affi[T'] = \C[\Lambda']$. Interpreting characters $\lambda: T \rightarrow \Cs$ which vanish on $T_x$ 
as elements of $\Lambda'$, thus as regular functions on the orbit $X' = T x$, the elements of $M$
are those characters that can be extended to a regular function on the whole of $X \supset X'$. \\
Since $T_x$ acts trivially on the dense subset $X'$ of $X$, it acts trivially on the whole of $X$
and the action of $T$ on $X$ factors into an action of $T'$ on $X$. Identifying $T'$ with its
image $X'$, we see that the general situation of toric varieties boils down to the case of toroidal
embeddings.

% 1.3.3

\subsubsection{Important examples of toroidal embeddings}
\label{subsubsection:examplestor}

\paragraph{Two basic examples.}
\begin{enumerate}
\item Recall from \ref{subsubsection:NSGQ} the identification of $\Pr^1(\C)$ with
  $\hat{\C} := \C \cup \{\infty\}$ by $[a:b] \mapsto a/b$, whence an embedding of $\Cs$
  into $\Pr^1(\C)$, and then an embedding of $\Tn = \Cs^n$ into $\Pr^1(\C)^n$. The left
  action of $\Tn$ on itself extends to an action on $\Pr^1(\C)^n$ by:
  $$
  (\lambda_1,\ldots,\lambda_n).([a_1:b_1],\ldots,[a_n:b_n]) :=
  ([\lambda_1 a_1:b_1],\ldots,[\lambda_n a_n:b_n]),
  $$
  making it a toroidal embedding.
\item We shall meet later the quadric hypersurface $XY = ZT$ in $\Pr^3(\C)$. There is
  a natural action of the torus $\T_3 \simeq \{(a,b,c,d) \in \Cs^4 \tq ab = cd\}$ on
  that hypersurface, defined by:
  $$
  (a,b,c,d) [X:Y:Z:T] := [aX:bY:cZ:dT],
  $$
  again a toroidal embedding.
\end{enumerate}

\paragraph{A significant example.}
The following lemma is more or less obvious:
\begin{lem}
  \label{lem:quotienttoroidal}
  Let $T \hookrightarrow X$ a toroidal embedding and let $S \subset T$ a substorus
  such that $T/S$ is itself a torus. Assume that the quotient $X/S$ is geometric.
  Then $T/S \hookrightarrow X/S$ is a toroidal embedding.
\end{lem}
It follows that the variety $\mathcal{K}_n$ described in \ref{subsubsection:NSGQ} is a
non separated toric variety. Another example will be provided by corollary \ref{cor:V_S^{(*)}/H}.

% 1.4

\subsection{Invariants and quotients related to the Riemann-Hilbert-Birkhoff correspondence:
  general case ($n \geq 2$ arbitrary)}
\label{subsection:invariantsquotientsforgeneraln}

All the subsection \ref{subsection:invariantsquotientsforgeneraln}, with general $n \geq 2$ and
$\mu \geq 1$, is intended to provide vocabulary and basic tools for a future detailed study of
the spaces $\F_{R,S,\x}$. Only from section \ref{section:invariantsquotientsJScase} on where we
assume that $n = \mu = 2$, do we present substantial results; so this subsection may be skipped
at first reading. \\

Let $V_{i,j}$, $1 \leq i,j \leq n$, be non trivial finite dimensional complex vector spaces and,
according to our general conventions:
$$
V := \prod_{1 \leq i,j \leq n} V_{i,j} \text{~and~} V^{(*)} := \prod_{1 \leq i,j \leq n} V_{i,j}^*.
$$
We write their elements as $M = (m_{i,j})$ (since they are here as models of the context described
in \ref{subsubsection:thecontext}). Let the $2n$ dimensional torus $\Dnc \times \Dnc$ act on $V$ by:
$$
(\Gamma,\,\Delta).M := \Gamma M \Delta^{-1} =
\left(\frac{\gamma_i}{\delta_j} m_{i,j}\right)_{i,j = 1,\ldots,n},
\text{~where~}
\Gamma =: \Diag(\gamma_1,\ldots,\gamma_n) \text{~and~} \Delta =: \Diag(\delta_1,\ldots,\delta_n).
$$
The kernel of the action is clearly $\Cs (I_n,I_n) = \{(\gamma I_n,\gamma I_n) \tq \gamma \in \Cs\}$
so we set:
$$
H := \dfrac{\Dnc \times \Dnc}{\Cs (I_n,I_n)},
$$
which acts faithfully on $V$. The group $H$ is itself a $(2n-1)$-dimensional torus: identifying
the obvious way $\Dnc \times \Dnc$ to $\Cs^{2n}$, it is the direct product of its algebraic subgroup
$\Cs (I_n,I_n)$ by anyone of the $(2n-1)$-dimensional subtori $\Cs^i \times \{1\} \times \Cs^j$,
$i + j = 2n-1$, so any of the corresponding mappings $\Cs^i \times \{1\} \times \Cs^j \rightarrow H$
is an isomorphism.

% 1.4.1

\subsubsection{Stabilizers}
\label{subsubsection:stabilizers}

(The class of) an element $(\underline{\lambda},\underline{\mu}) \in \Dnc \times \Dnc$ is in the
stabilizer $H_M$ of $M := (m_{i,j}) \in V$ \iff\ $m_{i,j} \neq 0 \Rightarrow \lambda_i = \mu_j$. Those
equations define a subtorus of $\Dnc \times \Dnc$, with dimension the number of degrees of freedom
left by those equations among the $2 n$ coefficients $\lambda_i$, $\mu_j$, and $H_M$ is a torus of
dimension that number minus $1$.

\begin{exas}
(i) Let $n = 2$. If $M$ has at most one zero component, $H_M$ is trivial. Otherwise it is not. \\
(ii) Let $n = 3$. If $M$ has zeroes at most at positions $(1,1)$, $(1,2)$, $(2,1)$, $(2,2)$, then
$H_M$ is trivial. If $M$ has a null column or a null line, $H_M$ is not trivial.
\end{exas}

So we define the support of $M$ as $S(M) := \{(i,j) \in \{1,\ldots,n\}^2 \tq m_{i,j} \neq 0\}$. Then,
for any subset $S \subset \{1,\ldots,n\}^2$, consider a bipartite graph with set of vertices
$\{1,\ldots,n\} \sqcup \{1,\ldots,n\}$ (disjoint union), where $i$ on the left is connected to $j$
on the right \iff\ $(i,j) \in S$; and let $\chi(S)$ the number of connected components of that graph.

\begin{prop}
  The dimension of the torus $H_M$ is $\chi(S(M)) - 1$.
\end{prop}
\Proof
The class of $(\Gamma,\Delta) \in \Dnc \times \Dnc$ is in $H_M$ \iff\ $\gamma_i = \delta_j$ for all
$(i,j) \in S(M)$. This leaves one degree of freedom (in choosing all $\gamma_i,\delta_j$) for each
connected component, so the set of all such pairs (\ie\ the  preimage of $H_M$ in $\Dnc \times \Dnc$)
has dimension $\chi(S(M))$. 
\Finpr

\begin{exas}
(i) Let $n = 2$. Then $\chi(S) = \max(1, 4 - \vert S \vert)$. \\
(ii) Let $n = 3$. Then $\chi(S) > 1$ if $S \subset \{i_1,i_2\} \times \{1,2,3\}$ for some $i_1,i_2$ or
if $S \subset \{1,2,3\} \times \{j_1,j_2\}$ for some $j_1,j_2$, while $\chi(S) = 1$ if the complement
of $S$ is contained in $\{i_1,i_2\} \times \{j_1,j_2\}$ for some $i_1,i_2,j_1,j_2$.
\end{exas}

% 1.4.2

\subsubsection{Orbits}

Let $M \in V$ and $S := S(M)$. We identify $V_S := \prod\limits_{(i,j) \in S} V_{i,j}$ with the subspace
of $V$ defined by a zero component at all $(i,j) \not\in S$. We denote (as in our general conventions)
$V_S^{(*)} := \prod\limits_{(i,j) \in S} V_{i,j}^*$, which we identify to the subset of $V_S$ (thus $V_S^{(*)}$
is open in $V_S$ which is closed in $V$). Clearly the orbit $H.M$ of $M$ is contained in $V_S^{(*)}$.
More precisely, let us set:
$$
L_{i,j}^*(M) := \Cs m_{i,j} \subset V_{i,j}^* \text{~and~}
L_S^{(*)}(M) := \prod\limits_{(i,j) \in S} L_{i,j}^* \subset V_S^{(*)}.
$$
Then:
$$
H.M \subset \prod_{(i,j) \in S} L_{i,j}^*(M) = L_S^{(*)}(M) \subset \prod_{(i,j) \in S} V_{i,j}^* = V_S^{(*)}.
$$
The factors $\lambda_i/\mu_j$ effectively involved in the action on $H.M$ are those such that
$(i,j) \in S$. Under the obvious isomorphism $L_S^{(*)}(M) \simeq \Cs^S$, the orbit $H.M$ goes to:
$$
T_S := \{(\lambda_i/\mu_j)_{(i,j) \in S} \tq \text{~all~} \lambda_i,\mu_j \in \Cs\} \subset \Cs^S.
$$

\begin{prop}
  $T_S$ is a closed subgroup of $\Cs^S$ and a torus of dimension $2 n - \chi(S)$.
\end{prop}
\Proof
  It is the image of the map:
  $$
  \Dnc \times \Dnc \rightarrow \Cs^S, \;
  (\underline{\lambda},\underline{\mu}) \mapsto (\lambda_i/\mu_j)_{(i,j) \in S}.
  $$
  That map being a morphism of algebraic groups, the image is closed (Borel, chap 1, \S 1, cor. 1.4).
  Since the kernel has dimension $\chi(S)$, the image has dimension $2 n - \chi(S)$. Also, being a
  connected closed subgroup of the torus $\Cs^S$, it is a torus.
\Finpr

We make this more precise as follows. Let $\pi_{i,j}: V_{i,j}^* \rightarrow \Pr(V_{i,j})$ the natural
projection and (with a slight abuse of the notation $\Pr$):
$$
\pi := \prod_{(i,j) \in S} \pi_{i,j}: V_S^{(*)} \rightarrow \Pr(V_S^{(*)}) := \prod_{(i,j) \in S} \Pr(V_{i,j}),
$$
so that $L_S^{(*)}(M) = \pi^{-1}(\pi(M))$ is closed in $V_S^{(*)}$ and $H.M$ is closed in $L_S^{(*)}(M)$.
We have commutative diagrams (in which horizontal arrows are canonical inclusions, vertical arrows are
non canonical isomorphisms):
$$
\xymatrix{
H.M \ar@<0ex>[r]  \ar@<0ex>[d]  & L_S^{(*)}(M) \ar@<0ex>[d]  \\
T_S \ar@<0ex>[r]  & \Cs^S
}
$$

\begin{cor}
  The orbit $H.M$ is closed in $V_S^{(*)}$.
\end{cor}

% 1.4.3

\subsubsection{Affine covering of $V_S^{(*)}$}

The problem is that $V_S^{(*)}$ is, in general, not an affine variety: if $W$ is a vector space, $W^*$
is an affine variety \iff\ $\dim W = 1$. However, for every linear form $e \in W^\vee$ (the dual),
$e \neq 0$, the open subset $W(e) := W \setminus \Ker\ e$ is an affine variety, with affine algebra
$\Af[W(e)] = \Af[W][1/e]$. Moreover, each $W(e)$ is stable under the natural $\Cs$-action and the
$W(e)$, $e \in W^\vee \setminus \{0\}$, cover $W^*$: the corresponding quotients $W(e)/\Cs$ are
the natural affine charts $\Pr(W)_e \subset \Pr(W)$. (Actually, taking the elements $e$ in a basis
of $W^\vee$ is enough.) \\

So, writing for short $V_S^{\vee *} := \prod\limits_{(i,j) \in S} (V_{i,j}^\vee \setminus \{0\})$, we define:
$$
\forall \e := (e_{i,j})_{(i,j) \in S} \in V_S^{\vee *} \;,\;
V_S(\e) := \prod\limits_{(i,j) \in S} V_{i,j}(e_{i,j}) \subset V_S^{(*)}.
$$
Those open subsets cover $V_S^{(*)}$, they are stable under the action of $H$ and they are affine.
Now, the really useful result is:

\begin{cor}
  \label{cor:V_S^{(*)}/H}
  (i) If $M \in V_S(\e)$, the orbit $H.M$ is closed in $V_S(\e)$. \\
  (ii) If moreover $\chi(S) = 1$, then all $M \in V_S(\e)$ are stable (in the sense of Brion,
  definition 1.25). As a consequence, there is a geometric quotient of $V_S^{(*)}$ by $H$. \\
  (iii) $V_S^{(*)}/H$ is a toric variety.
\end{cor}

The last statement flows from lemma \ref{lem:quotienttoroidal}.

\begin{exa}
This applies in particular to $n = 2$ and $S = \{1,2\} \times \{1,2\}$ or the same minus one ordered
pair $(i,j)$, \ie\ to the formats (each $\star$ stands for a non zero coefficient):
$$
\begin{pmatrix} \star & \star \\ \star & \star \end{pmatrix}, \quad
\begin{pmatrix} 0 & \star \\ \star & \star \end{pmatrix}, \quad
\begin{pmatrix} \star & 0 \\ \star & \star \end{pmatrix}, \quad
\begin{pmatrix} \star & \star \\ 0 & \star \end{pmatrix}, \quad
\begin{pmatrix} \star & \star \\ \star & 0 \end{pmatrix}.
$$
\end{exa}

We write down explicitly some affine algebras; for that, we endow each dual $V_{i,j}^\vee$ with a
particular basis $(u_{i,j}^{(k)})_{1 \leq k \leq d_{i,j}}$, where $d_{i,j} := \dim_\C V_{i,j}$. We let
$\Pr(V_{i,j})_{e_{i,j}}$ the image of $V_{i,j}(e_{i,j})$ in $\Pr(V_{i,j})$ and $\Pr(V_S^{(*)})_\e$ their product
for $(i,j) \in S$, \ie\ the image of $V_S(\e)$ in $\Pr(V_S^{(*)})$. Then:
$$
  \Af[V_{i,j}] = \C[(u_{i,j}^{(k)})_{1 \leq k \leq d_{i,j}}], \qquad
  \Af[V_{i,j}(e_{i,j})] = \C[(u_{i,j}^{(k)})_{1 \leq k \leq d_{i,j}}][1/e_{i,j}],
$$
$$
  \Af[V_S] = \bigotimes_{(i,j) \in S} \C[(u_{i,j}^{(k)})_{1 \leq k \leq d_{i,j}}] = \C[\text{all~} u_{i,j}^{(k)}],
$$
$$
  \Af[V_S(\e)] = \bigotimes_{(i,j) \in S} \C[(u_{i,j}^{(k)})_{1 \leq k \leq d_{i,j}}][1/e_{i,j}] =
  \Af[\text{all~} u_{i,j}^{(k)}][\text{all~} 1/e_{i,j}],
 $$
$$
  \Af[\Pr(V_{i,j})_{e_{i,j}}] = \C[(u_{i,j}^{(k)}/e_{i,j})_{1 \leq k \leq d_{i,j}}], \qquad
  \Af[\Pr(V_S^{(*)})_\e] = \C[\text{all~} u_{i,j}^{(k)}/e_{i,j}].
$$

If $\e$ is fixed, we shall write $x_{i,j}^{(k)} := u_{i,j}^{(k)}/e_{i,j}$, so that:
$$
\Af[\Pr(V_{i,j})_{e_{i,j}}] = \C[(x_{i,j}^{(k)})_{1 \leq k \leq d_{i,j}}] \text{~~and~~}
\Af[\Pr(V_S^{(*)})_\e]  = \C[\text{all~} x_{i,j}^{(k)}].
$$
Since each family $(u_{i,j}^{(k)})_{1 \leq k \leq d_{i,j}}$ is a basis of $V_{i,j}^\vee$, there are linear
relations:
$$
e_{i,j} = \sum_{1 \leq k \leq d_{i,j}} \alpha_{i,j}^{(k)} u_{i,j}^{(k)} \Longrightarrow
\sum_{1 \leq k \leq d_{i,j}} \alpha_{i,j}^{(k)} x_{i,j}^{(k)} = 1.
$$
Note that, although each $u_{i,j}^{(k)}$ is a function on $V_{i,j}$ and each $x_{i,j}^{(k)}$ is a function
on $V_{i,j}(e_{i,j})$, we respectively consider them (in an obvious way) as functions on $V_S^{(*)}$, resp.
on $V_S(\e)$.

\begin{exa}
  In the JS case (for ``Jimbo-Sakai''), \ie\ for $n = 2$, $S = \{1,2\} \times \{1,2\}$ and all
  $\dim V_{i,j} = 2$, we shall rather denote $(u_{i,j},v_{i,j})$ the chosen basis of $V_{i,j}^\vee$;
  and $x_{i,j} := u_{i,j}/e_{i,j}$, $y_{i,j} := v_{i,j}/e_{i,j}$. Then, a fixed $\e$ will be expressed as
  $e_{i,j} = \alpha_{i,j} u_{i,j} + \beta_{i,j} v_{i,j}$, so $\alpha_{i,j} x_{i,j} + \beta_{i,j} y_{i,j} = 1$.
\end{exa}

\begin{exa}
  In degenerate JS case, which is the same except that $S \subset \{1,2\} \times \{1,2\}$ misses one
  particular pair of indices (so that one coefficient is $0$, the three other ones are $\neq 0$), we
  keep those notations for all $(i,j) \in S$.
\end{exa}

% 1.4.4

\subsubsection{Equations}
\label{subsubsection:equations}

Since $T_S$ is a closed subgroup of the torus $\Cs^S$, it is defined by a set of monomial equations,
actually a subgroup of the character group $\mathbf{X}(\Cs^S) = \Z^S$. This subgroup is free, so it
admits a basis and there is a basis of monomial equations for $T_S$ of the form
$M(X_{i,j})_{(i,j) \in S} = 1$.
We shall determine in a moment such a basis. \\

At any rate, letting:
$$
\Xf(\Cs^S) = \left\{\prod_{(i,j) \in S} X_{i,j}^{r_{i,j}} \tq \text{~all~} r_{i,j} \in Z\right\} \simeq \Z^S,
$$
the group of characters of $\Cs^S$, let $\M_S \subset \Xf(\Cs^S)$ such a set of monomials, \ie\ a basis
for $\Xf(\Cs^S/T_S)$. Then, for each $M(X_{i,j}) \in \M_S$ and for each
$\e := (e_{i,j})_{(i,j) \in S} \in V_S^{\vee *})$, the substitutions $\e^* M := M(e_{i,j})$ define regular
functions on $V_S(\e)$ constant on each orbit.

\begin{rmk}
  One can prove, \eg\ by using the precise form of the basic cycles, that $\Cs^S/T_S$ is itself
  a torus.
\end{rmk}

\begin{lem}
  Relations $\e^* M = \text{constant}$, $M \in \M_S$, make up a complete set of equations for each orbit
  $H.M$ in $L^*(M)$.
\end{lem}
\Proof
  Indeed, $M = (m_{i,j})_{(i,j) \in S}$ being fixed, we may choose the isomorphism $L^*(M) \simeq \Cs^S$
  sending each $N = (n_{i,j})_{(i,j) \in S}$ to $(e_{i,j}(n_{i,j})/e_{i,j}(m_{i,j}))_{(i,j) \in S}$. Then the
  pullbacks of equations $M_c = 1$ are equations $\e^* M_c(N) = \e^* M_c(M)$.
\Finpr

Since the dimension of $T_S$ is $2n - \chi(S)$, the basis $\M_S$ has
$c(S) := \vert S \vert - 2 n + \chi(S)$ éléments. This number is the \emph{circuit rank} or the
\emph{cyclomatic number} of our bipartite graph (see either \cite[chap. 4]{Berge}, or Wikipedia,
heading ``circuit rank''). Here is a way to generate equations from cycles. We write $C(S)$ a fixed
basis of elementary cycles in the above bipartite graph. Every $c \in C(S)$ has the form of a loop
$i_1 \to j_1 \to i_2 \to j_2 \to \cdots \to i_k \to j_k \to i_1$, where $i_1,\ldots,i_k$, resp.
$j_1,\ldots,j_k$ are pairwise distinct, and where the cycle $c$ does not depend on the particular
selected origin $i_1$ of the loop. Then the coordinates $x_{i,j}$ of an element of $T_S$ satisfy:
$$
x_{i_1,j_1} x_{i_2,j_2} \cdots x_{i_k,j_k} =
\dfrac{\lambda_{i_1} \cdots \lambda_{i_k}}{\mu_{j_1} \cdots \mu_{i_k}} =
x_{i_1,j_k} x_{i_2,j_1} \cdots x_{i_k,j_{k-1}}
$$
for any antecedent $(\underline{\lambda},\underline{\mu})$. We then set:
$$
M_c := \dfrac{X_{i_1,j_1} X_{i_2,j_2} \cdots X_{i_k,j_k}}{X_{i_1,j_k} X_{i_2,j_1} \cdots X_{i_k,j_{k-1}}} \in \M_S.
$$
This is almost well defined from $c$, independent of the point of departure (chosen among the $i$-indices,
\ie\ in the left factor of $\{1,\ldots,n\} \times \{1,\ldots,n\} \supset S$). The only freedom left comes
from the change of orientation on the cycle, which can transform $M_c$ to $M_c^{-1}$. We can freeze it by
deciding for instance that $i_1$ is the smallest possible in the cycle, then $j_1$ is the smallest
possible for this $i_1$. \\

That being set, almost all statements until the end of \ref{subsubsection:equations}, as well as in
\ref{subsubsection:algebrasofinvariants}, flow from elementary combinatorial arguments and we leave
their proof to the reader.

\begin{prop}
  $T_S$ is exactly the closed subset of $\Cs^S$ defined by these equations, which are independent.
\end{prop}

\begin{cor}
  $H.M$ is exactly the closed subset of $L^*(M)$ defined by equations $\e^* M_c = \text{constant}$,
  $c \in C(S)$, which are independent.
\end{cor}

\begin{exa}
  In the JS case, there is only one equation for $T_S$, namely $x_{1,1} x_{2,2} = x_{1,2} x_{2,1}$, so the
  equations of orbits are $x_{1,1} x_{2,2}/x_{1,2} x_{2,1} = \text{constant}$.
\end{exa}

\begin{exa}
  In the degenerate JS case, there are no cycles and no equations: $H.M = L^*(M)$.
\end{exa}

We have the equations of $H.M$ in $L^*(M)$, now we look for the equations of $L^*(M)$ in $V_S(\e)$.
Since $L^*(M) = \pi^{-1}(\pi(M))$, using that:
$$
\pi(N) = \pi(M) \Longleftrightarrow
\forall (i,j) \in S \;,\; \forall u \in V_{i,j}^\vee \;,\; u(M) = u(N),
$$
so, resorting to the selected bases:
$$
\pi(N) = \pi(M) \Longleftrightarrow
\forall (i,j) \in S \;,\; \forall k = 1,\ldots, d_{i,j} \;,\; u_{i,j}^{(k)}(M) = u_{i,j}^{(k)}(N).
$$

\begin{lem}
  A complete set of equations of $L^*(M)$ in $V_S(\e)$ is given by the $u_{i,j}^{(k)} = \text{constant}$.
\end{lem}

\begin{cor}
  $H.M$ is exactly the closed subset of $V_S(\e)$ defined by equations $\e^* M_c = \text{constant}$,
  $c \in C(S)$, and $u_{i,j}^{(k)} = \text{constant}$, $(i,j) \in S$, $k = 1,\ldots, d_{i,j}$, which are
  independent.
\end{cor}

The above statements will be made more precise by computing algebras of invariants.

% 1.4.5

\subsubsection{Algebras of invariants}
\label{subsubsection:algebrasofinvariants}

We know that $\Af\left[\Cs^S\right] = \C\left[(X_{i,j},X_{i,j}^{-1})_{(i,j) \in S}\right]$. From now on,
we take the set $\M_S := \{M_c \tq c \in C(S)\}$ as a basis of equations of $T_S$ in $\Cs^S$.

\begin{lem}
The affine algebra of $T_S$ is:
$$
\Af[T_S] = \C\left[\{M,M^{-1} \tq M \in \M_S\}\right] = \C\left[\{M_c,M_c^{-1} \tq c \in \C(S)\}\right].
$$
\end{lem}

\begin{prop}
  The algebra of invariants on $V_S(\e)$ is:
  $$
  \Af[V_S(\e)/H] = \Af[V_S(\e)]^H = \C[\text{all~} x_{i,j}^{(k)}, \text{~all~} \e^*M_c,\e^*M_c^{-1}],
  $$
  with $\vert S \vert$ relations $\sum\limits_k \alpha_{i,j}^{(k)} x_{i,j}^{(k)} = 1$.
\end{prop}

\begin{cor}
  The quotient $V_S(\e)/H$ is isomorphic to $\Pr(V_S^{(*)})_\e \times T_S$.
\end{cor}

\begin{thm}
  The quotient $V_S/H$ is a fibered space over $\Pr(V_S^{(*)})$ with fibers isomorphic to $T_S$.
\end{thm}
\Proof
  The affine spaces $\Pr(V_S^{(*)})_\e$, $\e \in V_S^{\vee *}$, patch up into $\Pr(V_S^{(*)})$.
  For $\e,\e' \in V_S^{\vee *}$, the functions $\phi_{i,j} := e'_{i,j}/e_{i,j}$ from
  $\Pr(V_S^{(*)})_\e \cap \Pr(V_S^{(*)})_{\e'}$ to $\C^*$ define the transition functions on the affine
  atlas. They yield the transition functions on the fibers defined by the $M_c(\phi_{i,j})$.
  \Finpr

\begin{exa}
  In the JS case (recall that this means $n = \mu = 2$), we get a $\C^*$-fibration on $\Pr^1(\C)^4$.
  The transition functions are the
  $\dfrac{(e'_{1,1}/e_{1,1})(e'_{2,2}/e_{2,2})}{(e'_{1,2}/e_{1,2})(e'_{2,1}/e_{2,1})} \cdot$ \\
  In the degenerate JS case, we get $\Pr^1(\C)^3$.
\end{exa}

%%%%%%%%%%%%%%%%%%%%%%%%%%%%%%%%%%%%%%%%%%%%%%%%%%%%%%%%%%%%%%%%%%%%%%%%%%%%%

% 2

\section{Invariants and quotients in the (possibly degenerate) JS case}
\label{section:invariantsquotientsJScase}

% 2.1

\subsection{Some general notations}

Here $n = 2$ and each $V_{i,j}$ has dimension $2$. The group $H$ is
$\dfrac{\Ddc \times \Ddc}{\Cs (I_2,I_2)} \cdot$ \\

For $i,j \in \{1,2\}$, we fix a basis $(u_{i,j},v_{i,j})$ of $V_{i,j}^\vee$. \\

For $e_{i,j} = \alpha_{i,j} u_{i,j} + \beta_{i,j} v_{i,j} \in V_{i,j}^\vee \setminus \{0\}$, we set
$V_{i,j}(e_{i,j}) := V_{i,j} \setminus \Ker\ e_{i,j}$. Those form an affine covering of $V_{i,j}$.
The corresponding affine algebras are:
$$
\Af\left(V_{i,j}(e_{i,j})\right) = \C[u_{i,j},v_{i,j}][e_{i,j}^{\pm 1}].
$$
The associated projective space $\Pr(V_{i,j})$ is covered by the affine charts $\Pr(V_{i,j})_{e_{i,j}}$
images of the $\Af(V_{i,j}(e_{i,j}))$; the corresponding affine algebras are:
$$
\Af\left(\Pr(V_{i,j})_{e_{i,j}}\right) = \C[u_{i,j},v_{i,j}][e_{i,j}^{\pm 1}] = \C[x_{i,j},y_{i,j}][e_{i,j}^{\pm 1}],
$$
where we have set $x_{i,j} := u_{i,j}/e_{i,j}$ and $y_{i,j} := v_{i,j}/e_{i,j}$ (the letters $x,y$ are fixed
although the variables $x_{i,j},y_{i,j}$ actually depend on the linear form $e_{i,j}$ and on the particular
affine chart). They are not independent variables:
$$
\alpha_{i,j} x_{i,j} + \beta_{i,j} y_{i,j} = \alpha_{i,j} u_{i,j}/e_{i,j} + \beta_{i,j} v_{i,j}/e_{i,j} = 1.
$$

\begin{rems}
  (i) Let $V := V_{1,1} \times V_{1,2} \times V_{2,1} \times V_{2,2}$. Then the only invariants in $\Af(V)$
  under the $H$-action are the constants (this will easily follow from the calculations to come). Thus
  the categorical quotient $V/H$ of the affine variety $V$ by the reductive group $H$ is trivial. \\
  (ii) On the other hand, $e_{1,1},e_{1,2},e_{2,1},e_{2,2}$ being chosen as above, all the corresponding
  $x_{i,j},y_{i,j}$ (considered as rational functions on $V$) are invariant under the $H$-action and so
  is $\phi := (e_{1,1} e_{1,2})/(e_{2,1},e_{2,2})$. It will also follow from the calculations to come that
  they together generate the $\C(V)^H$ (where $\C(V)$ is the fraction field of $\Af(V)$); so $\C(V)^H$
  has transcendence degree $5$.
\end{rems}

% 2.2

\subsection{Non degenerate JS case: no zeroes allowed}

Here $S = \{1,2\} \times \{1,2\}$. We have:
$$
V_S = V^{(*)} = V_{1,1}^* \times V_{1,2}^* \times V_{2,1}^* \times V_{2,2}^*.
$$
Let $\e := (e_{1,1},e_{1,2},e_{2,1},e_{2,2})$, each $e_{i,j} \in V_{i,j}^\vee \setminus \{0\}$. Then:
$$
V_S(\e) = V^{(*)}(\e) = V_{1,1}(e_{1,1}) \times V_{1,2}(e_{1,2}) \times V_{2,1}(e_{2,1}) \times V_{2,2}(e_{2,2}).
$$
Those open subsets of $V^{(*)}$ form an affine covering. The corresponding affine algebras are:
$$
\Af\left(V^{(*)}(\e)\right) = \bigotimes_{i,j} \C[x_{i,j},y_{i,j}][e_{i,j}^{\pm 1}] =
\C[\text{all~} x_{i,j},y_{i,j}][\text{all~} e_{i,j}^{\pm 1}].
$$
To study the invariants under $H$, we set the torus:
$$
T_S = \Cs^4 \text{~with affine algebra~} \Af(T_S) = \C[\text{all~} X_{i,j}^{\pm 1}, i,j = 1,2] =
\bigoplus_{M \in \M^*} \C M,
$$
where $\M^* = \M_S$ is the free abelian group over the $X_{i,j}$, \ie\ the set of all monomials
$M := X_{1,1}^{r_{1,1}} X_{1,2}^{r_{1,2}} X_{2,1}^{r_{2,1}} X_{2,2}^{r_{2,2}}$, all $r_{i,j} \in \Z$. In particular
we distinguish $\Phi := \dfrac{X_{1,1} X_{2,2}}{X_{1,2} X_{2,1}} \cdot$ For $M \in \M^*$ and for $\e$
as above, we set:
$$
\e^* M := M(e_{i,j}) = e_{1,1}^{r_{1,1}} e_{1,2}^{r_{1,2}} e_{2,1}^{r_{2,1}} e_{2,2}^{r_{2,2}}.
$$
In particular $\e^* \Phi = \phi := \dfrac{e_{1,1} e_{2,2}}{e_{1,2} e_{2,1}}$, which is $H$-invariant.
(Of course ambiguous notation $\phi$ depends on the particular $\e$, we shall try to avoid confusions.)

\begin{lem}
  (i) $\Af\left(V^{(*)}(\e)\right)^H = \C[\text{all~} x_{i,j},y_{i,j}][\phi^{\pm 1}]$. \\
  (ii) There is a geometric quotient:
  $$
  V^{(*)}(\e)/H = \Pr(V^{(*)})_\e \times \Cs \text{~where~} \Pr(V^{(*)})_\e := \prod \Pr(V_{i,j})_{e_{i,j}}.
  $$
\end{lem}
\Proof
  (i) Every $P \in \Af\left(V^{(*)}(\e)\right)$ admits a unique expansion as:
  $$
  P = \sum_{M \in \M^*} P_M(\text{all~} x_{i,j},y_{i,j}) \, \e^* M.
  $$
  The factors $P_M$ are $H$-invariants, while the factors $\e^* M$ are semi-invariants, the only
  invariant ones being the powers of $\e^* \Phi$. (This is actually the decomposition of a linear
  representation flowing from the semi-simplicity of $H$.) \\
  (ii) We know that the categorical quotient has algebra $\Af\left(V^{(*)}(\e)\right)^H$ (affine by
  reductive). It is a geometric quotient by \cite[p 9]{Brion}, because orbits are closed and
  stabilizers are trivial (this was proved before).
\Finpr

\begin{prop}
  \label{prop:quotientV*/H}
  The quotient $V^{(*)}/H$ is (the total space of) a fibration with fiber $\Cs$ over the base
  $\Pr(V^{(*)}) := \prod \Pr(V_{i,j}) \simeq \Pr^1(\C)^4$. It is a geometric quotient.
\end{prop}
\Proof
The affine charts $\Pr(V^{(*)})_\e$, $\Pr(V^{(*)})_{\f}$ on $\Pr(V^{(*)})$ are patched along their
intersection $\Pr(V^{(*)})_\e \cap \Pr(V^{(*)})_{\f}$ by the transition functions $(e_{i,j}/f_{i,j})_{i,j}$.
The corresponding transition functions on the $\Cs$ factors are the $\dfrac{\e^*\Phi}{\f^*\Phi} \cdot$
The possibility of patching geometric quotients is guaranteed by
\cite[prop. 3.10 (b), p 71]{NewsteadTata}. One just has to check first that there is
a well defined morphism $V^{(*)} \rightarrow Y$, $Y$ the candidate quotient (the fiber space described
above); and then to check that it localizes to geometric quotients on a covering of $Y$. The first step
is easy, the second step comes from the lemma.
\Finprcourt

\begin{cor}
  \label{cor:quotientV*/H}
  The space $V^{(*)}/H$ is separated.
\end{cor}

% 2.3
 
\subsection{Degenerate JS case: one zero \emph{required}}

We fix (for instance) the position of the zero coefficient at $(1,1)$, so that here:
$$
S = \left(\{1,2\} \times \{1,2\}\right) \setminus \{(1,1)\} = \{(1,2),(2,1),(2,2)\}.
$$
We thus look at the $H$-action on:
$$
V_S = V^0 = \{0\} \times V_{1,2}^* \times V_{2,1}^* \times V_{2,2}^*.
$$
We shall write $\e' := (e_{1,2},e_{2,1},e_{2,2})$, each relevant $e_{i,j} \in V_{i,j}^\vee \setminus \{0\}$
(here and later, ``relevant'' means that $(i,j) = (1,1)$ is omitted). Then:
$$
V_S(\e') = V^0(\e') =  \{0\} \times V_{1,2}(e_{1,2}) \times V_{2,1}(e_{2,1}) \times V_{2,2}(e_{2,2}).
$$
Those open subsets of $V^0$ form an affine covering. The corresponding affine algebras are (same
notations as in the non degenerate JS case):
$$
\Af\left(V^0(\e')\right) = 
\C[x_{1,2},y_{1,2},x_{2,1},y_{2,1},x_{2,2},y_{2,2}][e_{1,2}^{\pm 1},e_{2,1}^{\pm 1},e_{2,2}^{\pm 1}].
$$

\begin{lem}
  (i) $\Af\left(V^0(\e')\right)^H = \C[x_{1,2},y_{1,2},x_{2,1},y_{2,1},x_{2,2},y_{2,2}]$. \\
  (ii) There is a geometric quotient:
  $$
  V^0(\e')/H = \Pr(V^0)_{\e'} :=
  \Pr(V_{1,2})_{e_{1,2}} \times \Pr(V_{2,1})_{e_{2,1}} \times \Pr(V_{2,2})_{e_{2,2}}.
  $$
\end{lem}
\Proof
  (i) Every $P \in \Af\left(V^0(\e')\right)$ admits a unique expansion as a sum
  $\sum P_M(\text{relevant~} x_{i,j},y_{i,j}) M(\text{relevant~} e_{i,j})$, but here the only invariant
  monomial $M$ is the trivial one (because the relevant $\lambda_i/\mu_j$ are three independent complex
  numbers). \\
  (ii) Follows as in the non degenerate case.
\Finpr

\begin{prop}
  There is a geometric quotient:
  $$
  V^0/H = \Pr(V^0) := \Pr(V_{1,2}) \times \Pr(V_{2,1}) \times \Pr(V_{2,2}) \simeq \Pr^1(\C)^3.
  $$
\end{prop}

% 2.4

\subsection{A candidate craddle for partial patching: one zero \emph{allowed}}

We fix again the position of the allowed zero coefficient at $(1,1)$, so that we are looking at the
$H$-action on:
$$
V' := V_{1,1} \times V_{1,2}^* \times V_{2,1}^* \times V_{2,2}^* = V^{(*)} \sqcup V^0.
$$
For $\e' := (e_{1,2},e_{2,1},e_{2,2})$, each relevant $e_{i,j} \in V_{i,j}^\vee \setminus \{0\}$, we set:
$$
V'(\e') = V_{1,1} \times V_{1,2}(e_{1,2}) \times V_{2,1}(e_{2,1}) \times V_{2,2}(e_{2,2}).
$$
Those open subsets of $V'$ form an affine covering. To describe the corresponding affine algebras, we
use the same notations as before, plus a new one:
$$
x'_{1,1} = u_{1,1} \dfrac{e_{2,2}}{e_{1,2} e_{2,1}}, \quad
y'_{1,1} = v_{1,1} \dfrac{e_{2,2}}{e_{1,2} e_{2,1}} \cdot
$$
They are $H$-invariant and, whenever $\e$ is involved, related by equalities:
$$
\begin{cases} x'_{1,1} = x_{1,1} \, \e^* \Phi, \\ y'_{1,1} = y_{1,1} \, \e^* \Phi, \end{cases}
\Longrightarrow \alpha_{1,1} x'_{1,1} + \beta_{1,1} y'_{1,1} = \e^* \Phi.
$$
The affine algebras are then given by:
$$
\Af\left(V'(\e')\right) = 
\C[x'_{1,1},y'_{1,1},x_{1,2},y_{1,2},x_{2,1},y_{2,1},x_{2,2},y_{2,2}][e_{1,2}^{\pm 1},e_{2,1}^{\pm 1},e_{2,2}^{\pm 1}].
$$
Then by the same arguments as before:

\begin{lem}
  (i) $\Af\left(V'(\e')\right)^H = \C[x'_{1,1},y'_{1,1},x_{1,2},y_{1,2},x_{2,1},y_{2,1},x_{2,2},y_{2,2}]$. \\
  (ii) There is a categorical quotient:
  $$
  V'(\e')/H = \Pr(V')_{\e'} :=
  \C^2 \times \Pr(V_{1,2})_{e_{1,2}} \times \Pr(V_{2,1})_{e_{2,1}} \times \Pr(V_{2,2})_{e_{2,2}}.
  $$
\end{lem}
\Proof
  (i) Is easy as in the previous calculations. \\
  (ii) Follows by the general statement about a reductive group acting on an affine variety.
\Finpr

\begin{prop}
  There is a categorical quotient $V'/H$; it is a vector bundle of rank $2$ over $\Pr(V^0)$. It is
  actually the cartesian square of a line bundle with transition functions the
  $\dfrac{f_{2,2}/(f_{1,2} f_{2,1})}{e_{2,2}/(e_{1,2} e_{2,1})} \cdot$ 
\end{prop}
\Proof
  Indeed, those are the laws of transformation independently of variables $x'_{1,1},y'_{1,1}$.
  We again use the patching argument from \cite{NewsteadTata} already used in the proof of
  proposition \ref{prop:quotientV*/H}.
\Finpr

\begin{rmk}
  The above bundle is actually an ``external tensor product'' of line bundles $\mathcal{O}(\pm 1)$.
\end{rmk}

% 2.5

\subsection{The patching}

To see how the quotients $V^{(*)} \rightarrow V^{(*)}/H$ and $V^0 \rightarrow V^0/H$ embed into the
quotient $V' \rightarrow V'/H$, we use two different principles.

% 2.5.1

\subsubsection{Non degenerate part}

Let $\e'$ be extracted from $\e$ omitting $e_{1,1}$. We have an open immersion of affine varieties
$V^{(*)}(\e) \rightarrow V'(\e')$ induced by the dual morphism of their affine algebras: 
$$  
\C[x'_{1,1},y'_{1,1},\text{all other~} x_{i,j},y_{i,j}][e_{1,2}^{\pm 1},e_{2,1}^{\pm 1},e_{2,2}^{\pm 1}]
\rightarrow
\C[\text{all~} x_{i,j},y_{i,j}][\text{all~} e_{i,j}^{\pm 1}]
$$
From the equality 
$u_{1,1} \dfrac{e_{2,2}}{e_{1,2} e_{2,1}} =\dfrac{u_{1,1}}{e_{1,1}} \, \dfrac{e_{1,1} e_{2,2}}{e_{1,2} e_{2,1}}$
and the similar one with $v_{1,1}$, we see that it sends $x'_{1,1}$ to $x_{1,1} \phi$ and $y'_{1,1}$
to $y_{1,1} \phi$. It restricts to the $H$-invariant subalgebras as the morphism sending $x'_{1,1}$ to
$x_{1,1} \phi$ and $y'_{1,1}$ to $y_{1,1} \phi$:
$$
\C[x'_{1,1},y'_{1,1},\text{all other~} x_{i,j},y_{i,j}]
\rightarrow
\C[\text{all~} x_{i,j},y_{i,j}][\phi^{\pm 1}].
$$
The open immersion being $H$-equivariant, it induces an immersion of the quotients from
$\Pr(V^{(*)})_\e \times \Cs$ to
$\C^2 \times \Pr(V_{1,2})_{e_{1,2}} \times \Pr(V_{2,1})_{e_{2,1}} \times \Pr(V_{2,2})_{e_{2,2}}$ dual to the above
morphism of algebras. Patching up, we get an open immersion:
$$
V^{(*)}/H = \Pr(V^{(*)}) \times \Cs \rightarrow V'/H = \text{plane bundle over~} \Pr(V^0),
$$
described in coordinates by the above morphism. So $V^{(*)} \rightarrow V^{(*)}/H$ is the restriction of
the categorical quotient $V' \rightarrow V'/H$ to an invariant open subset.

% 2.5.2

\subsubsection{Degenerate part}

Let $\e'$ be as usual. We have a closed immersion of affine varieties $V^0(\e') \rightarrow V'(\e')$
induced by the dual morphism of their affine algebras sending $x'_{1,1}$ and $y'_{1,1}$ to $0$:
$$  
\C[x'_{1,1},y'_{1,1},\text{all other~}][e_{1,2}^{\pm 1},e_{2,1}^{\pm 1},e_{2,2}^{\pm 1}]
\rightarrow
\C[x_{1,2},y_{1,2},x_{2,1},y_{2,1},x_{2,2},y_{2,2}][e_{1,2}^{\pm 1},e_{2,1}^{\pm 1},e_{2,2}^{\pm 1}].
$$
It restricts to a similar morphism of the invariant subalgebras:
$$
\C[x'_{1,1},y'_{1,1},x_{1,2},y_{1,2},x_{2,1},y_{2,1},x_{2,2},y_{2,2}]
\rightarrow
\C[x_{1,2},y_{1,2},x_{2,1},y_{2,1},x_{2,2},y_{2,2}].
$$
The dual morphism is a closed immersion to the $\{0\} \subset \C^2$ section:
$$
\Pr(V_{1,2})_{e_{1,2}} \times \Pr(V_{2,1})_{e_{2,1}} \times \Pr(V_{2,2})_{e_{2,2}}
\rightarrow
\C^2 \times \Pr(V_{1,2})_{e_{1,2}} \times \Pr(V_{2,1})_{e_{2,1}} \times \Pr(V_{2,2})_{e_{2,2}}
$$
Patching up, we get a closed immersion to the null section:
$$
V^0/H = \Pr(V^0) \rightarrow V'/H = \text{plane bundle over~} \Pr(V^0).
$$

%%%%%%%%%%%%%%%%%%%%%%%%%%%%%%%%%%%%%%%%%%%%%%%%%%%%%%%%%%%%%%%%%%%%%%%%%%%%%

% 3

\section{Quadrics in the non degenerate JS case and functions on them}

% 3.1

\subsection{General notations and conventions}
\label{Generalnotationsandconventions}

This section aims at translating in terms of (bi)linear algebra the properties of the monodromy
matrices in $F_{R,S,\underline{x}}$ (\cf\ the paper \cite{ORS}). This yields some objects subject to
some axioms, for which we briefly recall each time the justification from \cite{ORS}. From now on
in this article, $n = \mu = 2$ so that $N = 4$. The following conditions will be systematically
assumed:
\begin{itemize}
\item{(NR)} \emph{Strong non resonancy:}
  \begin{align*}
  \forall i,j \in \{1,\ldots,2\} \;,\; & i \neq j \Longrightarrow
  \left(\rho_i \not\equiv \rho_j \text{~and~} \sigma_i \not\equiv \sigma_j \right), \\
  \forall k,l \in \{1,\ldots,4\} \;,\; & k \neq l \Longrightarrow x_k \not\equiv x_l.
  \end{align*}
\item{(FR)} \emph{Fuchs relation:}
  $$
  (-1)^N x_1 \cdots x_N = \dfrac{\rho_1 \ldots \rho_n}{\sigma_1 \ldots \sigma_n} \cdot
  $$
\end{itemize}
(Recall that we denote $\equiv$ the congruence modulo $q^\Z$ in $\Cs$.) Also the following conditions
will be frequently invoked:
\begin{itemize}
\item{(NS)} \emph{Non splitting:}
  $$
  \forall i,j \in \{1,2\} \;,\; \rho_i/\sigma_j \not\equiv x_1 x_2.
  $$
\item{(SC)} \emph{Special condition:}
  $$
  \dfrac{\rho_1}{\sigma_1} \not\equiv \dfrac{\rho_2}{\sigma_2}
  \quad \text{~and~} \quad 
  \dfrac{\rho_2}{\sigma_1} \not\equiv \dfrac{\rho_1}{\sigma_2}.
  $$
\end{itemize}

For comments around (NR), (FR) and (NS), see \cite[5.1]{ORS}; for (SC), see \cite[4.4]{ORS}.

% 3.1.1

\subsubsection{Evaluation and linear forms}
\label{subsubsection:evaluationandlinearforms}

We are given on each $V_{i,j}$ four particular non trivial linear forms respectively corresponding to
the evaluation maps $m \mapsto m(x_k)$, $k = 1,2,3,4$. They define four families in the $V_{i,j}^\vee$. \\
{\textsl{Axiom. - }} For each $i,j \in \{1,2\}$, the four given forms are pairwise linearly independent. \\
{\textsl{Explanation. - }} If $m_{i,j}(x_k) = m_{i,j}(x_l) = 0$ with $m_{i,j} \neq 0$, then
$x_k x_l \equiv \rho_i/\sigma_j$, contradicting (NS). \\

With the previous notations for bases $(u_{i,j},v_{i,j})$ of the $V_{i,j}$, we can decide for instance
that $u_{i,j}$ and $v_{i,j}$ respectively correspond to the evaluations at $x_1$ and $x_2$.
We denote $w_{i,j}$, resp. $w'_{i,j}$, the evaluation at $x_3$, resp. at $x_4$ and express them as:
$$
\begin{cases}
w_{i,j} = \lambda_{i,j} u_{i,j} + \mu_{i,j} v_{i,j}, \quad \lambda_{i,j}, \mu_{i,j} \neq 0, \\
w'_{i,j} = \lambda'_{i,j} u_{i,j} + \mu'_{i,j} v_{i,j}, \quad \lambda'_{i,j}, \mu'_{i,j} \neq 0,
\end{cases}
\quad \text{with moreover the conditions~} \lambda_{i,j} \mu'_{i,j} - \lambda'_{i,j} \mu_{i,j} \neq 0.
$$

Through the natural projections from $V := V_{1,1} \times V_{1,2} \times V_{2,1} \times V_{2,2}$ to each
$V_{i,j}$, and the dual injections $V_{i,j}^\vee \rightarrow V^\vee$, we may (and will) intepret the above
as linear forms on $V$. In this way, we obtain four linear maps $\u$, $\v$, $\w$, and $\w'$, from $V$
to $\Mat_2(\C)$. 

% 3.1.2

\subsubsection{Determinant and bilinear forms}
\label{subsubsection:3.1.2}

Corresponding to the evaluations
$(m_{1,1} m_{2,2} - m_{1,2} m_{2,1})(x_k) = m_{1,1} m_{2,2}(x_k)  - m_{1,2} m_{2,1}(x_k)$ of the determinant,
we define on $V$ twelve bilinear forms as follows:
$$
\begin{matrix}
A_+ := u_{1,1} u_{2,2}, & A_- := u_{1,2} u_{2,1}, & A := A_+ - A_- = \det \u, \\
B_+ := v_{1,1} v_{2,2}, & B_- := v_{1,2} v_{2,1}, & B := B_+ - B_- = \det \v, \\
C_+ := w_{1,1} w_{2,2}, & C_- := w_{1,2} w_{2,1}, & C := C_+ - C_- = \det \w, \\
C'_+ := w'_{1,1} w'_{2,2}, & C'_- := w'_{1,2} w'_{2,1}, & C' := C'_+ - C'_- = \det \w'.
\end{matrix}
$$
{\textsl{Axiom. - }} There exists $\alpha,\beta,\gamma \in \Cs$ such that
$C_+ = \alpha A_+ + \beta B_+ + \gamma C_+$ and $C_- = \alpha A_- + \beta B_- + \gamma C_-$.
It follows that $C = \alpha A + \beta B + \gamma C$. \\
{\textsl{Explanation. - }} Fuchs relations $x_1 x_2 x_3 x_4 \equiv \dfrac{\rho_1 \rho_2}{\sigma_1 \sigma_2}$
imply that, if $m \in V_{4,\frac{\rho_1 \rho_2}{\sigma_1 \sigma_2}}$ vanishes at three of the $x_k$ then it vanishes
at the fourth. So this is really a relation among linear forms on $V_{4,\frac{\rho_1 \rho_2}{\sigma_1 \sigma_2}}$.

% 3.1.3

\subsubsection{Application to $F_{R,S,\x}$ and $\F_{R,S,\x}$}
\label{subsubsection:FtocalFquotientgeom}

We want to translate the condition (from \cite{ORS}) that $\det\ M$ vanishes at the $x_k$, $k = 1,\ldots,4$,
but does not vanish identically; equivalenly, because of the Fuchs relations and the fact that
$\det\ M \in V_{4,\frac{\rho_1 \rho_2}{\sigma_1 \sigma_2}}$: the matrices $M(x_k)$ all have rank $1$ (otherwise
$\det\ M$ would have a multiple zero). So this means that $\det\ M(x_k) = 0$ but no $M(x_k) = 0$. Also
it implies that $M$ has neither a null line nor a null column. At any rate, the subset $F \subset V$
is locally closed ($\det \neq 0$, open condition, within $\det\ M(x_k) = 0$, closed conditions) and
also stable under the action of $H$. We write the space of interest (in which we recognize $F_{R,S,\x}$):
$$
F := A^{-1}(0) \cap B^{-1}(0) \cap C^{-1}(0) \setminus W, \quad W := {\det}^{-1}(0).
$$

\begin{lem}
  $F \subset V^{(*)} = V_{1,1}^* \times V_{1,2}^* \times V_{2,1}^* \times V_{2,2}^*$.
\end{lem}
\Proof
  If for instance $m_{1,1} = 0$, then $m_{1,2} m_{2,1} \neq 0$ (lest $\det\ M$ be $0$) but it vanishes
  at all $x_k$ (because $\det\ M$ does), so at least one of $m_{1,2}$, $m_{2,1}$, vanishes at at least
  two of the $x_k$, contradicting (NS). A similar argument applies to all $m_{i,j}$.
\Finpr

The image $\F$ of $F$ under the geometric quotient map $V^{(*)} \rightarrow V^{(*)}/H$ (in which we
recognize $\F_{R,S,\x}$) is a subvariety of $V^{(*)}/H$. 

\begin{thm}
  \label{thm:LaStructureDeF}
  (i) The map $F_{R,S,\x} \rightarrow \F_{R,S,\x}$ is a geometric quotient. \\
  (ii) The monodromy data space $\F_{R,S,\x}$ is a separated variety.
\end{thm}
\Proof
(i) Using the lemma, it is an application of corollary \ref{cor:V_S^{(*)}/H}. \\
(ii) This follows from corollary \ref{cor:quotientV*/H}.
\Finpr

Another useful consequence of the lemma is:

\begin{cor}
  For $M \in F$, and $k \neq l$, matrices $M(x_k)$ and $M(x_l)$ cannot have a zero at the same
  position.
\end{cor}
\Proof
  If for instance $m_{1,1}$ vanishes at $x_k$ and $x_l$, $k \neq l$, then, by (NS), $m_{1,1} = 0$,
  which was just excluded.
\Finprcourt

% 3.2

\subsection{Rank $1$ matrices in $\Mat_2(\C)$ and the Segre embedding}

% 3.2.1

\subsubsection{Projective ``coordinates'' for rank $1$ matrices}
\label{subsubsection:projectivecoordinatesrank1matrices}

The image of the map $(C,L) \mapsto C L$ from $\Mat_{2,1}(\C) \times \Mat_{1,2}(\C)$ to $\Mat_2(\C)$ is
the subset of singular matrices. The restriction to $\Mat_{2,1}(\C)^* \times \Mat_{1,2}(\C)^*$ arrives
in $\Mat_2(\C)^*$, its image is the set $\Mat_2(\C)_1$ of matrices of rank $1$. Note that the latter
is locally closed in $\Mat_2(\C)$, so it is a quasi-affine variety, actually a $\Cs$-cone, so the
notation $\Pr\left(\Mat_2(\C)_1\right)$ is licit (it denotes a quasi-projective variety). Going to
the projective quotients, we obtain:
$$
\Pr\left(\Mat_{2,1}(\C)\right) \times \Pr\left(\Mat_{1,2}(\C)\right) \rightarrow
\Pr\left(\Mat_2(\C)\right),
$$
which can be identified to the \emph{Segre embedding} $\Pr^1(\C) \times \Pr^1(\C) \rightarrow \Pr^3(\C)$,
an isomorphism of the source with a smooth projective quadric surface, the \emph{Segre surface}.
$$
\xymatrix{
  \Mat_{2,1}(\C)^* \times \Mat_{1,2}(\C)^* \ar@<0ex>[r] \ar@<0ex>[d] &
  \Mat_2(\C)_1 \ar@<0ex>[r] \ar@<0ex>[d]  \ar@{.>}[ddl]|\hole  & \Mat_2(\C)^* \ar@<0ex>[d]  \\
  \Pr\left(\Mat_{2,1}(\C)\right) \times \Pr\left(\Mat_{1,2}(\C)\right) \ar@<0ex>[r] \ar@<0ex>[d] &
  \Pr\left(\Mat_2(\C)_1\right) \ar@<0ex>[r] \ar@<0ex>[d]^\simeq &
  \Pr\left(\Mat_2(\C)\right) \ar@<0ex>[d] \\
  \Pr^1(\C) \times \Pr^1(\C) \ar@<0ex>[r]^\simeq & \text{Segre surface} \ar@<0ex>[r]  & \Pr^3(\C)
}
$$

The Segre surface corresponds isomorphically to the image $\Pr\left(\Mat_2(\C)_1\right)$ of
$\Mat_2(\C)_1$ in $\Pr\left(\Mat_2(\C)\right)$. We write $(\rho,\rho')$ the composite map (the dotted
arrow in the previous diagram):
$$
\Mat_2(\C)_1 \rightarrow \Pr\left(\Mat_2(\C)_1\right) \rightarrow \text{Segre surface}
\rightarrow \Pr^1(\C) \times \Pr^1(\C), \quad
C L \mapsto (\text{class of~} C,\text{class of~} L).
$$
So $\rho$ and $\rho'$ are both regular maps from the variety $\Mat_2(\C)_1$ to $\Pr^1(\C)$. \\

Concretely, let $A := \begin{pmatrix} a_{1,1} & a_{1,2} \\ a_{2,1} & a_{2,2} \end{pmatrix} \in \Mat_2(\C)_1$.
Then $\rho(A),\rho'(A) \in \Pr^1(\C)$ are given by the formulas:
\begin{equation}
\label{eqn:projectivecoordinatesrank1matrices}
\rho(A) = [a_{1,1} : a_{1,2}] = [a_{2,1} : a_{2,2}] \text{~and~}
\rho'(A) = [a_{1,1} : a_{2,1}] = [a_{1,2} : a_{2,2}].
\end{equation}
In the first formula, one of the expressions $[a_{1,1} : a_{1,2}]$, $[a_{2,1} : a_{2,2}]$ may be undefined
(if one line of $A$ vanishes): then just dismiss the corresponding equality; same thing for the second
formula.

% 3.2.2

\subsubsection{Special points and mixed projective coordinates}

It is immediate by inspection of $C L$ that:
\begin{itemize}
\item The closed subset $\rho^{-1}(0)$ of $\Mat_2(\C)_1$ consists of matrices with trivial first line;
\item The closed subset $\rho^{-1}(\infty)$ of $\Mat_2(\C)_1$ consists of matrices with trivial second line;
\item The closed subset $\rho'^{-1}(0)$ of $\Mat_2(\C)_1$ consists of matrices with trivial first column;
\item The closed subset $\rho'^{-1}(\infty)$ of $\Mat_2(\C)_1$ consists of matrices with trivial second
  column.
\end{itemize}
Now, let $M_i := C_i L_i \in \Mat_2(\C)_1$, $i = 1,2$, be such that neither $\rho(M_1) = \rho(M_2) = 0$
nor $\rho(M_1) = \rho(M_2) = \infty$. Denoting $C_i := \begin{pmatrix} f_i \\ g_i \end{pmatrix} \neq 0$,
one easily checks that $f_1 g_2 = f_2 g_1 = 0$ would imply either $f_1 = f_2 = 0$ or $g_1 = g_2 = 0$,
which are both excluded by assumption. Therefore $[f_1 g_2:f_2 g_1] \in \Pr^1(\C)$ is well defined.
Using the fact that $C L = C' L' \Rightarrow C',C$ proportional, we see that $[f_1 g_2:f_2 g_1]$ depends
on $M$ only. We call it $\Pi(M_1,M_2)$. Using lines instead of columns we define similarly $\Pi'(M_1,M_2)$.
This way, we define two regular maps:
$$
\begin{cases}
  \Pi: \Mat_2(\C)_1 \times \Mat_2(\C)_1 \setminus
  \left((\rho^{-1}(0) \times \rho^{-1}(0)) \cup (\rho^{-1}(\infty) \times \rho^{-1}(\infty))\right)
  \rightarrow \Pr^1(\C), \\
  \Pi': \Mat_2(\C)_1 \times \Mat_2(\C)_1 \setminus
  \left((\rho'^{-1}(0) \times \rho'^{-1}(0)) \cup (\rho'^{-1}(\infty) \times \rho'^{-1}(\infty))\right)
  \rightarrow \Pr^1(\C).
\end{cases}
$$

% 3.3

\subsection{The maps $\rho_k$, $k = 1,\ldots,4$, and $[\rho]$ on $\F$}
\label{subsection:Therhokand[rho]}

We try to adapt the first approach in \cite{JR}. Independently of the more or less ``axiomatic''
presentation adopted here (unification will come later), we summarize some basic facts about $F$
and its quotient $\F$:
\begin{enumerate}
\item $F$ is a subset of the $8$-dimensional linear space $V$, defined by three homogeneous quadratic
  equations $\det M(x_k) = 0$, $k = 1,2,3$ (the fourth one is a consequence by Fuchs relation (FR)) and
  one inequation $\det M \neq 0$. Therefore $F$ is a quasi-affine variety of dimension $\geq 5$. 
\item Because of the non splitting condition (NS), we know that $F \subset V^{(*)}$. It is also clearly
  stable under the action of $H$. Since we know that there is a geometric quotient $V^{(*)}/H$, we see
  that $\F$ is the image of $F$ in $V^{(*)}/H$ and that it is a quasi-affine variety of dimension $\geq 2$.
\end{enumerate}

% 3.3.1

\subsubsection{Definition of the $\rho_k$ and $[\rho]$}

Let $\Gamma,\Delta \in \Ddc$. Then, with the usual notations $M = C L$:
$$
\rho(\Gamma M) = \rho(\Gamma C L) = \rho((\Gamma C) L) = \text{~class of~} L = \rho(C L) = \rho(M),
$$
while:
$$
\rho(M \Delta) = \rho(C L \Delta) = \rho(C (L \Delta)) = \text{~class of~} L  \Delta =
\chi(\Delta).\rho(C L) = \chi(\Delta).\rho(M),
$$
where\footnote{Of course this character $\chi$ is in no way related to the function $\chi$ on graphs
introduced in \ref{subsubsection:stabilizers}.} $\chi(\Diag(d,d')) := d/d'$ and we have $\Cs$ act
on $\Pr^1(\C)$ by $c.[x:y] := [c x:y]$. Symmetrically:
$$
\rho'(\Gamma M) = \chi(\Gamma).\rho'(M), \quad \rho'(M \Delta) = \rho'(M).
$$

We have defined $F$ so that it is sent into $\Mat_2(\C)_1$ by either one of $\u$, $\v$, $\w$, $\w'$.
Therefore we have a well defined regular map:
$$
(\rho_1,\rho_2,\rho_3,\rho_4) :=
(\rho \circ \u,\rho \circ \v,\rho \circ \w,\rho \circ \w'): F \rightarrow \left(\Pr^1(\C)\right)^4.
$$
It is equivariant under the action of $H$ in the following sense:
$$
M \mapsto (r_1,r_2,r_3,r_4) \in \left(\Pr^1(\C)\right)^4 \Longrightarrow
\Gamma M \Delta^{-1} \mapsto (c.r_1,c.r_2,c.r_3,c.r_4),
\text{~where~} c := \dfrac{\chi(\Gamma)}{\chi(\Delta)} \cdot
$$
Recall $\F$, the image of $F$ under $V^{(*)} \rightarrow V^{(*)}/H$. The above induces a map:
$$
[\rho]: \F \rightarrow \left(\Pr^1(\C)\right)^4/\Cs.
$$

\begin{prop}
  The map $[\rho]$ is injective.
\end{prop}
\Proof
  We want to show that if $M,N \in F$ are such that $\rho_k(N) = c \rho_k(M)$ for some $c \in \Cs$
  and $k = 1,\ldots,4$, then $N = \Gamma M \Delta$ for some $\Gamma,\Delta \in \Ddc$. First, changing
  $M$ to $M \Delta$ for some $\Delta \in \Ddc$ (\emph{e.g.} $\Diag(1,c)$), we can assume that $c = 1$.
  Then, from the lemma that follows applied in turn to the four equalities $\rho_k(M) = \rho_k(N)$, \ie\
  $\rho(M(x_k)) = \rho(N(x_k))$, we draw that for all $i,j = 1,2$:
  $$
  m_{i,1} n_{j,2} - m_{i,2} n_{j,1} \in V_{4,\rho_i \rho_j/(\sigma_1 \sigma_2)} \text{~vanishes at all~} x_k.
  $$
  For $i = j$, using properties (NR) (non resonancy) and (FR), we have:
  $$
  \rho_i^2/(\sigma_1 \sigma_2) \not\equiv \rho_1 \rho_2/(\sigma_1 \sigma_2) \Longrightarrow
  \rho_i^2/(\sigma_1 \sigma_2) \not\equiv x_1 x_2 x_3 x_4,
  $$
  so an element of $V_{4,\rho_i^2/(\sigma_1 \sigma_2)}$ which vanishes at all $x_k$ must be $0$. Therefore we
  have $m_{i,1} n_{i,2} = m_{i,2} n_{i,1}$ for $i = 1,2$. Let $\gamma_i := n_{i,1}/m_{i,1} = n_{i,2}/m_{i,2}$
  (recall that all $m_{i,j}$ and $n_{i,j}$ are $\neq 0$). Both $\gamma_i$ are non zero elliptic functions
  of order $\leq 2$. But, because of (NR), $m_{i,1}$ and $m_{i,2}$ cannot have two common zeroes. So
  actually $\gamma_1,\gamma_2 \in \Cs$ and $\Gamma := \Diag(\gamma_1,\gamma_2) \in \Ddc$ is such that
  $N = \Gamma M$.  
\Finprcourt

\begin{lem}
  Let $A := (a_{i,j}), B  := (b_{i,j}) \in \Mat_2(\C)_1$. Then:
  $$
  \rho(A) = \rho(B) \Longleftrightarrow \forall i,j = 1,2 \;,\; a_{i,1} b_{j,2} - a_{i,2} b_{j,1} = 0.
  $$
\end{lem}
\Proof
Note that, generally speaking, $[a_1:a_2] = [b_1:b_2] \Leftrightarrow a_1 b_2 = a_2 b_1$ and apply
the ``concrete description'' \eqref{eqn:projectivecoordinatesrank1matrices} of $\rho(A),\rho(B)$
at the end of \ref{subsubsection:projectivecoordinatesrank1matrices}.
\Finpr

We would like to consider $[\rho]$ as a regular map, but it is not clear in what sense the quotient
$\left(\Pr^1(\C)\right)^4/\Cs$ is defined beyond its obvious set-theoretic sense. For instance, generic
orbits under the $\Cs$ action on $\left(\Pr^1(\C)\right)^4$ have dimension $1$ while there are $16$
invariant points $\Theta_4 := \{0,\infty\}^4$, etc. However, we saw in \ref{subsubsection:NSGQ} that,
extracting the set $\Theta_4$, we obtain a (non separated) geometric quotient
$\mathcal{K}_4 := \left(\left(\Pr^1(\C)\right)^n \setminus \Theta_n\right)/\Cs$, which is moreover
a toric variety (\cf\ \ref{subsubsection:examplestor}).

% 3.3.2

\subsubsection{The image of $[\rho]$}
\label{subsubsection:3.3.2}

Instead of $\rho_1,\rho_2,\rho_3,\rho_4$, we shall write $\rho,\sigma,\tau,\tau'$ in harmony with our
notations $\u,\v,\w,\w'$ (and in slight contradiction with the previous use of $\rho$, but it does not
seem to matter much). We do the case of $F_1$, so $(\rho,\sigma,\tau,\tau')$ sends it to
$Y'_1 = \Cs \times \Cs \times \C \times \C_\infty$. We have defining equalities:
$$
\begin{cases} u_{1,2} = \rho u_{1,1}, \\ u_{2,2} = \rho u_{2,1}, \end{cases}
\begin{cases} v_{1,2} = \sigma v_{1,1}, \\ v_{2,2} = \sigma v_{2,1}, \end{cases}
\begin{cases} w_{1,2} = \tau w_{1,1}, \\ w_{2,2} = \tau w_{2,1}, \end{cases}
\begin{cases} w'_{1,2} = \tau' w'_{1,1}, \\ w'_{2,2} = \tau' w'_{2,1}. \end{cases}
$$
Only in the special case that $\tau' = \infty$ must we rather use $\tau'' := \tau'^{-1}$ (to be considered
at $\tau'' = 0$) and replace the last pair of equalities by
$\begin{cases} \tau'' w'_{1,2} = w'_{1,1}, \\ \tau'' w'_{2,2} = w'_{2,1}, \end{cases}$. We also recall that
$\begin{cases}
  w_{i,j} = \lambda_{i,j} u_{i,j} + \mu_{i,j} v_{i,j}, \\ w'_{i,j} = \lambda'_{i,j} u_{i,j} + \mu'_{i,j} v_{i,j},
\end{cases}$
for $i,j = 1,2$, whence:
$$
\begin{cases}
  \lambda_{1,2} \rho u_{1,1} + \mu_{1,2} \sigma v_{1,1} = \tau (\lambda_{1,1} u_{1,1} + \mu_{1,1} v_{1,1}), \\
  \lambda_{2,2} \rho u_{2,1} + \mu_{2,2} \sigma v_{2,1} = \tau (\lambda_{2,1} u_{2,1} + \mu_{2,1} v_{2,1}),
\end{cases}
\text{and~}
\begin{cases}
  \lambda'_{1,2} \rho u_{1,1} + \mu'_{1,2} \sigma v_{1,1} = \tau'(\lambda'_{1,1} u_{1,1} + \mu'_{1,1} v_{1,1}), \\
  \lambda'_{2,2} \rho u_{2,1} + \mu'_{2,2} \sigma v_{2,1} = \tau'(\lambda'_{2,1} u_{2,1} + \mu'_{2,1} v_{2,1}).
\end{cases}
$$
Of course, the last two equalities must be modified in an obvious way at $\tau' = \infty$ (\ie\ put a
$\tau''$ factor at left instead of a $\tau'$ factor at right). \\
We translate the above into linear systems; first we do that in the case $\tau' \neq \infty$:
$$
\begin{cases}
  (\lambda_{1,2} \rho - \lambda_{1,1} \tau) u_{1,1} + (\mu_{1,2} \sigma - \mu_{1,1} \tau) v_{1,1} = 0, \\
  (\lambda'_{1,2} \rho - \lambda'_{1,1} \tau') u_{1,1} + (\mu'_{1,2} \sigma - \mu'_{1,1} \tau') v_{1,1} = 0,
\end{cases}
\text{and~}
\begin{cases}
  (\lambda_{2,2} \rho - \lambda_{2,1} \tau) u_{2,1} + (\mu_{2,2} \sigma - \mu_{2,1} \tau) v_{2,1} = 0, \\
  (\lambda'_{2,2} \rho - \lambda'_{2,1} \tau') u_{2,1} + (\mu'_{2,2} \sigma - \mu'_{2,1} \tau') v_{2,1} = 0.
\end{cases}
$$
Noting that we can have neither $(u_{1,1},v_{1,1}) = (0,0)$ nor $(u_{2,1},v_{2,1}) = (0,0)$ (again because
of (NR) and (NS)), we deduce the vanishing of the determinants:
$$
\begin{cases}
(\lambda_{1,2} \rho - \lambda_{1,1} \tau) (\mu'_{1,2} \sigma - \mu'_{1,1} \tau') -
(\lambda'_{1,2} \rho - \lambda'_{1,1} \tau') (\mu_{1,2} \sigma - \mu_{1,1} \tau) = 0, \\
(\lambda_{2,2} \rho - \lambda_{2,1} \tau) (\mu'_{2,2} \sigma - \mu'_{2,1} \tau') -
(\lambda'_{2,2} \rho - \lambda'_{2,1} \tau') (\mu_{2,2} \sigma - \mu_{2,1} \tau) = 0.
\end{cases}
$$
Those can be written in the form of quadratic equations; for $i = 1,2$:
$$
A_i \rho \sigma + B_i \tau \tau' + C_i \rho \tau - D_i \rho \tau' - E_i \sigma \tau + F_i \sigma \tau' = 0,
$$
where the coefficients are the following (note that they are all $\neq 0$):
$$
A_i := \lambda_{i,2} \mu'_{i,2} - \lambda'_{i,2} \mu_{i,2},
B_i := \lambda_{i,1} \mu'_{i,1} - \lambda'_{i,1} \mu_{i,1},
C_i := \lambda'_{i,2} \mu_{i,1}, D_i := \lambda_{i,2} \mu'_{i,1}, 
E_i := \lambda'_{i,1} \mu_{i,2}, F_i := \lambda_{i,1} \mu'_{i,2}.
$$
We shall see soon (see \ref{subsubsection:3.4.3}) that the quadratic equations for $i = 1$ and for
$i = 2$ are actually proportional, so they define the same quadrics. \\

Near $\tau' = \infty$, \ie\ near $\tau'' = 0$, the equations become:
$$
A_i \rho \sigma \tau'' + B_i \tau + C_i \rho \tau \tau'' - D_i \rho -
E_i \sigma \tau \tau'' + F_i \sigma = 0,
$$
which, near $\tau'' = 0$, is a graph:
$$
\rho = \dfrac{B_i \tau + F_i \sigma - E_i \sigma \tau \tau''}{D_i - (A_i \sigma + C_i \tau) \tau''},
$$
whence the quadrics is nonsingular at points such that $\tau'' = 0$.

% 3.4

\subsection{Some explicit formulas}

They are based on the following particular choices\footnote{Note that our theta functions
differ from those in \cite{JR} by the change of variable $x \leftrightarrow -x$.} of bases
for $V_{2,\rho_i/\sigma_j}$. We set:
$$
\forall i,j = 1,2 \;,\; \forall k = 1,\ldots,4 \;,\;
e_{i,j}^k := \thq(-x/x_k,-x x_k \sigma_j/\rho_i).
$$
(We use the standard convention $\thq(a,b) := \thq(a) \thq(b)$.) Then the linear forms
$u_{i,j},v_{i,j},\ldots$ on $V_{i,j}$ (evaluations at $x_1,x_2,\ldots$) satisfy:
$$
\begin{cases} u_{i,j}(e_{i,j}^1) = 0, \\
  u_{i,j}(e_{i,j}^2) = \thq(-x_1/x_2,-x_1 x_2 \sigma_j/\rho_i), \end{cases}
\begin{cases}  v_{i,j}(e_{i,j}^1) = \thq(-x_2/x_1,-x_1 x_2 \sigma_j/\rho_i), \\
  v_{i,j}(e_{i,j}^2) = 0, \end{cases}
$$
$$
\begin{cases} w_{i,j}(e_{i,j}^1) = \thq(-x_3/x_1,-x_1 x_3 \sigma_j/\rho_i), \\
  w_{i,j}(e_{i,j}^2) = \thq(-x_3/x_2,-x_2 x_3 \sigma_j/\rho_i), \end{cases}
\begin{cases} w'_{i,j}(e_{i,j}^1) = \thq(-x_4/x_1,-x_1 x_4 \sigma_j/\rho_i), \\
  w'_{i,j}(e_{i,j}^2) = \thq(-x_4/x_2,-x_2 x_4 \sigma_j/\rho_i). \end{cases}
$$

Now, applying equalities $w_{i,j} = \lambda_{i,j} u_{i,j} + \mu_{i,j} v_{i,j}$ and
$w'_{i,j} = \lambda'_{i,j} u_{i,j} + \mu'_{i,j} v_{i,j}$ to all $e_{i,j}^k$, we readily get:
\begin{align*}
  \lambda_{i,j} &= \thq^{-1}(-x_1/x_2,x_1 x_2 \sigma_j/\rho_i) \thq(-x_3/x_2,x_2 x_3 \sigma_j/\rho_i), \\
  \mu_{i,j} &= \thq^{-1}(-x_2/x_1,x_1 x_2 \sigma_j/\rho_i) \thq(-x_3/x_1,x_1 x_3 \sigma_j/\rho_i), \\
  \lambda'_{i,j} &= \thq^{-1}(-x_1/x_2,x_1 x_2 \sigma_j/\rho_i) \thq(-x_4/x_2,x_2 x_4 \sigma_j/\rho_i), \\
  \mu'_{i,j} &= \thq^{-1}(-x_2/x_1,x_1 x_2 \sigma_j/\rho_i) \thq(-x_4/x_1,x_1 x_4 \sigma_j/\rho_i).
\end{align*}

% 3.4.1

\subsubsection{Notational conventions}
\label{subsubsection:3.4.1}

For convenience, we shall use the following abreviations:
$$
T(x) := \thq(-x) \text{~and~}
T\left(\dfrac{x,y,\dots}{z,\dots}\right) := \dfrac{T(x) T(y) \dots}{T(z) \dots};
$$
also:
$$
[kl.ji] := x_k x_l \sigma_j/\rho_i, \; k,l = 1,2,3,4, \; i,j = 1,2.
$$
Note that, since $x \thq(1/x) = \thq(x)$, we have $xy = 1 \Rightarrow T(y) = - y T(x)$. In order
to combine that with Fuchs relation $x_1 x_2 x_3 x_4 = \dfrac{\rho_1 \rho_2}{\sigma_1 \sigma_2}$,
we introduce the \emph{complementarity abreviations}:
$$
i,j \in \{1,2\} \underset{def}{\Longrightarrow} \{i,i'\} = \{j,j'\} = \{1,2\} \text{~and~}
k,l \in \{1,2,3,4\}, k \neq l \underset{def}{\Longrightarrow} \{k,l,k',l'\} = \{1,2,3,4\}.
$$
Then $[kl.ji] [k' l'.j' i'] = 1$, whence the \emph{complementarity relations}:
$$
T\left(\dfrac{[kl.ji]}{[k'l'.j'i']}\right) = - [kl.ji] = - x_k x_l \sigma_j/\rho_i.
$$

From \emph {loc. cit.} we immediately draw:
\begin{align*}
  \dfrac{\lambda_{2,j}}{\lambda_{1,j}} = T\left(\dfrac{[23.j2],[12.j1]}{[12.j2],[23.j1]}\right), \qquad &
  \dfrac{\mu_{2,j}}{\mu_{1,j}} = T\left(\dfrac{[13.j2],[12.j1]}{[12.j2],[13.j1]}\right), \\
  \dfrac{\lambda'_{2,j}}{\lambda'_{1,j}} = T\left(\dfrac{[24.j2],[12.j1]}{[12.j2],[24.j1]}\right), \qquad &
  \dfrac{\mu'_{2,j}}{\mu'_{1,j}} = T\left(\dfrac{[14.j2],[12.j1]}{[12.j2],[14.j1]}\right),
\end{align*}
and:
$$
\dfrac{\lambda'_{i,j}}{\lambda_{i,j}} =
T\left(\dfrac{x_4/x_2}{x_3/x_2}\right) \; T\left(\dfrac{[24.ji]}{[23.ji]}\right), \qquad
\dfrac{\mu'_{i,j}}{\mu_{i,j}} =
T\left(\dfrac{x_4/x_1}{x_3/x_1}\right) \; T\left(\dfrac{[14.ji]}{[13.ji]}\right).
$$

% 3.4.2

\subsubsection{The factor $\gamma$}
\label{subsubsection:3.4.2}

As an immediate application, let\footnote{This looks like a cross ratio, and it actually degenerates
into one when $q \to 1$.} $\gamma_0 := T\left(\dfrac{x_4/x_2,x_4/x_1}{x_3/x_2,x_3/x_1}\right)$. Then:
\begin{align*}
  \dfrac{\lambda'_{1,1} \mu'_{2,2}}{\lambda_{1,1} \mu_{2,2}} =
  \gamma_0 T\left(\dfrac{[24.11], [14.22]}{[23.11], [13.22]}\right), \qquad &
  \dfrac{\lambda'_{2,2} \mu'_{1,1}}{\lambda_{2,2} \mu_{1,1}} =
  \gamma_0 T\left(\dfrac{[24.22], [14.11]}{[23.22], [13.11]}\right), \\
  \dfrac{\lambda'_{1,2} \mu'_{2,1}}{\lambda_{1,2} \mu_{2,1}} =
  \gamma_0 T\left(\dfrac{[24.21], [14.12]}{[23.21], [13.12]}\right), \qquad &
  \dfrac{\lambda'_{2,1} \mu'_{1,2}}{\lambda_{2,1} \mu_{1,2}} =
  \gamma_0 T\left(\dfrac{[24.12], [14.21]}{[23.12], [13.21]}\right).
\end{align*}
All second factors of the right hand sides have the form
$T\left(\dfrac{u,v}{v^{-1},u^{-1}}\right) = uv =
x_1 x_2 x_4^2 \dfrac{\sigma_1 \sigma_2}{\rho_1 \rho_2} = \dfrac{x_4}{x_3}$,
so all the left hand sides have the same value $\gamma = \gamma_0 \dfrac{x_4}{x_3} \cdot$ \\

This $\gamma$ is actually the one which appears in the relations of \ref{subsubsection:3.1.2}. Indeed,
by the argument sketched there, we have more generally a linear relation among the evaluations at the $x_k$:
$$
\exists \alpha, \beta, \gamma \in \Cs \;:\; \forall f \in V_{4,(\rho_1 \rho_2)/(\sigma_1\sigma_2)} \;,\;
f(x_4) = \alpha f(x_1) + \beta f(x_2) + \gamma f(x_3).
$$
(This is in essence what was presented as an axiom in \ref{subsubsection:3.1.2}). In particular:
$$
w'_{1,1} w'_{2,2} = \alpha u_{1,1} u_{2,2} + \beta v_{1,1} v_{2,2} + \gamma w_{1,1} w_{2,2} \text{~and~}
w'_{1,2} w'_{2,1} = \alpha u_{1,2} u_{2,1} + \beta v_{1,2} v_{2,1} + \gamma w_{1,2} w_{2,1}.
$$
Substituting
$$
\lambda_{i,j} u_{i,j} + \mu_{i,j} v_{i,j} \text{~for~} w_{i,j} \text{~and~}
\lambda'_{i,j} u_{i,j} + \mu'_{i,j} v_{i,j} \text{~for~} w'_{i,j},
$$
and assuming moreover (SC), so that the families:
$$
(u_{1,1} u_{2,2}, u_{1,1} v_{2,2}, v_{1,1} u_{2,2}, v_{1,1} v_{2,2}) \text{~and~}
(u_{1,2} u_{2,1}, u_{1,2} v_{2,1}, v_{1,2} u_{2,1}, v_{1,2} v_{2,1})
$$
are free in $W$, we get by identification:
\begin{align*}
  \lambda'_{1,1} \lambda'_{2,2} = \alpha + \gamma \lambda_{1,1} \lambda_{2,2}, & \qquad
  \lambda'_{1,2} \lambda'_{2,1} = \alpha + \gamma \lambda_{1,2} \lambda_{2,1}, \\
  \mu_{1,1} \mu_{2,2} = \beta + \gamma \mu_{1,1} \mu_{2,2}, & \qquad
  \mu'_{1,2} \mu'_{2,1} = \beta + \gamma \mu_{1,2} \mu_{2,1}, \\
  \lambda'_{1,1} \mu'_{2,2} = \gamma \lambda_{1,1} \mu_{2,2}, & \qquad
  \lambda'_{1,2} \mu'_{2,1} = \gamma \lambda_{1,2} \mu_{2,1}, \\
  \lambda'_{2,2} \mu'_{1,1} = \gamma \lambda_{2,2} \mu_{1,1}, & \qquad
  \lambda'_{2,1} \mu'_{1,2} = \gamma \lambda_{2,1} \mu_{1,2}.
\end{align*}

\begin{rmk}
  It is perhaps possible to relax here assumption (SC) using an argument of analytic continuation.
\end{rmk}

% 3.4.3

\subsubsection{The two quadrics are the same}
\label{subsubsection:3.4.3}

We prove here that $A_2/A_1 = B_2/B_1 = C_2/C_1 = D_2/D_1 = E_2/E_1 = F_2/F_1$, so that the two quadrics
obtained in \ref{subsubsection:3.3.2} are actually identical. We begin with the four latter quotients,
which are simpler. Let us denote
$\phi := T\left(\dfrac{[12.11],[12.21]}{[12.12],[12.22]}\right) = \phi_1 \phi_2$,
where $\phi_j := T\left(\dfrac{[12.j1]}{[12.j2]}\right)$, $j = 1,2$; and
$\phi' := x_1 x_2 x_3 x_4 \dfrac{\sigma_2 \sigma_1}{\rho_2^2} = \dfrac{\rho_1}{\rho_2} \cdot$
The following equalities are readily checked (using formulas from \ref{subsubsection:3.4.1}):
\begin{align*}
  \dfrac{C_2}{C_1} = \dfrac{\lambda'_{2,2} \mu_{2,1}}{\lambda'_{1,2} \mu_{1,1}} &=
  T\left(\dfrac{[24.22],[12.21],[13.12],[12.11]}{[12.22],[24.21],[12.12],[13.11]}\right) =
  \phi T\left(\dfrac{[24.22],[13.12]}{[24.21],[13.11]}\right) = \phi \phi' \\
  \dfrac{D_2}{D_1} = \dfrac{\lambda_{2,2} \mu'_{2,1}}{\lambda_{1,2} \mu'_{1,1}} &=
  T\left(\dfrac{[23.22],[12.21],[14.12],[12.11]}{[12.22],[23.21],[12.12],[14.11]}\right) =
  \phi T\left(\dfrac{[23.22],[14.12]}{[23.21],[14.11]}\right) = \phi \phi' \\
  \dfrac{E_2}{E_1} = \dfrac{\lambda'_{2,1} \mu_{2,2}}{\lambda'_{1,1} \mu_{1,2}} &=
  T\left(\dfrac{[24.12],[12.11],[13.22],[12.21]}{[12.12],[24.11],[12.22],[13.21]}\right) =
  \phi T\left(\dfrac{[24.12],[13.22]}{[24.11],[13.21]}\right) = \phi \phi' \\
  \dfrac{F_2}{F_1} = \dfrac{\lambda_{2,1} \mu'_{2,2}}{\lambda_{1,1} \mu'_{1,2}} &=
  T\left(\dfrac{[23.12],[12.11],[14.22],[12.21]}{[12.12],[23.11],[12.22],[14.21]}\right) =
  \phi T\left(\dfrac{[23.12],[14.22]}{[23.11],[14.21]}\right) = \phi \phi'.
\end{align*}

Now we go for $A_2/A_1$ and $B_2/B_1$. Let (temporarily) $u := \sigma_j/\rho_i$; then:
\begin{align*}
  \lambda_{i,j} \mu'_{i,j} - \lambda'_{i,j} \mu_{i,j}
  &= T\left(\dfrac{x_3/x_2,x_2 x_3 u}{x_1/x_2,x_1 x_2 u}\right)
  T\left(\dfrac{x_4/x_1,x_1 x_4 u}{x_2/x_1,x_1 x_2 u}\right) -
  T\left(\dfrac{x_4/x_2,x_2 x_4 u}{x_1/x_2,x_1 x_2 u}\right)
  T\left(\dfrac{x_3/x_1,x_1 x_3 u}{x_2/x_1,x_1 x_2 u}\right) \\
  &= \dfrac{\Phi(x_4)}{T\left(x_1/x_2,x_2/x_1,x_1 x_2 u,x_1 x_2 u\right)},
\end{align*}
where $\Phi(x) := T\left(x_3/x_2,x_2 x_3 u\right) T\left(x/x_1,x_1 u x\right) -
T\left(x_3/x_1,x_1 x_3 u\right) T\left(x/x_2,x_2 u x\right)$.

\begin{lem}
  $$
  \Phi(x) = \dfrac{x_3}{x_1} T\left(x_1/x_2,x_1 x_2 u\right) T\left(x/x_3,x_3 u x\right).
  $$
\end{lem}
\Proof
Although this follows from a general three terms relation for Theta products (see remark below), we give
a direct argument. The function $\Phi$ is holomorphic over $\Cs$ and vanishes at $x = x_3$. It moreover
satisfies the functional equation $\sq \Phi = u^{-1} x^{-2} \Phi$, so it takes the form
$\Phi(x) = C T(x/x_3) T(x_3 u x)$ for some $C \in \C$, which can be determined by evaluation at (say)
$x = x_2$:
$$
C = \dfrac{\Phi(x_2)}{T(x_2/x_3) T(x_2 x_3 u)} =
\dfrac{T\left(x_3/x_2,x_2 x_3 u\right) T\left(x_2/x_1,x_1 x_2 u\right)}{T(x_2/x_3) T(x_2 x_3 u)},
$$
because the second term in $\Phi(x_2)$ includes the vanishing factor $T\left(x/x_2,x_2 u x\right)$.
Thus:
$$
C = \dfrac{T\left(x_3/x_2\right)}{T(x_2/x_3)} T\left(x_2/x_1,x_1 x_2 u\right) =
- \dfrac{x_3}{x_2} T\left(x_2/x_1,x_1 x_2 u\right) = \dfrac{x_3}{x_1} T\left(x_1/x_2,x_1 x_2 u\right),
$$
by the (now) customary transformation rules. The stated formula follows.
\Finprcourt

\begin{prop}
  $$
  \lambda_{i,j} \mu'_{i,j} - \lambda'_{i,j} \mu_{i,j} = 
  \dfrac{x_3}{x_1} T\left(\dfrac{x_4/x_3}{x_2/x_1}\right) T\left(\dfrac{x_3 x_4 u}{x_1 x_2 u}\right) =
  \dfrac{x_3}{x_1} T\left(\dfrac{x_4/x_3}{x_2/x_1}\right) T\left(\dfrac{[34.ji]}{[12.ji]}\right).
  $$
\end{prop}
\Proof
Combining the lemma with the formula that precedes it, we get:
$$
\lambda_{i,j} \mu'_{i,j} - \lambda'_{i,j} \mu_{i,j} =
\dfrac{x_3}{x_1} T\left(\dfrac{x_1/x_2,x_1 x_2 u,x_4/x_3,x_3 x_4 u}{x_1/x_2,x_2/x_1,x_1 x_2 u,x_1 x_2 u}\right)
= \dfrac{x_3}{x_1} T\left(\dfrac{x_4/x_3,x_3 x_4 u}{x_2/x_1,x_1 x_2 u}\right) \cdot
$$
\Finprcourt

\begin{cor}
  \begin{align*}
    A_i = \dfrac{x_3}{x_1} T\left(\dfrac{x_4/x_3}{x_2/x_1}\right) T\left(\dfrac{[34.2i]}{[12.2i]}\right),
    & \quad
    B_i = \dfrac{x_3}{x_1} T\left(\dfrac{x_4/x_3}{x_2/x_1}\right) T\left(\dfrac{[34.1i]}{[12.1i]}\right), \\
    A_2/A_1 = T\left(\dfrac{[34.22],[12.21]}{[34.21],[12.22]}\right), & \quad
    B_2/B_1 = T\left(\dfrac{[34.12],[12.11]}{[34.11],[12.12]}\right).
  \end{align*}
\end{cor}

Using the complementarity relations, we eventually get the identity of our two quadrics:

\begin{cor}
  $$
  A_2/A_1 = B_2/B_1 =
  \dfrac{\rho_1}{\rho_2} T\left(\dfrac{[12.11],[12.21]}{[12.12],[12.22]}\right) =
  C_2/C_1 = D_2/D_1 = E_2/E_1 = F_2/F_1.
  $$
\end{cor}

\begin{rmk}
  The formula of the lemma can also be thought as reflecting the fact that $V_{2,u^{-1}}$ has dimension $2$.
  Along the same lines one can prove a general three terms relation; we exhibit here three different forms
  of the said relation. \\
  $$
  c T(a) T(b/c) T(x/a) T(x/bc) + a T(b) T(c/a) T(x/b) T(x/ac) + b T(c) T(a/b) T(x/c) T(x/ab) = 0.
  $$
  $$
  T(b) T(a/c) T(x/b) T(x/ac) - T(a) T(b/c) T(x/a) T(x/bc) = \dfrac{b}{c} T(c) T(a/b) T(x/c) T(x/ab).
  $$
  $$
  c T(a/b) T(d/c) T(x/ab) T(x/cd) + a T(a/c) T(d/b) T(x/ac) T(x/bd) + d T(c/b) T(a/d) T(x/ad) T(x/bc) = 0.
  $$
  We are not sure whether or not those relations are related to ``Fay's trisecant formula'' (\cf\
  \cite[chapter IIIb]{MumfordTataII}).
\end{rmk}

% 3.4.4

\subsubsection{The discriminant and the singular locus}
\label{subsubsection:3.4.4}

The discriminant of the quadratic form is (for $i = 1,2$):
\begin{align*}
\det \begin{pmatrix} 0 & A_i & C_i & -D_i \\ A_i & 0 & -E_i & F_i \\
  C_i & -E_i & 0 & B_i \\ -D_i & F_i & B_i & 0  \end{pmatrix} &=
A_i^2 B_i^2 + C_i^2 F_i^2 + D_i^2 E_i^2 - 2 A_i B_i C_i F_i - 2 A_i B_i D_i E_i - 2 C_i D_i E_i F_i \\
&= (A_i B_i - C_i F_i - D_i E_i)^2 - 4 C_iD_iE_iF_i
\end{align*}
Using the previously obtained expressions, we find:
\begin{align*}
\text{the above} &=
\left((\lambda_{i,2} \mu'_{i,2} - \lambda'_{i,2} \mu_{i,2})(\lambda_{i,1} \mu'_{i,1} - \lambda'_{i,1} \mu_{i,1}) 
- \lambda'_{i,2} \mu_{i,1} \lambda_{i,1} \mu'_{i,2} - \lambda_{i,2} \mu'_{i,1} \lambda'_{i,1} \mu_{i,2}\right)^2 -
4 \lambda'_{i,2} \mu_{i,1} \lambda_{i,2} \mu'_{i,1} \lambda'_{i,1} \mu_{i,2} \lambda_{i,1} \mu'_{i,2} \\
&= (- \lambda'_{i,2} \mu_{i,2} \lambda_{i,1} \mu'_{i,1} - \lambda_{i,2} \mu'_{i,2} \lambda'_{i,1} \mu_{i,1})^2 -
4 \lambda'_{i,2} \mu_{i,1} \lambda_{i,2} \mu'_{i,1} \lambda'_{i,1} \mu_{i,2} \lambda_{i,1} \mu'_{i,2} \\
&= (\lambda'_{i,1} \lambda_{i,2} \mu_{i,1} \mu'_{i,2} - \lambda_{i,1} \lambda'_{i,2} \mu'_{i,1} \mu_{i,2})^2 
= \Delta_i^2,
\end{align*}
where we have set
$\Delta_i := \lambda'_{i,1} \lambda_{i,2} \mu_{i,1} \mu'_{i,2} - \lambda_{i,1} \lambda'_{i,2} \mu'_{i,1} \mu_{i,2}$.

\begin{rmk}
  \label{rem:quadraticequationcanbewritten}
  The quadratic equation can be written:
  $$
  (A_i \rho - E_i \tau + F_i \tau') (A_i \sigma + C_i \tau - D_i \tau') +
  (E_i C_i \tau^2 + D_i F_i \tau'^2 + (A_i B_i - C_i F_i - D_i E_i) \tau \tau') = 0,
  $$
  \ie\ $X Y + Z Z' = 0$, where
  $X := A_i \rho - E_i \tau + F_i \tau'$, $Y := A_i \sigma + C_i \tau - D_i \tau'$,
  $Z := \lambda_{i,1} \lambda'_{i,2} \tau - \lambda'_{i,1} \lambda_{i,2} \tau'$ and
  $Z' := \mu_{i,1} \mu'_{i,2} \tau - \mu'_{i,1} \mu_{i,2} \tau'$. We shall see next
  (see corollary \ref{cor:Delta-i-nonnul} further below) that $\Delta_i \neq 0$,
  so that $X,Y,Z,Z'$ can be taken as coordinates.
\end{rmk}

Now we use again the explicit formulas. We write $\Delta_i$ as $t_1 - t_2$ (``first and second term'')
where:
\begin{align*}
  t_1 &= \lambda'_{i,1} \lambda_{i,2} \mu_{i,1} \mu'_{i,2} \\
  &=
  T\left(\dfrac{x_3/x_2,x_4/x_2,x_3/x_1,x_4/x_1}{x_1/x_2,x_1/x_2,x_2/x_1,x_2/x_1}\right)
  T\left(\dfrac{[23.1i],[24.2i],[13.2i],[14.1i]}{[12.1i],[12.1i],[12.2i],[12.2i]}\right), \\
  t_2 &= \lambda_{i,1} \lambda'_{i,2} \mu'_{i,1} \mu_{i,2} \\
  &= 
  T\left(\dfrac{x_3/x_2,x_4/x_2,x_3/x_1,x_4/x_1}{x_1/x_2,x_1/x_2,x_2/x_1,x_2/x_1}\right)
  T\left(\dfrac{[23.2i],[24.1i],[13.1i],[14.2i]}{[12.1i],[12.1i],[12.2i],[12.2i]}\right).
\end{align*}
Setting $u := \sigma_1/\rho_i$ and $v := \sigma_2/\rho_i$, we see that $t_1$ and $t_2$ have a common factor:
$$
\text{common factor} =
T\left(\dfrac{x_3/x_2,x_4/x_2,x_3/x_1,x_4/x_1}{x_1/x_2,x_1/x_2,x_2/x_1,x_2/x_1}\right)
T\left(\dfrac{1}{x_1 x_2 u,x_1 x_2 u,x_1 x_2 v,x_1 x_2 v}\right),
$$
so we compute:
$$
\Delta_i = \text{(common factor)} \times
\left(T(x_2 x_3 u,x_1 x_4 u,x_2 x_4 v,x_1 x_3 v) - T(x_2 x_4 u,x_1 x_3 u,x_2 x_3 v,x_1 x_4 v)\right).
$$
Again, instead of the general three terms relation, we favor a direct argument. The second factor on
the right hand side of $\Delta_i$ is $\Psi(x_4)$, where:
$$
\Psi(x) := T(x_2 x_3 u,x_1 x_3 v) T(x_1 u x,x_2 v x) - T(x_1 x_3 u,x_2 x_3 v) T(x_2 u x,x_1 v x).
$$
This function is holomorphic on $\Cs$, vanishes at $x = x_3$ and satisfies $\sq \Psi = c \Psi$,
where $c := \dfrac{1}{x_1 x_2 u v} \cdot$ From this, $\Psi(x) = C T(x/x_3) T(x_1 x_2 x_3 u v x)$ for
some $C \in \C$ which can be determined through evaluation at $x = 1/(x_2 u)$ (so that the second
term in $\Psi$ vanishes), thus yielding:
\begin{align*}
  \Psi(x)
  &= \dfrac{\Psi(1/(x_2 u))}{T(1/(x_2 x_3 u)) T(x_1 x_3 v)} T(x/x_3) T(x_1 x_2 x_3 u v x) \\
  &= \dfrac{T(x_2 x_3 u,x_1 x_3 v) T(x_1/x_2,v/u)}{T(1/(x_2 x_3 u)) T(x_1 x_3 v)}
  T(x/x_3) T(x_1 x_2 x_3 u v x) \\
  \Longrightarrow \Psi(x_4)
  &= \dfrac{T(x_2 x_3 u,x_1 x_3 v) T(x_1/x_2,v/u)}{T(1/(x_2 x_3 u)) T(x_1 x_3 v)}
  T(x_4/x_3) T(x_1 x_2 x_3 u v x_4) \\
  &= x_2 x_3 v T\left(x_1/x_2,\sigma_1/\sigma_2,x_4/x_3,\rho_{i'}/\rho_i\right).
\end{align*}
For the last equality, we applied the usual transformation rules for $T$ along with the obvious equalities
$u/v = \sigma_1/\sigma_2$ and
$x_1 x_2 x_3 u v x_4 = x_1 x_2 x_3 x_4 \sigma_1 \sigma_2/\rho_i^2 = \rho_{i'}/\rho_i$
(by Fuchs relation). In the end, we get one of the possible totally factored forms for (the square root of)
the discriminant:

\begin{prop}
  \label{prop:Delta-i-explicite}
  $$
  \Delta_i =
        [23.2i] \, T\left(\dfrac{x_3/x_2,x_4/x_2,x_3/x_1,x_4/x_1,x_4/x_3}{x_1/x_2,x_2/x_1,x_2/x_1}\right)
  T\left(\dfrac{\sigma_1/\sigma_2,\rho_{i'}/\rho_i}{[12.1i],[12.1i],[12.2i],[12.2i]}\right).
$$
\end{prop}

\begin{cor}
  \label{cor:Delta-i-nonnul}
  $\Delta_i \neq 0$.
\end{cor}

% 3.4.5

\subsubsection{The locus $\det M = 0$}

We describe the image by $[\rho]$ of the locus $\det M = 0$ inside $V^{(*)}$. So let
$M := \begin{pmatrix} m_{1,1} & m_{1,2} \\ m_{2,1} & m_{2,2} \end{pmatrix}$ such that $m_{i,j} \neq 0$
and $\sq m_{i,j}/m_{i,j} = (\rho_i/\sigma_j) x^{-2}$, $i,j = 1,2$; and $m_{1,1} m_{2,2} = m_{1,2} m_{2,1}$.

\begin{lem}
  (i) There are $C_{i,j} \in \Cs$ and $a_i,b_j \in \Cs$, $i,j = 1,2$, such that
  $m_{i,j} = C_{i,j} \thq(x/a_i) \thq(x/b_j)$. \\
  (ii) One has $a_i b_j = \rho_i/\sigma_j$, $i,j = 1,2$, and $C_{1,1} C_{2,2} = C_{1,2} C_{2,1}$.
\end{lem}
\Proof
  (ii) flows directly from (i) by the usual arguments, plus the fact that here:
$$
\det M = (C_{1,1} C_{2,2} - C_{1,2} C_{2,1}) \thq(x/a_1) \thq(x/a_2) \thq(x/b_1) \thq(x/b_2).
$$
  (i) We first see that each $m_{i,j}$ has a common zero with his line neighbour $m_{i,j'}$ and his
  column neighbour $m_{i',j}$ (recall the ``complementarity conventions'' $1' := 2$, $2' := 1$).
  Indeed, the relation $m_{1,1} m_{2,2} = m_{1,2} m_{2,1}$ entails (for instance)
  $m_{1,2}/m_{1,1} = m_{2,2}/m_{2,1}$. The left hand side has poles among the zeroes of $m_{1,1}$, the
  product of which must be $\equiv \rho_1/\sigma_1$, so they cannot be the same as the poles of the
  right hand side, so there must be simplifications, \ie\ $m_{1,2},m_{1,1}$ have a common zero (no more
  for exactly the same reason) and likewise $m_{2,2},m_{2,1}$ have a common zero. \\
  Then we see that there cannot be a common zero to all $m_{i,j}$; indeed, if $a$ was such a zero, the
  same argument applied to $\dfrac{1}{\thq(-x/a)} M$ would lead (among the same lines) to a contradiction.
  The conclusion easily follows.
\Finpr

\begin{prop}
  \label{prop:locusdetnul}
  The image by $[\rho]$ of the locus $\det M = 0$ inside $V^{(*)}$ can be parameterized (up to the
  $\Cs$ action on $\left(\Pr^1(\C)\right)^4$, by $\Cs \ni t \mapsto (f_1(t),f_2(t),f_3(t),f_4(t))$,
  where $f_k(t) := \dfrac{\thq(\sigma_2 x_k t)}{\thq(\sigma_1 x_k t)}$, $k = 1,2,3,4$.
\end{prop}
\Proof
In the direct sense, it is easy to check that any such $(f_1(t),f_2(t),f_3(t),f_4(t))$ can indeed be
realized as some $[\rho](M)$. Conversely, first note that $M$ as described in the lemma is equivalent,
modulo the $\Ddc \times \Ddc$-action, to the same with all $C_{i,j} = 1$. So we take it in that form.
Then $m_{1,2}/m_{1,1} = m_{2,2}/m_{2,1} = \thq(x/b_2)/\thq(x/b_1)$ has to be a solution of
$\sq f/f = \sigma_1/\sigma_2$, so $b_2/b_1 = \sigma_1/\sigma_2$, so we can define:
$$
t := \dfrac{1}{\sigma_1 b_1} = \dfrac{1}{\sigma_2 b_2} \Longrightarrow
m_{1,2}/m_{1,1} = m_{2,2}/m_{2,1} = \dfrac{\thq(\sigma_2 x t)}{\thq(\sigma_1 x t)},
$$
hence the corresponding $(f_1(t),f_2(t),f_3(t),f_4(t))$ is the image of $M$.
\Finpr

\begin{thm}
  \label{thm:locusdetnul}
  The image by $[\rho]$ of the locus $\det M = 0$ inside $V^{(*)}$ lays inside all the quadrics of
  a two-parameters pencil.
\end{thm}
\Proof
First note that each $f_k(x) := \thq(\sigma_2 x_k x)/\thq(\sigma_1 x_k x)$ satisfies
$\sq f_k = (\sigma_1/\sigma_2) f_k$, so each
$g_{k,l} := f_k f_l \thq(\sigma_1 x_1 x) \thq(\sigma_1 x_2 x) \thq(\sigma_1 x_3 x) \thq(\sigma_1 x_4 x)$
satisfies $\sq g_{k,l} = c x^{-4} g_{k,l}$ with:
$$
c := (\sigma_1/\sigma_2)^2 \dfrac{1}{\sigma_1^4 x_1 x_2 x_3 x_4} = \dfrac{1}{\rho_1 \rho_2 \sigma_1 \sigma_2},
$$
since, by (FR), $x_1 x_2 x_3 x_4 = (\rho_1 \rho_2)/(\sigma_1 \sigma_2)$. For $1 \leq k < l \leq 4$, the
denominator of $f_k f_l$ is chased by the theta product and $g_{k,l} \in \Rwg$. Therefore we get six such
functions $g_{k,l}$ in $V_{4,c}$, which has dimension $4$. Thus, chasing denominators, there is (at least)
a two-dimensional space of linear relations of the form:
$$
A \rho \sigma + B \tau \tau' + C \rho \tau - D \rho \tau' - E \sigma \tau + F \sigma \tau' = 0.
$$
Among them, of course, are those found in \ref{subsubsection:3.3.2} (which amount to the same by
\ref{subsubsection:3.4.3}).
\Finpr

%%%%%%%%%%%%%%%%%%%%%%%%%%%%%%%%%%%%%%%%%%%%%%%%%%%%%%%%%%%%%%%%%%%%%%%%%%%%%

% 4

\section{Algebraic threefold and algebraic surface associated to a quadratic form}

In all this part $\underline{A}=(A,B,C,D,E) \in (\C^*)^4$, We return to the coordinates $\rho_i$
and set:
\[
\Phi_{\underline{A}}:=A\rho_1\rho_2+B\rho_3\rho_4+C\rho_1\rho_3-D\rho_1\rho_4
-E\rho_2\rho_3+F\rho_2\rho_4. 
\]
on $\C^4=U_0$.
We will say that ${\underline{A}}$ is \emph{regular} if the discriminant of the quadratic form 
$\Phi_{\underline{A}}$ is not identically $0$ and \emph{singular} on the contrary. \\

The notations $U_p$ and $\mathcal{U}_p$ come from \ref{subsubsection:NSGQ}.

% 4.1

\subsection{Generalities}
\label{subsection:GeneralitiesAlgebraicthreefold}

\begin{lem}
The rank of the quadratic form $\Phi_{\underline{A}}$ is $4$ or $3$. More precisely: 
\begin{itemize}
\item[(i)]
$\underline{A}$ is regular if and only if the rank of 
  $\Phi_{\underline{A}}$ is $4$, then, up to a linear change of coordinates, the equation of
  $C_{\underline{A}}$ is $XY=ZT$ and the only singular point of $C_{\underline{A}}$ is $(0,0,0,0)$;
  the quadric $Q_{\underline{A}}$ is \emph{smooth}.
\item[(ii)]
$\underline{A}$ is singular if and only if the rank of 
  $\Phi_{\underline{A}}$ is $3$ and, up to a linear change of coordinates, the equation of
  $C_{\underline{A}}$ is $XY=Z^2$.  The singular points are the
  $\lambda (a,b,c,d)\vert \, \lambda\in \C$ for some $(a,b,c,d)\in (\C^*)^4$;
  the quadric $Q_{\underline{A}}$ has a \emph{unique singular point}.
\end{itemize}
\end{lem}

\Proof
\begin{itemize}
\item[(i)] See remark \ref{rem:quadraticequationcanbewritten} and corollary \ref{cor:Delta-i-nonnul}.
\item[(ii)]
  We suppose that $\underline{A}$ is singular, then, up to a linear change of coordinates,
  the equation of $S_{[A]}$ is $XY=Z^2$.  More precisely $Z=\alpha \rho_3+\alpha' \rho_4$, with 
  $$
  Z^2=EC\rho_3^2+DF\rho_4^2+(AB-CF-DE)\rho_3\rho_4.
  $$
  As $EC \neq 0$ and
$DF\neq 0$, we have $\alpha\alpha'\neq 0$. We verify easily that $X,Y,Z$ are independent.
\end{itemize}
\Finpr

We consider the homogeneous form of $\Phi_A$:
\[
\Phi_{\underline{A}_{hom}}=A\rho_1^x\rho_2^x\rho_3^y\rho_4^y+B\rho_1^y\rho_2^y\rho_3^x\rho_4^x+C\rho_1^x\rho_2^y\rho_3^x\rho_4^y-D\rho_1^x\rho_2^y\rho_3^y\rho_4^x
-E\rho_1^y\rho_2^x\rho_3^x\rho_4^y+F\rho_1^y\rho_2^x\rho_3^y\rho_4^x
\]
on $\left((\C^2)^*\right)^4$.
The quadratic form $\Phi_{\underline{A}}$ defines an affine quadratic cone $C_{\underline{A}}$ in 
$\C^4=U_0$, invariant under the action of $\C^*$ defined by
$\rho \mapsto \lambda \rho$. The geometric quotient
$C_{\underline{A}}^*/\C^*$ is a quadric $Q_{\underline{A}}$ in
$\mathcal{U}_0\approx \Pr^3(\C)$. 
% In particular, we have a cone $C_{\underline{A}'}$ and a quadric
% $Q_{\underline{A}'}$.
Similarly, the quadratic form:
\[
\widetilde{\Phi}_{\underline{A}}:=A\tilde\rho_1\tilde\rho_2+B\tilde\rho_3\tilde\rho_4+C\tilde\rho_1\tilde\rho_3-D\tilde\rho_1\tilde\rho_4
-E\tilde\rho_2\tilde\rho_3+F\tilde\rho_2\tilde\rho_4 
\]
defines an affine quadratic cone $\tilde C_{\underline{A}}$ in 
$\C^4=U_\infty$. The geometric quotient
$\tilde C_{\underline{A}}^*/\C^*$ is a quadric $\tilde Q_{\underline{A}}$ in
$\mathcal{U}_\infty\approx \Pr^3(\C)$. \\

The homogeneous form $\Phi_{\underline{A}_{hom}}$ defines a closed hypersurface 
$S_{\underline{A}}$ in $\left(\Pr^1(\C)^4\right)$ invariant by $\C^*$. This hypersurface contains
the quadratic cones $C_{\underline{A}}$ and $\tilde{C}_{\underline{A}}$, more precisely it is the
Zariski closure of each cone. 

The image of $S_{\underline{A}}\setminus \Theta_4$ in the quotient $\mathcal{K}_{\, 4}$ is denoted 
$\mathcal{S}_{\underline{A}}$. It is an algebraic surface. It contains the quadrics
$Q_{\underline{A}}$ and $\tilde{Q}_{\underline{A}}$ . More precisely, it is the Zariski closure of
each quadric. The quadrics are compact and not equal, therefore $\mathcal{S}_{\underline{A}}$
is \emph{not separated}. \\

For $i=1,2,3,4$, we denote $S_{\underline{A},i}^0=S_{\underline{A}} \cap \{\sigma_i=0\}$,
$S_{\underline{A},i}^\infty=S_{\underline{A}} \cap \{\sigma_i=\infty\}$, and, for $i\neq j$,
$S_{\underline{A},i,j}^{0,\infty}=S_{\underline{A},i}^0 \cap S_{\underline{A},j}^\infty$.

We have:
\[
S_{\underline{A}}=C_{\underline{A}} \cup \bigcup_{i=1,2,3,4} S_{\underline{A},i}^\infty
=\tilde{C}_{\underline{A}} \cup \bigcup_{i=1,2,3,4} S_{\underline{A},i}^0
=C_{\underline{A}} \cup  \tilde{C}_{\underline{A}} \cup \bigcup_{i,j=1,2,3,4; i\neq j} S_{\underline{A},i,j}^{0,\infty},
\]
\[
\mathcal{S}_{\underline{A}}=C_{\underline{A}} \cup \bigcup_{i=1,2,3,4} \mathcal{S}_{\underline{A},i}^\infty
=\tilde{C}_{\underline{A}} \cup \bigcup_{i=1,2,3,4} \mathcal{S}_{\underline{A},i}^0
=Q_{\underline{A}} \cup  \tilde{Q}_{\underline{A}} \cup \bigcup_{i,j=1,2,3,4; i\neq j} \mathcal{S}_{\underline{A},i,j}^{0,\infty}.
\]
We will see later that the $\mathcal{S}_{\underline{A},i,j}^{0,\infty}$ are \emph{points} which do not
belong to $Q_{\underline{A}} \cup  \tilde{Q}_{\underline{A}}$.

We will now describe the $S_{\underline{A},i}^0$,
$S_{\underline{A},i}^\infty$, $\mathcal{S}_{\underline{A},i}^0$, $S_{\underline{A},i}^\infty$ as union of
``simple" pieces. We will see in particular that the $\mathcal{S}_{\underline{A},i}^0$ (resp.
$\mathcal{S}_{\underline{A},i}^\infty$) are the Zariski closures in $\mathcal{S}_{\underline{A}}$
of some projective lines. (They are non separated unions of smooth rational curves.)

We will describe the set of the points of $S_{\underline{A}}$ admitting at least a coordinate
equal to $\infty$. We consider all the possible cases.
\begin{enumerate}
\item 
Exactly one coordinate equal to $\infty$. We can suppose that it is $\rho_4$.
We consider:
\[
\Psi_1:=A\rho_1\rho_2\tilde\rho_4+B\rho_3+C\rho_1\rho_3\tilde\rho_4-D\rho_1
-E\rho_2\rho_3\tilde\rho_4+F\rho_2.
\]
As $\tilde{\rho}_4=0$, we have 
$\Psi_1=B\rho_3-D\rho_1+F\rho_2=0$. This equation defines a plane of 
$\left(\Pr^1(\C)\right)^3\times \{\infty\}$.
The set of the points whose the only $\infty$ coordinate is $\rho_4$ is a plane of
$\C^3\times \{\infty\}$
 
\item \label{enum:twocoordinates}
  Exactly two coordinates equal to $\infty$. We will prove that this cannot happen. We can suppose
  that it is $\rho_3$ and $\rho_4$. We consider:
\[
\Psi_2:=A\rho_1\rho_2\tilde\rho_3\tilde\rho_4+B+C\rho_1\tilde\rho_4-D\rho_1\tilde\rho_3
-E\rho_2\rho_3\tilde\rho_4+F\rho_2\tilde\rho_3.
\]
As $\tilde{\rho}_4=\tilde{\rho}_3=0$, we have $\Psi_2=B=0$. This is impossible.
 \item 
   Exactly three coordinates equal to $\infty$. We can suppose that it is $\rho_2$, $\rho_3$
   and $\rho_4$. We consider:
\[
\Psi_3=A\rho_1\tilde\rho_3\tilde\rho_4+B\tilde\rho_2+C\rho_1\tilde\rho_2\tilde\rho_4-D\rho_1\tilde\rho_2\tilde\rho_3
-E\tilde\rho_4+F\tilde\rho_3.
\]
As $\tilde{\rho}_4=\tilde{\rho}_3=\tilde{\rho}_3=0$, $\Psi_3=0$ for all 
$\rho_1\in \C$.
\item 
Four coordinates equal to $\infty$: this is the point $(\infty,\infty,\infty,\infty) \in \Theta_4$  
\end{enumerate}

We remark that the $4$ lines $L_{i,j,k}^0=\{\sigma_i=\sigma_j=\sigma_k=0\}$, $i<j<k$, and
the $4$ lines $L_{i,j,k}^\infty=\{\sigma_i=\sigma_j=\sigma_k=\infty\}$, $i<j<k$, are contained
in $S_{\underline{A}}$ for all the values of $\underline{A}$. We denote $p_{i,j,k}^0$ (resp.
$p_{i,j,k}^\infty$ the images of the $L_{i,j,k}^0 \setminus \Theta_4$
(resp. $L_{i,j,k}^\infty \setminus \Theta_4$) in the quotient $\mathcal{S}_{\underline{A}}$.
They are points. We have
$p_{i,j,k}^0 \in Q_{\underline{A}}$ and  $ p_{i,j,k}^\infty \in \tilde Q_{\underline{A}}$.

% 4.2

\subsection{Smoothness properties}

We will describe the singular points of $S_{\underline{A}}$ and $\mathcal{S}_{\underline{A}}$
respectively in the cases  $\underline{A}$ \emph{non-singular} and \emph{singular}.
We recall that we supposed $\underline{A} \in (\C^*)^4$.

\begin{thm}
\label{thsmooth}
We suppose that $\underline{A}$ is \emph{non-singular}. Then: 
\begin{itemize}
\item[(i)]
the only singularities of $S_{\underline{A}}$ are the points $(0,0,0,0)$ and 
$(\infty,\infty,\infty,\infty)$. They are \emph{normal} of type $A_1$;
\item[(ii)]
$\mathcal{S}_{\underline{A}}$ is a (non separated) \emph{smooth} algebraic surface.
\end{itemize}
\end{thm}

\Proof

\begin{itemize}
\item[(i)]
The proof is a variant, in an abstract setting, of a proof in \cite[5.2]{JR}.
We have $S_{\underline{A}} \cap U_0=C_{\underline{A}}$ and we know that the only singularity of 
$C_{\underline{A}}$ is $\{(0,0,0,0)\}$. (The situation is similar for the singularities in $U_\infty$,
the only singularity is $\{(\infty,\infty,\infty,\infty)\}$.)
 
We will search singular points with one (at least) coordinate equal to 
$\infty$. We return to the coordinates $\rho_i$ and we check each case (as in \cite[5.2, p. 37]{JR}).

\begin{enumerate}
\item 
Exactly one coordinate equal to $\infty$. We can suppose that it is $\rho_4$.
We consider:
\[
\Psi_1:=A\rho_1\rho_2\tilde\rho_4+B\rho_3+C\rho_1\rho_3\tilde\rho_4-D\rho_1
-E\rho_2\rho_3\tilde\rho_4+F\rho_2.
\]
We compute the gradient $\nabla \Psi_1$ on $\{\tilde\rho_4=0\}$:
\[
\nabla {\Psi_1}_{\{\tilde\rho_4=0\}}=(-D,F,B,A\rho_1\rho_2+C\rho_1\rho_3-E\rho_2\rho_3).
\]
As $D,F,B\neq 0$, $\nabla \Psi_1\neq 0$ on $\{\tilde\rho_4=0\}$, there are no singularity on 
$\{\tilde\rho_4=0\}$ (i.e. $\rho_4=\infty$).
\item 
  Exactly two coordinates equal to $\infty$. We know that this does not happen (\cf\ page
  \pageref{enum:twocoordinates}, at the end of \ref{subsection:GeneralitiesAlgebraicthreefold},
  item 2 of the enumeration). 

\item 
Exactly three coordinates equal to $\infty$. We can suppose that it is $\rho_2$, $\rho_3$ and 
$\rho_4$. We consider:
\[
\Psi_3=A\rho_1\tilde\rho_3\tilde\rho_4+B\tilde\rho_2+C\rho_1\tilde\rho_2\tilde\rho_4-D\rho_1\tilde\rho_2\tilde\rho_3
-E\tilde\rho_4+F\tilde\rho_3.
\]
We compute the gradient $\nabla \Psi_3$ on 
$\{\tilde\rho_2=\tilde\rho_3=\tilde\rho_4=0\}$:
\[
{\nabla {\Psi_3}}_{\{\tilde\rho_2=\tilde\rho_3=\tilde\rho_4=0\}}
=(0,B,F,-E).
\]
As $B,F,E\neq 0$, $\nabla \Psi_3\neq 0$ on $\{\tilde\rho_2=\tilde\rho_3=\tilde\rho_4=0\}$.
There are no singularity on
$\{\tilde\rho_2=\tilde\rho_3=\tilde\rho_4=0\}$. In particular $(0,\infty,\infty,\infty)$
is not singular.
\item 
  Four coordinate equal to $\infty$. The point $(\infty,\infty,\infty,\infty)$ is singular
  (it is clear using the coordinates $\tilde\rho_i$).
\end{enumerate}
\item[(ii)]
  We can consider the smooth algebraic threefold $S_{\underline{A}} \setminus \Theta_4$ as an
  \emph{analytic manifold}. The group $\C^*$ operates on it without fixed point, therefore
  the corresponding action of its its Lie group $\C$ defines a foliation and we can put on
  the quotient 
  $\left(S_{\underline{A}} \setminus \Theta_4\right)/\C^*$ (interpreted as a set) the analytic
  structure of the space of leaves of the foliation. We get an analytic manifold but a priori
  it could be non Haussdorf\footnote{In general spaces of leaves are non Haussdorf and therefore
  quite pathological.}, and in fact it  is the case. This is very briefly sketched in \cite{JR}
  (\cf\ 3.2, page 38), without reference to the problem of separation\footnote{In the usual results
  on the quotient of a manifold by a free $G$-action, the group $G$ is supposed \emph{compact},
  but $\C^*$ is not compact.}. On the other side, there is a structure of analytic space on
  $\mathcal{S}_{\underline{A}}$ deduced from the algebraic quotient structure by GAGA \cite{SerreGAGA}. 

We denote this space by $\mathcal{S}^{an}_{\underline{A}}$. If we could prove that 
$\mathcal{S}^{an}_{\underline{A}}$ and the analytic manifold of leaves coincide, then 
$\mathcal{S}^{an}_{\underline{A}}$ would be smooth and it would follow that 
$\mathcal{S}_{\underline{A}}$ is smooth, using the proposition \ref{propsmoothalgan} recalled below.
Unfortunately we cannot prove \emph{a priori} the equality of the two analytic structures.
Therefore we will use another approach and prove that the (non separated) algebraic surface 
$\mathcal{S}_{\underline{A}}$ is smooth directly "by hand". We will also give an abstract proof
which will allow us to prove \emph{a posteriori} the equality of the two analytic structures. 

The surface $\mathcal{S}_{\underline{A}}$ is the union of the two quadrics 
$Q_{\underline{A}}$ and $\tilde{Q}_{\underline{A}}$ and of a finite set of points. The two quadrics
are smooth, therefore it remains only to look at the points. These points correspond to the images
of the points of $S_{\underline{A}}$ admitting one (and only one) $0$ coordinate and
one (and only one) $\infty$ coordinate. Up to a coordinate permutation, all the cases are
the same, therefore it is sufficient to consider only one case. Such cases already appeared above 
(at the end of \ref{subsubsection:3.3.2}).

We reformulate in $\rho_i$ coordinates:

Near $\rho_4 = \infty$, \ie\ near $\tilde\rho_4= 0$, the equation~
\[
A \rho_1\rho_2 + B \rho_3\rho_4 + C \rho_1\rho_3 - D \rho_1\rho_4 - E \rho_2\rho_3 
+ F \rho_2\rho_4= 0,
\]
becomes:
\[
A \rho_1\rho_2 \tilde\rho_4 + B \rho_3 + C \rho_1 \rho_ 3\tilde\rho_4 
- D \rho_1 - E \rho_2\rho_3  \tilde\rho_4  + F  \rho_2 = 0,
\]
which, near $\tilde\rho_4= 0$, is a graph:
\[
\rho_1 = \dfrac{B  \rho_3 + F \rho_2 - E \rho_2\rho_3  \tilde\rho_4}{D - (A \rho_2 + C_i \rho_3) \tilde\rho_4 },
\]
whence the surface $S_{\underline{A}}$ is non singular at points such that 
$\rho_4 = \infty$, and in particular at the points defined by $\rho_4 = \infty$ and $\rho_2 = 0$.
If $(\rho_1,0,\rho_3,\infty) \in S_{\underline{A}} \setminus \Theta_4$, then
$B \rho_3- D \rho_1=0$, $\rho_3 \neq 0$ and we can write:
\[
\rho_1/\rho_3 = \dfrac{B + F \rho_2/\rho_3 - E \rho_2  \tilde\rho_4}{D - (A \rho_2/\rho_3 + C)\rho_3 \tilde\rho_4 },
\]
Let $p\in \mathcal{K}_4$ be the image of $S_{\underline{A},2,4}^{0,\infty} \setminus \Theta_4$, then,
in a neighborhood of $p$ in $\mathcal{K}_4$,  we can choose as coordinates 
$(\xi_1=\rho_1/\rho_3,\xi_2=\rho_2/\rho_3,\tilde\xi_4=\rho_3\tilde\rho_4)$ and, near the point
$p=(\xi_1=B/D,\xi_2=0,\tilde\xi_4=0)$, the surface
$\mathcal{S}_{\underline{A}}$ is a graph:
\[
\xi_1=  \dfrac{B + F\xi_2 - E \xi_2  }{D - (A \xi_2 + C) \tilde\xi_4 },
\]
therefore it is non singular at $p$.

\end{itemize}
\Finpr

We will give now an ``abstract" proof of the smoothness of $\mathcal{S}_{\underline{A}}$, using
the Luna's slice theorem (\cf\ \cite{Luna})
It will give more information and allow us to compare the structure of analytic manifold associated
by GAGA to the algebraic structure of $\mathcal{S}_{\underline{A}}$ and the structure of analytic
manifold used in \cite{JR} (\cf\ proof of theorem 2.18).

\begin{prop}
\label{proplunaslice}
Let $G$ be a reductive group acting on an affine variety $X$. Let $X' \subset X$ be the subset
of stable points. We suppose that $G$ operates \emph{freely} on $X'$  (i. e. without fixed point).
Then we have a geometric quotient $\pi: X' \rightarrow Y'=X'/G$ and $X'$ is a $G$-principal bundle
over $Y'$. Moreover, if $x\in X'$ is a smooth point, then $\pi (x)$ is a smooth point of $Y'$.
\end{prop}

\Proof
\Cf\ \cite[Proposition 5.7]{Drezet}.
\Finpr

We can apply the above proposition to the $X_p=S_{\underline{A}} \cap U_p$ 
(\cf\ \ref{subsubsection:NSGQ}). 
We have $X'_p=S_{\underline{A}} \cap (U_p \setminus \{p\})$. The $Y'_p$ are smooth and cover
$\mathcal{S}_{\underline{A}}$. Therefore $\mathcal{S}_{\underline{A}}$ is smooth. 

We can consider the analytic manifolds $(S_{\underline{A}} \setminus \Theta_4)^{an}$ and 
$\mathcal{S}_{\underline{A}}^{an}$  associated by GAGA to 
$S_{\underline{A}} \setminus \Theta_4$ and 
$\mathcal{S}_{\underline{A}}$. Then 
$S_{\underline{A}} \setminus \Theta_4)^{an} \rightarrow \mathcal{S}_{\underline{A}}^{an}$ is a locally
trivial analytic principal $G$-bundle and the analytic manifold structure on 
$\mathcal{S}_{\underline{A}}$ defined by the space of leaves (as in \cite{JR}) is the same than 
$\mathcal{S}_{\underline{A}}^{an}$.

We recall the following result. 

\begin{prop}
\label{propsmoothalgan}
Let $X$ be an algebraic variety. Let $X^{an}$ be the associated analytic space. Then $X$ is smooth
if and only if   $X^{an}$ is smooth
\end{prop}

\Proof
It follows from \cite[no 6, prop. 3, p 9]{SerreGAGA}, \cite[no 24, prop. 27, p 39]{SerreGAGA}
and \cite[chapter VIII, \S 5, no 2, cor. of prop. 1]{BAC89}.
\Finpr

\begin{prop}
We suppose that $S_{\underline{A}}$ is \emph{singular}. Then: 
%the singular points of  belongs to a $\C^*$-invariant rational line. More precisely
\begin{itemize}
\item[(i)]
The set of singular points of $S_{\underline{A}} \setminus \Theta_4$ is the orbit: 
$\C^* (a_1,a_2,a_3,a_4)$, for some 
$(a_1,a_2,a_3,a_4) \in (\C^*)^4$.
\item[(ii)]
The surface
$\mathcal{S}_{\underline{A}}$ is a \emph{non-separated} \emph{normal surface} with a \emph{unique
singular point}. This singular point is of type $A_1$.
\end{itemize}
\end{prop}

\Proof
\begin{itemize}
\item[(i)]
%We will work with coordinates $\sigma_i^x$ and  $\sigma_i^y$:
%$\sigma_i=\sigma_i^x/\sigma_i^y$, 
%$\tilde\sigma_i=\sigma_i^y/\sigma_i^x$.
We already know the singularities in $C_{\underline{A}}=S_{\underline{A}} \cap U_0$: 
$\{(0,0,0,0)\} \cup \C^* (a,b,c,d)$\footnote{The situation is similar for the singularities
in $U_\infty$.}.
The points of $S_{\underline{A}}$ which do not belong to 
$C_{\underline{A}} \cup \tilde C_{\underline{A}} \cup \Theta_4$ are smooth: the proof is the same
than in the regular case.
\item[(ii)]
  The quadric $Q_{\underline{A}}$ is a \emph{normal surface} with a \emph{unique singular point}.
  This singular point is of type $A_1$. It is also the unique singular point of 
$\tilde Q_{\underline{A}}$. The points of $\mathcal{S}_{\underline{A}}$ which do not belong to 
  $Q_{\underline{A}} \cup \tilde Q_{\underline{A}}$ are smooth: the proof is the same than in the
  regular case.

\end{itemize}
\Finpr

% 4.3

\subsection{Pencils of quadrics}

We recall basics on pencils of quadrics in $\Pr^3(\C)$. For details and complements, see
\cite{GriffithsHarris}.

We consider two (non colinear) quadratic forms $\Phi_1$ and 
$\Phi_2$ on $\C^4$ and the pencil of quadrics of $\Pr^3(\C)$ defined by: 
$\{\Phi_{\underline{\lambda}}:=\lambda_1\Phi_1+\lambda_2\Phi_2\vert \, (\lambda_1,\lambda_2)\in (C^2)^*\}$.
We set $C_i=\{\Phi_i=0\}$ and we denote $Q_i$ its image in $\Pr^3(\C)$.

The base of the pencil is the curve $\mathcal{B}:=Q_1\cap Q_2$. If $Q$ is a quadric vanishing on 
$\mathcal{B}$, then $Q$ belongs to the pencil.

By definition the \emph{rank} of a pencil is 
$\max_{i=1,2} rank\, \Phi_i$ (it is independant of the generators). If the rank of a quadric of
the pencil is the rank of the pencil, we will say that this quadric is \emph{generic}.

There are three possibilities.
\begin{enumerate}
\item 
  The rank of the pencil is $4$. Then all the quadrics  of the pencils, except a finite number
  (corresponding to a vanishing discriminant: $\text{discr}\, \Phi_{\underline{\lambda}}=0$), are smooth.
\item
  The rank of the pencil is $3$. Then all the quadrics of the pencil, except a finite number,
  have a unique singular point .
\item
The rank of the pencil is $\leq 2$.
\end{enumerate}

In the following, all the pencils of quadrics are in the first case.

%%%%%%%%%%%%%%%%%%%%%%%%%%%%%%%%%%%%%%%%%%%%%%%%%%%%%%%%%%%%%%%%%%%%%%%%%%%%%

% 5

\section{Projective embeddings of $\F$}

In all this part, we will \emph{fix} $\sigma_2$, $\mu_1$, $\mu_2$, $x_1$, $x_2$, $x_3$, $x_4$ and put 
$\sigma_1=\omega \sigma_2$, with $\omega\in \C^*$ arbitrary\footnote{In \cite{JR}, the parameter
is $\kappa_0$: $\omega=\sigma_1/\sigma_2=\kappa_0^2$.}. 
We will suppose (FR) and (NS). We will suppose (NR), except for $(\sigma_1,\sigma_2)$: we will
allow $\omega \in q^\Z$.

For us \emph{embedding} is allways in the sense of algebraic geometry. A morphism 
$f:X \rightarrow Y$ is an embedding if it is injective and if induces an isomorphism between $X$
and its image $f(X)$, the algebraic structure on $f(X)$ being the structure induced by $Y$.
\emph{Be careful:} it is different from the notion of embedding used in \cite{JR}\footnote{In
\cite{JR}, $X$ is a \emph{set}, $Y$ an algebraic variety, and $f$ an injective map whose image
is a locally closed surface of $Y$.}. 

% 5.1

\subsection{Embedding of $\F$ into 
$\mathcal{K}_{\, 4}=\left(\left(\Pr^1(\C)\right)^4\setminus \Theta_4\right)/\C^*$}
\label{subsecembk4}

We recall the two quadratic forms ($i=1,2$, \cf\ \ref{subsubsection:3.3.2}):
\[
A_i \rho \sigma \tau'' + B_i \tau + C_i \rho \tau \tau'' - D_i \rho -
E_i \sigma \tau \tau'' + F_i \sigma = 0,
\]
% mettre un label et une ref
In the following, we will consider only the case $i=1$ and set
$\underline{A}=(A,B,C,D,E,F)=(A_1,B_1,C_1,D_1,E_1,F_1)$.

Then $\underline{A}$ depends on $\omega=\sigma_1/\sigma_2$ and, if necessary, we will denote
$\underline{A}(\omega)$. We fix a value of 
$(\sigma_1,\sigma_2)$ such that $\sigma_1=\sigma_2$ and we denote $\underline{A}'$ the
corresponding value of $\underline{A}$. We denote $\underline{A}'=\underline{A}(1)$.

We return to the coordinates $\rho_i$ and set:
\[
\Phi_{\underline{A}}:=A\rho_1\rho_2+B\rho_3\rho_4+C\rho_1\rho_3-D\rho_1\rho_4
-E\rho_2\rho_3+F\rho_2\rho_4=
A\rho\sigma+B\tau\tau'+C\rho\tau-D\rho\tau'
-E\sigma\tau+F\sigma\tau'
\]

It is easy to compare with the definitions in \cite{JR}. To $T=(T_{i,j})$ of \cite{JR}, we associate 
$\underline{A}_T:=(T_{12},T_{34},T_{13},-T_{14},-T_{23},T_{24})$.

\begin{lem}
We have $Q_{\underline{A}_{T}} =Q_{\underline{A}}$. 
\end{lem}

\Proof
Let $R\subset \mathcal{S}_{\underline{A}}$ be the image of $\F$ in $\mathcal{K}_{\, 4}$.
We have $R \cap U_0 \subset Q_{\underline{A}}$ and, using 
\cite{JR}, $R \cap U_0 \subset Q_{\underline{A}_{T}}$. Moreover,  $R \cap U_0$ is a Zariski open subset of 
$Q_{\underline{A}}$ and, therefore
$\dim (R\cap U_0)=2$. 

We consider $R':=Q_{\underline{A}} \cap Q_{\underline{A}_{T}}$. The quadric 
$Q_{\underline{A}}$ is irreducible, therefore we have two possibilities:
\begin{enumerate}
\item 
$R'=Q_{\underline{A}}$,
\item
$\dim R'\leq 1$.
\end{enumerate}
We have $R\cap U_0\subset R'$ and therefore $\dim R'\geq 2$. Therefore the only possibility
is $1$ and $Q_{\underline{A}_{T}}=Q_{\underline{A}}$.
\Finpr

\begin{rmk}
Using a similar argument, we get another proof of $Q_{\underline{A}_1}=Q_{\underline{A}_2}$.
\end{rmk}

We denote $X$ the Zariski closure in $\left(\Pr^1(\C)\right)^4$ of $[\rho](W)$ (the image by
$[\rho]$ of $\{\det M=0\}$. We have:
\[
X\setminus \Theta_4=X'\cup \bigcup_{i,j,k;1\leq i<j<k\leq 4} (L_{i,j,k}^\infty \setminus \Theta_4)
\cup \bigcup_{i,j,k;1\leq i<j<k\leq 4} (L_{i,j,k}^0 \setminus \Theta_4), 
\]
where $X'$ is an irreducible affine surface of $\C^4$ (\cf\ \cite[5, p 34]{JR}).

We denote $\mathcal{X}$ the geometric quotient $\mathcal{X}:=(X\setminus \Theta_4)/\C^*$
and $\mathcal{X}'$ the image of $\mathcal{X}'$. Then 
$\mathcal{X}$ is the union of the $8$ points $p_{i,j,k}^\infty$, $p_{i,j,k}^0$ and of the
irreducible curve $\mathcal{X}'$. 
% We have $\mathcal{X}' \subset Q_{\underline{A}}$.
\emph{Be careful:} $\mathcal{X}$ is \emph{not separated}.

\begin{rmk}
  We have $\mathcal{X}'\subset Q_{\underline{A}(\omega)}$, for all $\omega \in \C^*$ (\cf\ the theorem
  \ref{thm:locusdetnul}).
The curve $\mathcal{X}$ is the closure of the curve $\mathcal{X}'$ into 
$\mathcal{K}_{\, 4}$. The $4$ points $p_{i,j,k}^0$ belong to the closure ot the curve $\mathcal{X}'$ into
$Q_{\underline{A}(\omega)}$. Similarly $\mathcal{X}'\subset \tilde{Q}_{\underline{A}(\omega)}$ and 
the $4$ points $p_{i,j,k}^\infty$ belong to the closure ot the curve $\mathcal{X}'$ into
$\tilde{Q}_{\underline{A}(\omega)}$.
\end{rmk}

\begin{prop}
\begin{itemize}
\item[(i)]
The surfaces $X$ and $X'$ are independent of $\omega$. We have a parametrization of the curve 
$\mathcal{X}'$ (\cf\ proposition \ref{prop:locusdetnul} and \cite{JR}). This curve is
\emph{irreducible}.
\item[(i)]
  The curve $\mathcal{X}$ is the union of $\mathcal{X}'$ and of the $8$ points $p_{i,j,k}^0$,
  $p_{i,j,k}^\infty$.
\item[(iii)]
For all $\omega \in \C^*$, we have $X\subset S_{\underline{A}(\omega)}$.
\item[(iv)]
For all $\omega \in \C^*$, we have
$X=S_{\underline{A}(\omega)} \cap S_{\underline{A}(1)}$
\item[(iv)]
Let $\omega_0 \in \C^*$, non resonant.
For all $\omega\in \C$, the quadric $Q_{\underline{A}(\omega)}$ belongs to the pencil of quadrics
generated by $Q_{\underline{A}(\omega_0)}$ and $Q_{\underline{A}(1)}$.
The basis of the pencil is the irreducible quadratic curve $\mathcal{X}' \cup 
\bigcup_{1\leq i<j<k\leq 4} \{p_{i,j,k}^0\}$.
\end{itemize}
\end{prop}

\Proof
It follows easily from \cite[theorem 2.15 and 5.1]{JR}.
\Finpr

We denote $\mathcal{S}^*_{\underline{A}(\omega)}:=\mathcal{S}_{\underline{A}(\omega)}\setminus \mathcal{X}$.
We consider the morphism
 $[\sigma]: F \rightarrow \left(\Pr^1(\C)\right)^4$. Composing with the projection, we get a morphism
 $F \rightarrow \mathcal{K}_{\, 4}$, and using the universal property of the categorical quotient 
 $F \rightarrow  \F$, we get a morphism  $\F \rightarrow \mathcal{K}_{\, 4}$.   

If the condition (NR) is satisfied, then $\omega \in \C^*$ is \emph{non-resonant}.

\begin{thm}
We suppose that the conditions (NR) and (NS) are satisfied.
\begin{itemize}
\item[(i)]
The map $[\rho]$ induces an \emph{injective morphism}
$\F \rightarrow \mathcal{S}_{\underline{A}(\omega)}$ and a \emph{bijection}:
\[
\F \rightarrow  \mathcal{S}^*_{\underline{A}(\omega)}
\]
\item[(ii)]
The morphism $\F \rightarrow  \mathcal{S}^*_{\underline{A}(\omega)}$ is an \emph{isomorphism}.
\item[(iii)]
The geometric quotient $\F$ is smooth.
\item[(iv)]
The algebraic variety $\mathcal{S}_{\underline{A}(\omega)}^*$ is \emph{separated}.
\end{itemize}
\end{thm}

\Proof
\begin{itemize}
\item[(i)]
It follows from \cite[theorem 2.15]{JR}.
\item[(ii)]
If follows from (i) and from the proposition \ref{prop:IsomorphyCriterion}.
\item[(iii)]
  It follows immediately from (ii) and the fact that  $\mathcal{S}_{\underline{A}(\omega)}$ is smooth
  (\cf\ the theorem \ref{thsmooth}).
\item[(iv)]
  It follows immediately from (ii) and the fact that $\F$ is quasi-projective and therefore
  separated (\cf\ proposition \ref{prop:quotientV*/H}).
\end{itemize}

\Finpr

The surface $\mathcal{S}_{\underline{A}}$ is not separated, but the surface 
$\mathcal{S}_{\underline{A}}^*$ is separated. It is interesting to understand \emph{why},
looking at the ``puzzle":
\begin{equation}
\label{equapuzzle}
\mathcal{S}^*_{\underline{A}(\omega)}:=\mathcal{S}_{\underline{A}(\omega)}\setminus \mathcal{X}=\left(Q_{\underline{A}(\omega)}\setminus \left(\mathcal{X}' \cup
\bigcup_{1\leq i<j<k\leq 4} \{p_{i,j,k}^0\}\right)\right) \cup 
\bigcup_{i=1,2,3,4} (S_{\underline{A}(\omega),i}^{\infty}\setminus \mathcal{X}). 
\end{equation}

The curves $S_{\underline{A},i}^{\infty}$ and $S_{\underline{A},i}^{0}$ are non separated. If we remove them
from $S_{\underline{A}}^{\infty}$, we get an open set of 
$Q_{\underline{A}}$ of $\mathcal{K}_{\, 4}$ and this set is separated. We set:
\[
\Delta_i^{0}:=  S_{\underline{A},i}^{0} \cap Q_{\underline{A}},  ~~
\Delta_i^{\infty}:=  S_{\underline{A},i}^{\infty} \cap \tilde{Q}_{\underline{A}}.
\]
The $\Delta_i^{0}$ are complete rational curves of $Q_{\underline{A}}$. We can identify 
$Q_{\underline{A}}$ to a quadric of $\Pr^3(\C)$ and then the $\Delta_i^{0}$ are projective lines
in $\Pr^3(\C)$. Similarly the $\Delta_i^{\infty}$ are complete rational curves of 
$\tilde{Q}_{\underline{A}}$.

We have: 
\[
S_{\underline{A},h}^{0}=\Delta_h^{0} \cup \bigcup_{j=1,2,3,4;j\neq h} \mathcal{S}_{\underline{A},h,j}^{0,\infty}.
\]
Therefore, removing $3$ points from the curve $S_{\underline{A},h}^{0}$, we get a separated complete
rational curve. We will see that there is another way, less evident, to remove $3$ points in order
to get a separated complete rational curve: it is to remove the $3$ points which belongs to the
family of the $p_{i,j,k}^0$.

We set $\mathcal{X}_0:=\mathcal{X}' \cup \bigcup_{1\leq i<j<k\leq 4} \{p_{i,j,k}^0\}$. It is a complete
quartic curve contained in $Q_{\underline{A}}$. The lines $\Delta_i^{0}$ cut this curve along $4$ points.
For each $i=1,2,3,4$, three of them belong to the family of the $p_{i,j,k}^0$, and another
\footnote{\emph{A priori} there could exist a double point, but, according to what will follow,
it would imply that $\rho_i=0$ defines a complete rational curve in $\F$ and it is
impossible.} to $\mathcal{X}'$. Similarly, we set 
$\mathcal{X}_\infty:=\mathcal{X}' \cup \bigcup_{1\leq i<j<k\leq 4} \{p_{i,j,k}^\infty\}$. It is a complete
quartic curve contained in $\tilde{Q}_{\underline{A}}$. The lines $\Delta_i^{\infty}$ cut this curve
along $4$ points. For each $i=1,2,3,4$, three of them belong to the family of the $p_{i,j,k}^\infty$,
and another to $\mathcal{X}'$.

Using some results of \cite{ORS}, we will prove later that $\rho_i=0$ and $\rho_i=\infty$ define
algebraic curves in $\F$ which are isomorphic to the affine line $\C$. Therefore
the curve $S^*_{\underline{A},i}$ defined in $\mathcal{S}^*_{\underline{A}}$ by $\rho_i=0$ is also
isomorphic to the affine line $\C$ and it is obtained  from the complete rational curve
$\Delta_i^{0}$ by removing $4$ points and adding $3$ points which does not belong to 
$Q_{\underline{A}} \cup \tilde{Q}_{\underline{A}}$.

% 5.2

\subsection{Embeddings of $\F$ into $\C^6$ and $\C^4$}

We will give an algebraized version of some results of \cite{JR} (\cf\ theorem 2.21 and 5.3).

% 5.2.1

\subsubsection{Embedding into $\C^6$}

We replace $\underline{A}_{T(\kappa_0)}$ and $\underline{A}_{T(1)}$, where $T(\kappa_0)$ and 
$T(1)$ are defined as in \cite{JR}, by $\underline{A}=\underline{A}(\omega)$ and
$\underline{A}'=\underline{A}(1)$. We suppose that $\underline{A}$ is non resonant.

Starting from the data 
$\underline{A},\underline{A}'$, we will define an \emph{affine algebraic variety} 
$\mathcal{Y}\subset \C^6$. We denote 
$\eta_{ij}$, $i,j=1,2,3,4$, $i\neq j$, the coordinates on $\C^6$:
$(\eta_{12},\eta_{13},\eta_{14},\eta_{23},\eta_{24},\eta_{34}) \in \C^6$.

Now, we define $\mathcal{Y}$ by the following equations (\cf\ \cite{JR},
 (2.26a), (2.26b), (2.26c), (2.26d), page $12$):
 
 \begin{equation}
\label{equac61}
 A\eta_{12}+C\eta_{13}-D\eta_{14}-E\eta_{23}+F\eta_{24}
 +B\eta_{34}=0,
 \end{equation}
 
\begin{equation}
\label{equac62}
 A'\eta_{12}+C'\eta_{13}-D'\eta_{14}-E'\eta_{23}+F'\eta_{24}+B'\eta_{34}=1,
 \end{equation}
 
\begin{equation}
\label{equac63}
\eta_{12}\eta_{34}-\eta_{13}\eta_{24}=0.
\end{equation}

\begin{equation}
\label{equac64}
 \eta_{13}\eta_{24}-\eta_{14}\eta_{23}=0.
\end{equation}

The variety $\mathcal{Y}$ is clearly a surface. The divisor at infinity of this surface is
an algebraic curve $\check{\mathcal X}$ which is independent of $\omega$. One can verify that
the curve $\mathcal{X}'$ defined in \ref{subsecembk4} is isomorphic to a Zariski dense open subset
of $\check{\mathcal X}$.

We consider what happens in the $\Pr^5(\C)$ at infinity of $\C^6$. We use coordinates
$\eta_{12}=X_1/T, ...$. The equations \ref{equac63} and \ref{equac64} give two quadrics
$C_1=\{X_1X_2-X_3X_4=0\}$ and $C_2=\{X_3X_4-X_5X_6=0\}$. 

We suppose $X_4\neq 0$. We use coordinates
$\xi_i=X_i/X_4$, $i=1,2,3,5,6$. Then the equations are $\{\xi_1\xi_2-\xi_3=0\}$ and
$\{\xi_3-\xi_5\xi_6=0\}$. The projective space $L$ of dimension $2$ defined by $\xi_1=\xi_3=\xi_5=0$
is contained in the intersectioncof the two quadrics. 
The equations \ref{equac61} and \ref{equac62} give $2$ hyperplanes $H_1$ and $H_2$:
\[  
 A\xi_1+C\xi_2-D\xi_3+F\xi_5
 +B\xi_6=E,
 \]
 \[  
 A'\xi_1+C'\xi_2-D'\xi_3+F'\xi_5
 +B'\xi_6=E',
 \]
 The quadrics $C_1$ and $C_2$ are tangent at the unique point $(0,0,0,0,0)$. As $EE'\neq 0$,
 this point does not belong to the curve $\check{\mathcal X}$. There are similar computations
 in the other charts. Therefore the curve $\check{\mathcal X}$ is \emph{smooth}. \\

Following \cite{JR}, (2.25), we introduce coordinates on
$\left(\Pr^1(\C)\right)^4$:
\[
\eta_{ij}:=\frac{\rho_i\rho_j}{\Phi_{\underline{A}'}(\rho)}, \quad 1\leq i<j\leq 4.
\]

We have, in homogeneous coordinates:
\[
\eta_{12}:=\frac{\rho^x_i\rho^y_j}{\Phi_{\underline{A}'_{hom}}(\rho^x\rho^y)},
\]
and similar expressions for the others $i,j$. Therefore, we have a \emph{morphism} of algebraic
varieties: 
$$
[\eta]: \left(\Pr^1(\C)\right)^6 \setminus S_{\underline{A}'} \rightarrow \C^6.
$$

The image of $[\sigma]: F \rightarrow \left(\Pr^1(\C)\right)^6$ is contained in
$\left(\Pr^1(\C)\right)^6 \setminus S_{\underline{A}'}$, therefore $[\eta]\circ [\sigma]$ defines
a morphism from $F$ to $\C^6$. It induces a morphism $\F \rightarrow \C^6$. By the
definition of the $\eta_{ij}$, the image of this morphism is contained into $\mathcal{Y}$.
Using \cite{JR} (\cf\ Theorem 2.21)), this morphism induces a \emph{bijective} morphism 
$\vartheta: \F \rightarrow \mathcal{Y}$. We will prove below\footnote{Under the
conditions (NR) and (NS).} that $\mathcal{Y}$ is a \emph{normal} algebraic surface, then
we see, using the propostion \ref{prop:IsomorphyCriterion}, that $\vartheta$ is an \emph{isomorphism}.

\begin{prop}
We suppose that the conditions (NR) and (NS)  are satisfied.
\begin{itemize}
\item[(i)]
  The ``natural" morphism $\theta: \F \rightarrow \mathcal{Y}$ (induced by
  $[\eta] \circ [\rho]$) .is an isomorphism.
\item[(ii)]
The algebraic surface $\mathcal{Y}$ is \emph{smooth}.
\end{itemize}
\end{prop}

% 5.2.2

\subsubsection{Embedding into $\C^4$}

For generic values of the parameters\footnote{For a precise formulation, \cf\ \cite{JR}.},
we can eliminate two coordinates and prove that the algebraic surface $\mathcal{Y}$ is isomorphic
to the intersection $\mathcal{Y}'$ of two quadrics in $\C^4$, an \emph{affine Segre surface.}

In coordinates $\eta_{12},\eta_{13},\eta_{14},\eta_{23}$, two equations are (\cf\ \cite{JR}, (2.27)):
\[
u_0 \eta_{12}^2+u_1 \eta_{12}\eta_{13}+u_2 \eta_{12}\eta_{14}+u_3 \eta_{12}\eta_{23}
+u_4 \eta_{14}\eta_{23}+u_5 \eta_{12}=0
\]
\[
v_0 \eta_{13}^2+v_1 \eta_{12}\eta_{13}+v_2 \eta_{13}\eta_{14}+v_3 \eta_{13}\eta_{23}
+v_4 \eta_{14}\eta_{23}+v_5 \eta_{13}=0.
\]
A Segre complete surface is \emph{normal} (\cf\ \ref{subsubdelpezzo} below), therefore $\mathcal{Y}$
and $\mathcal{Y}'$ are normal surfaces.

\begin{prop}
We suppose that the conditions (NR) and (NS)  are satisfied.
\begin{itemize}
\item[(i)]
  There exists an isomorphism $\theta: \F \rightarrow \mathcal{Y}'$ between $\F$
  and an affine Segre surface. 
\item[(ii)]
The algebraic surface $\mathcal{Y}'$ is \emph{smooth}.
\item[(iii)]
The algebraic surface $\F$ is \emph{rational}.
\item[(iv)]
The algebraic surface $\F$ is \emph{affine}.
\end{itemize}
\end{prop}

\smallskip

We remark that $u_0$, $u_1$, $u_2$, $u_3$, $u_4$ and $v_0$, $v_1$, $v_2$, $v_3$, $v_4$ does not
depend on $\omega$: in the corresponding formulas in \cite{JR}, they do not depend on $\kappa_0$. 

We denote $\tilde{X}$ the "divisor at infinity" of the affine surface $\tilde{Y}'$.
\label{divinfty}
It is the intersection of two quadrics of $\Pr^3(\C)$. Their equations are defined respectively by
$(u_0,u_1,u_2,u_3,u_4)$ and $(v_0,v_1,v_2,v_3,v_4)$ and they do not depend on $\omega$.
Then $\tilde{X}$ is a quartic curve independent of $\omega$.  Supposing (NR) and (NS) satisfied,
we conjecture that the two quadrics are tangent at a unique point. In that case, $\tilde{X}$
would be an \emph{irreducible nodal curve} of genus $2$ (\cf\ \ref{subopenprob}).

One can prove that the curve 
$\mathcal{X}'$ defined in \ref{subsecembk4} is isomorphic to a Zariski dense open set of 
$\tilde{X}$.

% 5.3

\subsection{The $\Pi_{k,l}$ and $\Pi'_{k,l}$, $k,l = 1,\ldots,4$, of \cite{ORS}}
\label{subpi}

% 5.3.1

\subsubsection{The $\Pi_{k,l}$ and $\Pi'_{k,l}$}

Let $M_1 = C_1 L_1, M_2 = C_2 L_2 \in \Mat_2(\C)_1$ such that $\Pi(M_1,M_2)$ is defined and write
as before $\begin{pmatrix} f_i \\ g_i \end{pmatrix}$ the columns $C_i$. If $\Gamma = \Diag(c,c')$,
we have:
$$
\Pi(\Gamma M_1,\Gamma M_2) = [c f_1 c' g_2: c f_2 c' g_1] = [f_1 g_2:f_2 g_1] = \Pi(M_1,M_2).
$$
On the other hand, right multiplication acts only on the $L_i$ so clearly
$\Pi(M_1 \Delta,M_2 \Delta) = \Pi(M_1,M_2)$. Similar statements hold for $\Pi'$ and we state:
$$
\forall \Gamma, \Delta \in \Ddc \times \Ddc \;,\; \Pi(\Gamma M_1 \Delta,\Gamma M_2 \Delta) = \Pi(M_1,M_2)
\text{~and~} \Pi'(\Gamma M_1 \Delta,\Gamma M_2 \Delta) = \Pi'(M_1,M_2).
$$
(It is understood implicitly that ``if one side is defined, etc''.) \\

We start by $\Pi_{1,2}$. We have seen that, if $M \in F$, then $\u(M)$ and $\v(M)$ belong to
$\Mat_2(\C)_1$ and they cannot both have a null first line; nor both have a null second line,
and similarly for columns. Therefore $\Pi\left(\u(M),\v(M)\right) \in \Pr^1(\C)$ and
$\Pi'\left(\u(M),\v(M)\right) \in \Pr^1(\C)$ are well defined. We call them $\Pi_{1,2}(M)$ and
$\Pi'_{1,2}(M)$. Since $u(\Gamma M \Delta) = \Gamma \u(M) \Delta$, we see that $\Pi_{1,2}$ and
$\Pi'_{1,2}$ are $H$-invariant on $F$. \\

Using evaluations for $1 \leq k < l \leq l$, we thus obtain twelve $H$-invariant regular maps
$\Pi_{k,l}$ and $\Pi'_{k,l}$ from $F$ to $\Pr^1(\C)$ and, by going to the quotient, twelve regular
maps $\Pi_{k,l}$ and $\Pi'_{k,l}$ from $\F$ to $\Pr^1(\C)$ (with a slight abuse of notations).

We can define maps $\rho'_i$, $i=1,2,3,4$, from the map $\rho'$ (\cf\
\ref{subsubsection:projectivecoordinatesrank1matrices}) as we define the $\rho_i$ from $\rho$.
We have $\Pi_{i,j}=\rho'_i/\rho'_j$ and $\Pi'_{i,j}=\rho_i/\rho_j$.

% 5.3.2

\subsubsection{Parametrizations of $\F$}
\label{subsubqpant}

In \cite{ORS}, we defined some $q$-pants parametrizations (\cf\ 7.2.3). They are injective maps
from a smooth algebraic surface (built using a punctured elliptic curve) to the set of monodromy
data. We can verify easily that if we put on the set of monodromy data the structure $\F$,
then the $q$-pant parametrizations are morphisms. \\

Using the $q$-pant parametrizations, we see easily that the $16$ algebraic curves defined on 
$\F$ by $\{\rho_i=0\}$, $\{\rho_i=\infty\}$, $\{\rho'_i=0\}$, $\{\rho'_i=\infty\}$ are
isomorphic to the affine line $\C$. Therefore the images $S_{\underline{A},i}^{0,*}$ and
$S_{\underline{A},i}^{\infty,*}$ of the eight curves defined by $\{\rho_i=0\}$ and $\{\rho_i=\infty\}$
are also isomorphic to the affine line $\C$.

If $\rho_i=0$ or $\rho_j=\infty$, then $\Pi_{i,j}=0$, therefore the intersection of the two lines
defined by $\Pi_{i,j}^{-1}(0)$ (\cf\ \cite[7.2.5]{ORS}) 
is the point of $\F$ whose image in $S^*_{\underline{A}}$ is $p_{i,j}^{0,\infty}$.

% 5.4

\subsection{The $16$ lines on the Segre surface}

In \cite{ORS}, using the $\Pi_{ij}$ and $\Pi'_{ij}$, we described $16$ particular curves on the
space of monodromy data (\cf\ 7.2.5), related to some properties of \emph{partial
reducibility}\footnote{It is a $q$-analog of a description of the $24$ lines on the Fricke surface
of $P_{VI}$ as partial reducibility loci, \cf\ \cite{ORS}, 7.1.3.}. Using the $\sigma_i$, $\sigma'_i$,
we get an equivalent definition of this set (\cf\ \ref{subpi}). 
On the other side, there exists $16$ lines on a (generic) complete surface. It is natural
to compare the two sets. It is evident that the $8$ curves defined by the $\sigma_i$ are lines
of the Segre surface. It is less clear for the $8$ curves by the $\sigma'_i$. In a first step,
we doubted that these $8$ curves are lines on the Segre surface, until Nalini Joshi and Pieter
Roffelsen announced a proof of this fact in an article in preparation,
``On q-Painlev\'e VI and the geometry of Segre surfaces'' (tentative title).
Then we looked more carefully and noticed that it follows in fact easily from a result of \cite{ORS},
as we will explain.

% 5.4.1

\subsubsection{The $16$ lines}

We recall some formulas from \cite{ORS} (\cf\ 7.2.5, page 1229). 
We set:
\begin{align*}
\mathbf{e}_q^{1;1,2;3}(\underline{\rho},\underline{\sigma},\underline{x}):=\dfrac
{\thq\left(\frac{\sigma_1}{\rho_1} x_1 x_3\right)
\thq\left(\frac{\sigma_1}{\rho_2} x_2 x_3\right)}
{\thq\left(\frac{\sigma_1}{\rho_1} x_2 x_3\right)
\thq\left(\frac{\sigma_1}{\rho_2} x_1 x_3\right)}, \\
\mathbf{e}_q^{2;1,2;3}(\underline{\rho},\underline{\sigma},\underline{x}) :=
\dfrac
{\thq\left(\frac{\sigma_2}{\rho_1} x_1 x_3\right)
\thq\left(\frac{\sigma_2}{\rho_2} x_2 x_3\right)}
{\thq\left(\frac{\sigma_2}{\rho_1} x_2 x_3\right)
\thq\left(\frac{\sigma_2}{\rho_2} x_1 x_3\right)}
\cdot
\end{align*}

If the conditions (NR), (NS) are satisfied, then 
$\mathbf{e}_q^{1;1,2;3}(\underline{\rho},\underline{\sigma},\underline{x}),\mathbf{e}_q^{2;1,2;3}(\underline{\rho},\underline{\sigma},\underline{x}) \neq 0$.
We have:
\begin{equation}
\label{equapiandpi}
%\Pi_{1,2}^{-1}(0) &=(\Pi')_{3,4}^{-1}\left(\mathbf{e}_q^{1;3,4;1}\right), ~~
%\Pi_{1,2}^{-1}(\infty)=(\Pi')_{3,4}^{-1}\left(\mathbf{e}_q^{2;3,4;1}\right) \\
(\Pi')_{1,2}^{-1}(0) =\Pi_{1,2}^{-1}\left(\mathbf{e}_q^{1;1,2;3}\right), ~~
(\Pi')_{1,2}^{-1}(\infty)=\Pi_{1,2}^{-1}\left(\mathbf{e}_q^{2;1,2;3}\right),
\end{equation}
and similar formulas for the others $\Pi_{ij}$ and $\Pi'_{ij}$, and exchanging $\Pi$ and $\Pi'$. \\

We consider the curve of $\F$ defined by $\rho'_1=0$. Using the formula obtained from 
(\ref{equapiandpi}) exchanging $\Pi$ and $\Pi'$, we see that it is contained in:
\[
\rho_1- \alpha \rho_2=0, ~~ \rho_1- \beta \rho_3=0, ~~ \rho_1- \gamma \rho_4=0,
\]
for some $\alpha,\beta,\gamma \in \C^*$.

We consider the image in $\mathcal{Y}'$ of the curve of $\F$ defined by 
$\rho'_1=0$. We denote it abusively by $\{\rho'_1\}=0$. This curve is contained in
\[
\eta_{13}- \alpha  \eta_{23}=0, ~~ \eta_{12}- \beta \eta_{23}=0, ~~  \eta_{13}- \gamma \eta_{34}=0.
\]
We can replace the last equation by the equivalent equation 
$\beta \eta_{13}- \gamma \eta_{14}=0$ and we get:
\[
\eta_{13}- \alpha  \eta_{23}=0, ~~ \eta_{12}- \beta \eta_{23}=0, ~~  \beta \eta_{13}- \gamma \eta_{14}=0.
\]
These hyperplane equations are independent and define a \emph{line} $L$ in $\C^4$. The curve
$\{\rho'_1=0\}$ is contained in $L$.

Using \cite{ORS}, we proved above that  $\{\rho'_1=0\} \subset \F$ is, as an abstract
algebraic curve, an affine line (\cf\ \ref{subsubqpant}), therefore its image in the affine Segre
surface is $L$. We can also verify using the equations.

The proofs are similar for $\{\rho'_i\}=0$, $i=2,3,4$, with small variations.

\begin{prop}
\label{proplinessegre}
The equations:
\[
\sigma_i=0,~\infty, ~~ \sigma'_i=0,~\infty, ~~i,j=1,2,3,4.
\]
define $16$ affine lines on the affine surface Segre surface $\overline{\mathcal{Y}}'$.
\end{prop}

\Cf\ also, for an independent announcement of this result, a lecture by Nalini Joshi at the workshop
``Asymptotics if beyond all-orders phenomena'' (Isaac Newton Institute, Cambridge, 3rd November, 2022):
``Asymptotics of discrete Painlev\'e equations''.

It is well known that there exists at most $16$ lines on a complete Segre surface. If the surface
is smooth, then the $16$ lines are two by two distinct.

\emph{Be careful:} the $16$ lines defined in the above proposition are not necessarily two by
two distinct. The following result (which follows easily from our description) is a first step
in the study of this question. We will return to the problem in \ref{subsubdelpezzo}.

\begin{prop}
  We suppose $\omega\in \C^*$ non-resonant. The $8$ lines of the family $\sigma_i=0,~\infty$ (resp.
  $\sigma'_i=0,~\infty$) on the affine Segre surface $\mathcal{Y}(\omega)$ are \emph{two by two
  distinct}.
\end{prop}

% 5.4.2

\subsubsection{Segre surfaces and del Pezzo surfaces of degree 4}
\label{subsubdelpezzo}

We follow the presentation of del Pezzo surfaces in \cite{DolgachevCAG}.
In that monograph, Dolgachev proposes various definitions of del Pezzo surfaces. We will use
the following definition (\cf\ \cite[8.1.3, Definition 8.1.2]{DolgachevCAG}).

\begin{defn}
\label{defdelpezzo}
A del Pezzo surface is a projective \emph{normal} algebraic surface $X$ satisfying the two
following conditions:
\begin{itemize}
\item 
all singularities are \emph{rational double points};
\item
the canonical sheaf $\omega_X$ is \emph{invertible} and $\omega_X^{-1}$ is \emph{ample}.
\end{itemize}
The \emph{degree} $d$ of a del Pezzo surface $X$ is, by definition\footnote{According to the
usual notations $K_X$ is the canonical divisor.} $d=K_X^2$. 

\end{defn}

A del Pezzo surface of degree $d\geq 3$ can be embedded as a surface of degree $d$ in 
$\Pr^d(\C)$ ($\omega_X^{-1}$ is \emph{very ample}). For $d=1,2$, the ample divisor is not very ample,
but it gives an embedding into a \emph{weighted projective space}.

The Segre surfaces, defined as a complete intersection of two quartics in $\Pr^4(\C)$ are del Pezzo
surfaces of degree $4$. In particular they are \emph{normal}. \emph{Be careful:} the Segre surfaces
are not necessarily smooth.

\smallskip

For an abstract \emph{complete} algebraic variety, we will call \emph{lines} the $(-1)$-rational
curves. For an affine surface in $\C^n$, the lines are the usual ones.

We recall the following result (\cf\ \cite{DolgachevCAG}, Proposition 8.1.18).

\begin{prop}
  The blow up of $N\leq 8$ points in $\Pr^2(\C)$ is a del Pezzo surface if and only if the points are
  in general position.
\end{prop}

The $N$ divisors above the $N$ point on the blow-up are \emph{lines}. 

\emph{Be careful.} A \emph{singular} del Pezzo surface $X$ cannot be described as a
blow-up\footnote{Such a blow-up is smooth.} But a \emph{minimal resolution}\footnote{That is,
there are no $(-1)$-curves in the fibers.} of $X$ can be described as a
blow-up $Y$\footnote{More precisely $X$ is obtained from $Y$ by blowing down some non
singular $(-2)$-rational curves.}. In this case the 
$N$ points are in \emph{almost general position} (\cf\ \cite{DolgachevCAG}.)

A smooth Segre surface is isomorphic to the blow up of $5$ points, in general position, in 
$\Pr^2(\C)$. In this case, ``in general position" is equivalent to ``of which no three are collinear".
This gives a good description of the $16$ lines. They are:
\begin{itemize}
\item 
the $5$ divisors above the $5$ points,
\item
the proper transforms of the $10$ lines of $\Pr^2(\C)$ joining pair of points,
\item
the proper transform of the conic of $\Pr^2(\C)$ through the $5$ points.
\end{itemize}

A nice picture of the incidence graph of the $16$ lines on a smooth Segre surface is in
\cite{DolgachevCAG}: Figure 8.5. (\Cf\ also a picture in the slides of the lecture of
Nalini Joshi quoted just after proposition \ref{proplinessegre}.) \\

We have the following result (\cf\ \cite{DolgachevCAG}, Table 8.6).

\begin{prop}
  The $16$ lines on a (complete) Segre surface are two by two distinct if and only if the surface
  is \emph{smooth}
\end{prop}
If the number of lines is $<16$, it is  $\leq 12$ (\cf\ \cite{DolgachevCAG}, Table 8.6). 

If the $16$ affine lines defined in the Proposition \ref{proplinessegre} are two by two distinct,
then they correspond to the $16$ lines on the complete Segre surface and this surface is smooth.
We conjecture that the lines are distinct for some generic conditions, which could be (NR), (NS)
plus the condition $\mathbf{Hyp}_{48}$ of \cite{ORS}, 7.2.3.  More generally, we conjecture that,
in all cases, all the lines in the complete Segre surfaces are defined by the equations of the
Proposition \ref{proplinessegre}.

The $8$ lines in each family are two by two distinct, but a line in a family could be equal to
a line of the other family. If two lines coincide, then a singularity at infinity
appears\footnote{It is at a nodal point of the quartic curve at infinity.} and there is
at most $12$ distinct lines. 

% 5.5

\subsection{Morphisms from $\F$ to $\left(\Pr^1(\C)\right)^6$ and to $\left(\Pr^1(\C)\right)^3$}

The six functions $\Pi'_{ij}=\rho_{ij}=\rho_i/\rho_j$, $1\leq i<j\leq 4$, define a morphism
from $\F$ to $\left(\Pr^1(\C)\right)^6$ (\cf\ \cite{ORS}). Therefore the three functions
$\rho_{12},~\rho_{23},\rho_{34}$ define a morphism $\phi$ from $\F$ to $\left(\Pr^1(\C)\right)^3$.
According to \cite{JR}, this morphism is injective. It is an \emph{immersion}. A natural question
is to ask if this immersion is an \emph{embedding}. We denote $\mathcal{Z}'$ the image of the
morphism $\phi$ and $\mathcal{Z}$ its Zariski closure in $\left(\Pr^1(\C)\right)^3$. As $\F$
is smooth, using the proposition \ref{prop:IsomorphyCriterion}, $\phi$ is an embedding if and only
if the points of $\mathcal{Z}'$ are smooth points of $\mathcal{Z}$.
The problem seems difficult and we will only give some indications.

In \cite{JR} the authors give an equation $\mathbf{F}=0$ satisfied by $\mathcal{Z}'$ (\cf\ 5.32).
It seems a priori easy to decide, using the gradient of $\mathbf{F}$ if $\{\mathbf{F}=0\}$ is smooth.
However, we will see that the smoothness involves some possible algebraic relations between
$A,B,C,D$ and it is not evident to see if these relations are true or not. With our notations
(that is replacing the $T_{ij}$ by $\underline{A}$), we have:
\begin{equation}
\label{equapimap3}
A \rho_{12}\rho_{23}^2\rho_{34}+ C \rho_{12}\rho_{23}\rho_{34}-D \rho_{12}\rho_{23}
-E\rho_{23}\rho_{34}+ F \rho_{23}+B=0.
\end{equation}
\[
\mathbf{F}(X,Y,Z)=AXY^2Z+ C XYZ-D XY-E YZ+ F Y+B=0.
\]
\[
\mathbf{F}_X=AY^2Z+ C YZ-D Y, ~~ \mathbf{F}_Y=2AXYZ+ C XZ-D X-E Y+ F, ~~
\mathbf{F}_Z=AXY^2+ C XY-E Y.
\]
As $Y=0$ is excluded, the conditions became:
\[
\mathbf{F}=0, ~~ AYZ+CZ-D=0, ~  \mathbf{F}_Y=0, ~~ A XY+ C X-E=0
\]
Using $\displaystyle Z=\frac{D}{AY+C}$ and $\displaystyle X=\frac{E}{AY+C}$, we can eliminate
$X$ and $Z$ from
$\mathbf{F}=0$ and $\mathbf{F}_Y=0$. We get two relations:
\[
CDE+(AY+C)^2(-EY+F)=0, ~~  DE-2DEY(AY+C)+(AY+C)^2(FY+B)=0,
\]
that is
\begin{align*}
P_{\underline{A}}(Y) &:=-A^2EY^3+A(AF-2CE)Y^2+C(2AF-CE)Y+C(DE+CF)=0, \\
Q_{\underline{A}}(Y) &:=AF^2Y^3+A(2CF-2DE+AB)Y^2+C(2AB-2DE+CF)Y+BC^2=0
\end{align*}
Using the resultant of the two cubic polynomials, we get an algebraic relation 
$\mathrm{Res}(P_{\underline{A}},Q_{\underline{A}})=0$
between $A,B,C,D,E$. It can be identically satisfied or only satisfied for some values of
$\underline{A}$. If it is not satisfied, then the surface is smooth in the domain of the
corresponding chart. If it is satisfied, there are a singular point in this domain.

There are similar relations for the other charts. If one relation is satisfied, then the surface
$\{\mathbf{F}=0\}$ is singular. It remains to check if $\mathcal{Z}'$ contains singularities.
It seems difficult to conclude with this method. It is perhaps possible to use some ideas of
\cite{ORS}, in relation with the $16$ ``lines" on $\F$, with explicit computations
allowed by (\ref{equapimap3}).

%%%%%%%%%%%%%%%%%%%%%%%%%%%%%%%%%%%%%%%%%%%%%%%%%%%%%%%%%%%%%%%%%%%%%%%%%%%%%

% 6

\section{The $\G$ model}
\label{section:GModel}

We return here to an alternative description of our monodromy data space, introduced under
the name of $\G$ in \cite[\S 4]{ORS}, see in particular \S 4.4. We begin with a slightly
more abstract point of view, as in section \ref{section:invariantsquotientsJScase} above.
To begin with, assume a dimension $4$ space $W$ together with two isomorphisms:
$$
\mu_+: V_{1,1} \otimes V_{2,2} \simeq W \text{~and~} \mu_-: V_{1,2} \otimes V_{2,1} \simeq W.
$$
In the case $V_{i,j} = V_{2,\rho_i/\sigma_j}$ and $W := V_{4,\frac{\rho_1 \rho_2}{\sigma_1 \sigma_2}}$,
both $\mu_+$ and $\mu_-$ come from the multiplication maps $V_{1,1} \times V_{2,2} \rightarrow W$
and $V_{1,2} \times V_{2,1} \rightarrow W$, defined as $(f,g) \mapsto fg$. However, the induced
map $f \otimes g \mapsto fg$ is injective \iff\ the following special condition (SC) (already
recalled at the beginning of \ref{Generalnotationsandconventions}) holds \cite[\S 4.4]{ORS}:
$$
\text{(SC)} \quad 
\dfrac{\rho_1}{\sigma_1} \not\equiv \dfrac{\rho_2}{\sigma_2}
\quad \text{~and~} \quad 
\dfrac{\rho_2}{\sigma_1} \not\equiv \dfrac{\rho_1}{\sigma_2}.
$$
Indeed, if for exemple $\rho_1/\sigma_1 = \rho_2/\sigma_2$, permuting $m_{1,1}$ and $m_{2,2}$ is
allowed and does not change the image in $W$. So here we work under assumption (SC). Also, for
simplicity of notations, we write $\mu_+(x_{1,1} \otimes x_{2,2})$ as $x_{1,1} x_{2,2}$ and
$\mu_-(x_{1,2} \otimes x_{2,1})$ as $x_{1,2} x_{2,1}$.

% 6.1

\subsection{A bijective morphism}

% 6.1.1

\subsubsection{The regular map}

The product of the composite maps:
$$
\begin{cases}
  V_{1,1} \times V_{2,2} \longrightarrow V_{1,1} \otimes V_{2,2} \overset{\mu_+}{\longrightarrow} W, \\
  (x_{1,1},x_{2,2}) \mapsto x_{1,1} x_{2,2},
\end{cases}
\quad \text{~and~} \quad
\begin{cases}
  V_{1,2} \times V_{2,1} \longrightarrow V_{1,2} \otimes V_{2,1} \overset{\mu_-}{\longrightarrow} W, \\
  (x_{1,2},x_{2,1}) \mapsto x_{1,2} x_{2,1},
\end{cases}
$$
defines a mapping\footnote{The idea is that, if $M = (x_{i,j})_{i,j = 1,2}$, then $\det\ M = y_+ - y_-$.}:
$$
\begin{cases} V := V_{1,1} \times V_{2,2} \times V_{1,2} \times V_{2,1} \rightarrow W \times W, \\
  (x_{1,1},x_{2,2},x_{1,2},x_{2,1}) \mapsto (y_+,y_-) := (x_{1,1} x_{2,2},x_{1,2} x_{2,1}),
\end{cases}
$$
which restricts to a regular morphism $V^{(*)} \rightarrow W^* \times W^*$ (since, in all generality,
$f,g \neq 0$ implies $f \otimes g \neq 0$). The natural $H$-action on the source $V^{(*)}$ (obvious
notations):
$$
(\lambda_i,\mu_j).(x_{i,j}) := (\lambda_i x_{i,j}/\mu_j)
$$
corresponds to the $\Cs$-homothety action by $\dfrac{\lambda_1 \lambda_2}{\mu_1 \mu_2} $ on the
target $W^* \times W^*$. Recalling that $W^* \times W^* \subset (W \times W)^*$, which admits a
geometric quotient $\Pr(W \times W)$ under the $\Cs$-action, we deduce:

\begin{lem}
  The above map induces a regular map $V^{(*)}/H \rightarrow \Pr(W \times W)$.
\end{lem}
\Proof
  This follows from the universal property of the geometric quotient $V^{(*)}/H$, \cf\ proposition
  \ref{prop:quotientV*/H}.
\Finprcourt

% 6.1.2

\subsubsection{The injectivity}

The map $V^{(*)} \rightarrow W^* \times W^*$ has the property that two elements of the source which
are in the same $H$-orbit have their images are in the same $\Cs$-orbit; this is is obvious, because:
$$
(x_{i,j})_{i,j = 1,2} \mapsto (y_+,y_-) \Longrightarrow
(\lambda_i x_{i,j}/\mu_j)_{i,j = 1,2} \mapsto \dfrac{\lambda_1 \lambda_2}{\mu_1 \mu_2} \, (y_+,y_-).
$$
The converse being true, we get:

\begin{lem}
  The above map $V^{(*)}/H \rightarrow \Pr(W \times W)$ is a bijective regular morphism
  from $V^{(*)}/H$ to a subset $\overline{\Sigma}$ of the image $(W^* \times W^*)/\Cs$ of
  $W^* \times W^*$
  in $\Pr(W \times W)$.
\end{lem}
\Proof
Indeed, the converse of the above implication follows readily from the general fact that, for any
two vector spaces $E,F$:
$$
\forall x,x' \in E \setminus \{0\} \;,\; \forall y,y' \in F \setminus \{0\} \;,\;
\left(x' \otimes y' = x \otimes y\right) \Longleftrightarrow
\left(\exists \lambda \in \Cs \;:\; (x',y') = (\lambda x,\lambda^{-1} y)\right).
$$
\Finprcourt

% 6.1.3

\subsubsection{The image $\overline{\Sigma}$ of $V^{(*)}/H$}
\label{subsubsection:TheImageSigmaBarre}

For any two dimensional spaces $E,F$ with respective bases $(e_1,e_2)$ and $(f_1,f_2)$, the space
$E \otimes F$ has a corresponding basis $(e_i \otimes f_j)_{i,j = 1,2}$ and the coordinates in that
basis of $x \otimes y$, $x = x_1 e_1 + x_2 e_2$, $y = y_1 f_1 + y_2 f_2$, are $(x_i y_j)_{i,j = 1,2}$.
The image of $E^* \times F^*$ in $(E \times F)^*$ is therefore the homogeneous quadric hypersurface
$X T = Y Z$ (in adequate coordinates). \\

In particular, call $\Sigma_+$, $\Sigma_-$ the respective images of $V_{1,1} \times V_{2,2}$ and
of $V_{1,2} \times V_{2,1}$ in $W$, so that the image of $V^{(*)}/H$ in $\Pr(W \times W)$ is
$\overline{\Sigma} = (\Sigma_+ \times \Sigma_-)/\Cs$. Note that
$\Sigma_+ \times \Sigma_- \rightarrow \overline{\Sigma}$, being a restriction of the geometric
quotient $(W \times W)^* \rightarrow \Pr(W \times W)$, is itself a geometric quotient. \\

So we may suppose homogeneous coordinates $[X:Y:Z:T:X':Y':Z':T']$ chosen on the projective space
$\Pr(W \times W)$ such that $\overline{\Sigma}$ is the closed subset of $(W^* \times W^*)/\Cs$
defined by equations:
$$
\overline{\Sigma}: (X T = Y Z) \wedge (X' T' = Y' Z').
$$

\begin{prop}
  The (geometric) quotient $V^{(*)}/H$ is isomorphic to $\overline{\Sigma}$, which is a smooth
  quasi-projective variety of dimension $5$.
\end{prop}
\Proof
Clearly $\Sigma$ is closed in $(W^* \times W^*)/\Cs$, whence a quasi-projective variety.
It is the union of $16$ affine charts defined by the non-vanishing of two coordinates:
$$
\left(\{X \neq 0\} \cup \{Y \neq 0\} \cup \{Z \neq 0\} \cup \{T \neq 0\}\right) \times
\left(\{X' \neq 0\} \cup \{Y' \neq 0\} \cup \{Z' \neq 0\} \cup \{T' \neq 0\}\right)
$$
In each of those charts, it is easily checked to be smooth; for instance, in the chart
$X, X' \neq 0$, it writes (in corresponding affine coordinates) $t = y z$, $t' = y' z'$.
The dimension is plainly $5$. \\
The isomorphy statement follows from the fact that a bijective regular morphism with smooth target
is an isomorphism \cite[pp 46-48]{MumfordCPV}.
\Finprcourt

% 6.2

\subsection{The image $\mathcal{G}$ of $\F$}

Here we forget the abstract point of view and take $W := V_{4,\frac{\rho_1 \rho_2}{\sigma_1 \sigma_2}}$
from start and we appeal to the general properties of the spaces $V_{k,c}$ explained in
\ref{subsubsection:thecontext}. \\

To express the conditions $\det\ M(x_k) = 0$, $\det M \neq 0$, defining $F_{R,S,\x}$, we first
introduce four linear forms $\lambda_k: f \mapsto f(x_k)$, $k = 1,\ldots,4$, on $W$ and their
kernels, the hyperplanes $H_k := \Ker\ \lambda_k$. From Fuchs relation and from (NR), one deduces
that the rank of the family $(\lambda_k)_{k = 1,\ldots,4}$ is $3$, so that the linear subspace
$L := H_1 \cap H_2 \cap H_3 \cap H_4$ of $W$ is a line. \\

On the other hand, if $x_0 \in \Cs$ is arbitrary non congruent modulo $q^\Z$ to any of
$x_1,\ldots,x_4$, the corresponding linear form $\lambda_0: f \mapsto f(x_0)$, defines a
hyperplane $H_0 := \Ker\ \lambda_0$ transverse to $L$, \ie\ such that $H_0 \cap L = \{0\}$. \\

We shall therefore consider the following sets:
$$
G := \{(y_+,y_-) \in \Sigma_+ \times \Sigma_- \tq y_+ - y_- \in L \setminus H_0\}
\subset \Sigma_+ \times \Sigma_- 
\text{~and~}
\mathcal{G} := G/\Cs \subset \overline{\Sigma}.
$$
Clearly, $G$ is a subvariety of $\Sigma_+ \times \Sigma_-$ and, by restriction of the geometric
quotient $\Sigma_+ \times \Sigma_- \rightarrow \overline{\Sigma}$ (see
\ref{subsubsection:TheImageSigmaBarre}), the map $G \rightarrow \mathcal{G}$ is itself a geometric
quotient; and $\mathcal{G}$ is a quasi-projective subvariety of $\overline{\Sigma}$. Now recall from
theorem \ref{thm:LaStructureDeF} that $\F_{R,S,\x}$ is a separated subvariety of $V^{(*)}/H$ and a
geometric quotient of $F_{R,S,\x}$.

\begin{prop}
  The morphism $V^{(*)}/H \rightarrow \Sigma$ induces a bijective morphism from $\F = \F_{R,S,\x}$
  to $\mathcal{G}$.
\end{prop}
\Proof
If $(x_{i,j})_{i,j = 1,2} \in V^{(*)}$, the fact that it corresponds to a point of $F_{R,S,\x}$ means:
vanishing of $y_+ - y_- = \det M$ at $x_1,\ldots,x_4$ but not identically, equivalently:
$y_+ - y_- \in L \setminus H_0$. So the image of $\F$ is $\mathcal{G}$.
\Finpr

Unhappily, we are not able to prove that the above is an isomorphism. This, of course, is
equivalent to the target being smooth (in one direction, tautologically; in the other direction,
thanks to the result fro \cite{MumfordCPV} used at the end of
\ref{subsubsection:TheImageSigmaBarre}).

%%%%%%%%%%%%%%%%%%%%%%%%%%%%%%%%%%%%%%%%%%%%%%%%%%%%%%%%%%%%%%%%%%%%%%%%%%%%%

% 7

\section{Conclusion}

% 7.1

\subsection{What we achieved in this article.}

In this article we improved some results of \cite{ORS} and \cite{JR}.

\begin{itemize}
\item 
  We interpreted the space of monodromy data of $q$-$P_{VI}$ as a geometric quotient (in the
  algebraic geometry sense).
\item
We put on the surface\footnote{With our notations, it is $\mathcal{S}^*(\omega)$.} 
$\mathcal{S}^*(\kappa,t_0)$ defined in \cite{JR} (\cf\ Theorem 2.15) a structure of algebraic variety.
We proved that (under conditions (NR) and (NS)) the structure of analytic manifold on this surface
defined in \cite{JR} is the structure of analytic space associated to our algebraic structure in
GAGA sense. We proved that (under conditions (NR) and (NS))  the algebraic surface
$\mathcal{S}^*=\mathcal{S}^*(\kappa,t_0)$ is smooth and that the natural map 
$\F \rightarrow \mathcal{S}^*$ is an isomorphism of algebraic varieties. It follows
that $\F$ is smooth and that $\mathcal{S}^*$ is separated.
\item
In \cite{JR}, the authors introduce an affine Segre surface\footnote{With our notations, it is 
$\mathcal{Y}'$.} $\F(\kappa,t_0)$ in $\C^4$ (an intersection of two quadrics) and define
a bijection between the \emph{set} of monodromy data and $\F(\kappa,t_0)$. We prove that
this bijection corresponds to an isomorphism of algebraic varieties:
$\F \rightarrow \F(\kappa,t_0)$. It follows that the affine Segre surface
$\mathcal{Y}'=\F(\kappa,t_0)$ is smooth and that $\F$ is a rational surface.
\item
  We proved that the $16$ lines on the complete Segre surface $\mathcal{Y}$ are contained in the
  affine Segre surface $\mathcal{Y}'$ and that they correspond to the $16$ ``lines" introduced
  in \cite{ORS} in relation with some questions of partial reducibility (\cf\ 7.2.5).
\item
  In \cite{ORS}, using the maps $\Pi_{i,j}$, we defined two maps from the set of monodromy data
  to respectively $\left(\Pr^1(\C)\right)^6$ and $\left(\Pr^1(\C)\right)^3$. It is proved in \cite{JR}
  that they are \emph{injective}. Then they induce injective morphisms
  $\F \rightarrow \left(\Pr^1(\C)\right)^6$ and
  $\mathcal {F} \rightarrow \left(\Pr^1(\C)\right)^3$. By definition, they are \emph{immersions}.
  We do not know if they are embeddings\footnote{The terminology used in the corresponding
  statement in \cite{JR} is misleading.}. We conjecture that it is the case for the first one
  but we think that it seems dubious for the second.
\item
There is a natural morphism from $\F$ to the surface $\G$ defined in \cite{ORS}. 
It is a bijection. We do not know if it is an isomorphism.
\end{itemize}

As we will explain in the next paragraph, there emerges from our detailed algebraic descriptions
of the space of monodromy data of $q$-$P_{VI}$ an interesting abstract structure (in terms of
symplectic algebraic geometry) which could exist for \emph{all} Painlev\'e equations, continuous
or discrete.

% 7.2

\subsection{Open problems}
\label{subopenprob}

Even if we made obvious progress, some questions about the space of monodromy data of 
$q$-$P_{VI}$ remain open. More generally, in a sharp contrast of what is known on the left side,
the study of the spaces of monodromy data of the discrete Painlev\'e equations is just beginning.
We propose some paths of exploration.

\begin{enumerate}

\item 
\emph{The curve at infinity on the Segre surface.}
In the differential case, the complete Fricke surface of $P_{VI}$  is (in the generic case) a 
smooth cubic surface. The plane at infinity is special, it is a \emph{tritangent plane}. It cuts
the complete surface along $3$ of the $27$ lines. In our case, the projective space  at infinity
of $\C^4$ (a $\Pr^3(\C)$) cut the complete Segre surface $\mathcal{Y}$ along a quartic curve 
$\tilde{\mathcal{X}}$, the intersection of two quadrics in this $\Pr^3(\C)$ (\cf\ the page
\pageref{divinfty}). There is a natural question: what is special in the choice of this $\Pr^3(\C)$
at infinity and what is the corresponding geometric property of the hyperplane section 
$\tilde{\mathcal{X}}$ ? We do the following conjectures (supposing (NR) and (NS) satisfied):
\begin{itemize}
\item 
the two quadrics defining $\tilde{\mathcal{X}}$ are \emph{tangent} at a \emph{unique} point 
$p\in \tilde{\mathcal{X}}$;
\item
  the quartic curve $\tilde{\mathcal{X}}$ is an irreducible nodal curve, it admits a unique node,
  which is a simple node\footnote{If the two quadrics are smooth, then a real model is the Viviani
  window.}, and  is at $p$;
\item
  the genus of $\tilde{\mathcal{X}}$ is $2$ and its Jacobian surface is, up to isogeny, the product
  of two copies of the elliptic curve\footnote{Modulo this conjecture, the curve $\tilde{\mathcal{X}}$
  would be an \emph{Humbert curve}.} $\Eq := \Cs/q^\Z$. (This could be related with the
  parametrization described in the proposition \ref{prop:locusdetnul}.)

\end{itemize}
Then $\tilde{\mathcal{X}}$ would be an \emph{anticanonical cycle} and
$(\mathcal{Y},\tilde{\mathcal{X}})$ would be a \emph{Looijenga pair} (\cf\ the Introduction
of \cite{Loo}), like in the Fricke surface case. Compare with the Conjecture 7.1 of \cite{CMR}. 

Another problem is to identify the Segre surfaces associated to $q$-$P_{VI}$ among all the Segre
surfaces. It could be related to the above conjectures.

\item
\emph{Immersions and bijective morphisms.}
The map from $\F$ to $\left(\Pr^1(\C)\right)^6$ defined by the $\Pi_{ij}$ is an
\emph{immersion}. We conjecture that it is an \emph{embedding}.

The ``natural map'' from $\F$ to the surface $\G$ (defined in \cite[\S 4]{ORS}, see section
\ref{section:GModel} here) is, if conditions (NR), (NS) and (SC) are satisfied, a bijective
morphism. We conjecture that it an isomorphism.

\item
  \emph{Elliptic fibrations.} In \cite{ORS}, we described some `elliptic fibrations" (by punctured
  elliptic curves) on the space of monodromy data, in relation with Mano decompositions (\cf\ 6.6.2).
  These fibrations remain mysterious. Il would be good to interpreted them as ``induced'' by some
  true elliptic fibrations on some surfaces. If one considers these ``elliptic fibrations'',
  it appears that some ``lines'' are missing to get a true elliptic fibration. We propose two
  methods to ``add lines''. The first one seems bettter.
\begin{itemize}
\item
To use a double covering of $\F$.
\item
To try to find a ``good completion" of $\F$, with some lines at infinity.
\end{itemize}

\begin{itemize}
\item[a)]
\emph{Double covering.}
We can consider the double covering of the complete Segre surface $\mathcal{Y}$ ramified above
the curve at infinity $\tilde{\mathcal{X}}$. (It could be in relation with an involution appearing
in \cite{ORS}, 5.1.4.) The preimages of the $16$ lines split into $32$ curves. We can hope that
this double covering is an elliptic surface. We can even be more optimistic\footnote{Compare
with \cite{DolgachevCAG}, Remark 8.6.2.},  and suppose moreover
that it is a K3 surface. 
Such a picture appeared in the heuristic considerations of part 7.3 of \cite{ORS} (with some
confusion beween the surface and a two covering\footnote{Using the results of the present article,
we know now that the surface is \emph{rational}, it is not a K3 surface.}).
\item[b)]
\emph{Good completion.}
We conjecture that the map from $\F$ to $\left(\Pr^1(\C)\right)^6$ defined by the $\Pi_{ij}$
is an \emph{embedding}. Admitting this fact, we denote by $Z$ the image and by $\overline{Z}$ its
Zariski closure in $\left(\Pr^1(\C)\right)^6$. We can hope that $\overline{Z} \setminus Z$ contains
some lines which can help to understand the "elliptic fibrations". Being more optimistic, it is
possible that $\overline{Z}$ is an elliptic surface.
\end{itemize}
\item
\emph{Exceptional cases for the parameters}
It would be interesting to study the cases of resonant values of the parameters and the
reducibility cases\footnote{For this last case, \cf\ \cite{JR}, 4.2.}. In the differential
case, the Riemann-Hilbert map is no longer an analytic isomorphism (it is only proper) and
there are exceptional fibers, above singular points of the characteristic varieties, which
correspond to Riccati solutions. The quotients are no longer good quotients. Parabolic
structures and stacks are needed. T. Mochizuki introduced  $q$-analogs of parabolic structures
in \cite{Mochizuki}. 
\item
  \emph{Equations with Parameters.} For the classical $P_{VI}$ case, there exists a Riemann-Hilbert
  correspondence in family, taking account of the parameters. We plan to extend \cite{ORS} and
  the present work in a Riemann-Hilbert-Birkhoff correspondence in family in a future work.

\item
\emph{Symplectic structures.}
On the left side, there is a symplectic structure on the space of initial conditions, defined
by the $2$-form\footnote{Which is invariant by the equation and whose polar set is the exceptional
divisor \cite{SakaiRat}.} $\frac{dp \wedge dq}{pq}$. We expect the existence of a ``natural"
symplectic structure on the
right side, on the space of monodromy data, defined by a regular non degenerated $2$-form 
$\omega_{\F}$ on $\F$. 

Admitting this, we can transport the symplectic form $\omega_{\F}$ to the affine Segre
surface $\mathcal{Y}'$. We get a $2$-form $\mathcal{Y}'$ denoted $\omega_{\mathcal{Y}'}$.
It extends into a rational $2$-form $\omega_{\mathcal{Y}}$ on the complete Segre surface
$\mathcal{Y}$. The divisor defined by $\omega_{\mathcal{Y}}$ (the canonical divisor) is the
opposite of the divisor defined by the curve at infinity 
$\tilde{\mathcal{X}}=\mathcal{Y} \setminus \mathcal{Y}'$ (more precisely, this curve is
the pole divisor of $\omega_{\mathcal{Y}}$).

We will give an idea of the definition of a non-degenerated regular $2$-form on the affine
Segre surface which could define the symplectic structure. (It imitates, in some sense,
the construction of the symplectic form by Poincar\'e residues in the case of the Fricke surface.)
To simplify we use the coordinates $(x_1,x_3,x_3,x_4)$ on $\C^4$, we denote $f_1$, $f_2$ quadratic
polynomials defining the two quadrics and we denote $S=V(f_1,f_2)$ the affine Segre surface.
We set $\Omega:=dx_1\wedge dx_2\wedge dx_3\wedge dx_4$. As 
$(f_1,f_2)$ is a complete intersection, we can associate to the symbol 
$\left[\begin{matrix}
\Omega \\
f_1f_2
\end{matrix}\right]$ a $(4,2)$ current
$R:=Res_{f_1,f_2} \left[\begin{matrix}
\Omega \\
f_1f_2
\end{matrix}\right]$ on $\C^4$ (\cf\ \cite{RamisRuget}), 
Lemme $17$\footnote{Such residues are generalizations of the residues studied in
\cite{GriffithsHarris}, $5$.}). We can write this current $R=\omega \wedge [S]$,
where $\omega$ is a rational $2$-form and $[S]$ the $(2,2)$ current of integration
on $S$. The restriction $\omega_S$ of $\omega$ to $S$ is regular and non degenerated.
This $2$-form is our candidate for the symplectic form on $S$\footnote{It is easy to
compute it algebraically, using $\Omega=\omega \wedge df_1 \wedge df_2$, on $S$.}. 
\item
  \emph{Generalization to equations in the Murata  list.} It is natural to try to
  extend\footnote{A first step was made by Anton Eloy in his 2016 thesis
  ``Classification et géométrie des équations aux q-différences : étude globale de q-Painlevé,
  classification non isoformelle et Stokes à pentes arbitraires'', to be found at URL
  \verb?https://www.theses.fr/2016TOU30223?,
  for the cases $q-P(A_4)=$ $q$-$P_V$, $q-P(A_5)$ and $q-P(A_6)$.} 
  some results of \cite{ORS} and of the present article to the equation $q-P(A_2)$ and to the
  equations under $q-P(A_3)$ in the Murata list \cite{Murata}\footnote{For the equations 
  $q-P(A_0^*)$ and $q-P(A_1)$ , Murata does not give Lax pairs.}. In the cases of the equations
  under $q-P(A_3)$, it appears an \emph{irregularity} at $0$ or $\infty$ and it is necessary
  to add the $q$-Stokes multipliers to the monodromy data.
  In the differential case all the Painlev\'e character varieties are affine cubic surfaces.
  By analogy, Yousuke Ohyama formulates recently\footnote{Lecture: \emph{Global analysis on
  the Painlev\'e equations}, Conference: Painlev\'e Equations: 
  From Classical to Modern Analysis, October $26$, $2022$, IRMA, Strasbourg.} the following
  problem: Are the other monodromy manifolds of $q$-Painlev\'e equations Segre surfaces ?
  We think that the answer could be negative. As noticed H. Sakai, in the $q$-difference case,
  $q$-$P_{VI}=q-P(A_3)$ ``is not the boss": there are $q$-difference equations above it in the
  Murata list. We quote \cite{JR}: ``We further note that Chekhov et al. \cite{CMR} conjectured
  explicit affine del Pezzo surfaces of degree three as the monodromy manifolds of the
  $q$-Painlev\'e equations 
higher up in Sakai classication scheme \cite{SakaiRat} than $q$-$P_{VI}$" (\cf\ Table $4$).
We think that, on the contrary, for equations \emph{under} $q-P(A_3)$, the spaces of monodromy
data are (singular) affine Segre surfaces.

The maximum number of lines on the spaces of monodromy data seems to be $21$. We propose some
guesses for the number of lines (for generic values of the parameters) for the surfaces above
$q-P(A_3)$: $q-P(A_0^*): 21$, $q-P(A_1):18$, $q-P(A_2):18$, $q-P(A_3):16$.

We conjecture that the spaces of monodromy data of all the $q$-Painlev\'e equations in Murata
list are affine and possess a symplectic form satisfying the conjectures in $7$ (\cf\ $13$ below). 
\item
\emph{Different Lax pairs.}
For a given Painlev\'e \'equation (and more generally for a given space of initial data) there
can exist several Lax pairs. In differential cases, most of Lax pairs are connected by middle
convolutions, Laplace transforms or change of variables.  A classical example is $P_{VI}$:
\begin{itemize}
\item 
  $Y'=A(x,t)Y$, $A$ is rank two and there are four regular singular points, $0, 1, t, \infty$,
  the system is Fuchsian;
\item
$Y' =(A+B/x)Y$,  A, B are rank three and  $A= \mathrm{diag} (0,1,t)$, the system is irregular.
\end{itemize}
These two Lax pairs are ``equivalent by a Laplace transform". For differential $P_V$, there is
a similar example in \cite{Klimes}.

For the $q$-difference case, the equation $q$-$P_{IV}$ studied in \cite{JR2}
is defined by a linear Fuchsian system. It is a particular case of the systems introduced in
\cite{ORS}, 1.3.1, with $A(x,t)$ of degree $3$ in $x$. But this equation is an equation $q-P(A_5)$
in the Murata list and, in the Murata Lax pair, $A(x,t)$ is of degree $2$ in $x$ and
\emph{irregular}. The two spaces of initial data are equal, but the Painlev\'e equations are
different and an equivalence is not obvious.

Each Lax pair gives a space of monodromy data. It is not clear that the ``natural algebraic
structures'' on these spaces are independent of the choice\footnote{The analytic spaces associated
by GAGA are isomorphic.}.
\item
\emph{Character varieties.}
It would be interesting to express the spaces of monodromy data in terms of representations of
some ``fundamental groups'' (or groupoids), up to equivalence. It is related to the understanding
of the ``intermediate singular points", \cf\ the conclusion of \cite{ORS}. This could help to
understand the symplectic structure.
\item
\emph{Confluences.}
There are several types of confluence.
\begin{enumerate}
\item 
  \emph{Confluences in the Murata list.} The degeneration pattern of Murata (\cf\ \cite{Murata},
  Table 2) is a $q$-analog of the Ohyama-Okumura confluence scheme \cite{OhyamaOkumura} and Murata
  describe the confluences from $q-P(A_3)$. 
  In the differential case, Martin Klimes discovered that the confluence $P_{VI} \rightarrow P_V$
  can be translated  into a birational symplectic map from a character variety to the other:
  $\chi_{VI} \rightarrow \chi_V$ (\cf\ \cite{Klimes} and the 2022 submitted article ``Dynamics
  of the fifth Painlev\'e foliation'' by E. Paul and J.P. Ramis) and Klimes-Paul-Ramis conjectured
  that it is the same for all the confluences the Ohyama-Okumura confluence scheme
  \cite{OhyamaOkumura}. It is natural to conjecture that something similar will happen for
the spaces of monodromy data in Murata list. A system of birational transformations in a family
of cubic surfaces and Segre surfaces would appear. 

In the differential case, it is possible to derive the birational maps from very simple natural
maps between the wild fundamental groupoids of the linearized Painlev\'e equations (\cf\ the above
quoted submitted article of E. Paul and J.P. Ramis).
Therefore the problem for the $q$-analogs can be related to the interpretation of the spaces of
monodromy data as spaces of representations of some ``$q$-wild groupoids". But it is perhaps
possible to use directly \cite{JSAIF}.

\item
\emph{Confluence of $q$-$P_{VI}$ towards $P_{VI}$.} Jimbo-Sakai described the confluence of
$q$-$P_{VI}$ towards $P_{VI}$, when $q\rightarrow 1$ (\cf\ \cite{JimboSakai}, $5$). It seems dubious
that it can be translated into a birational map: the number of lines increases by confluence.

\end{enumerate}
\item
\emph{Difference equations.}
In the Sakai list there are also \emph{difference} equations. Using Birkhoff \cite{Birkhoff1}, 
it is possible to define spaces of monodromy data for such equations. An exemple is the difference
family similar to the family (1.1) of \cite{ORS}, considered by Birkhoff. One can try to put a
structure of geometric quotient on such spaces (at least in some cases as the equations in
\cite{ArinkinBorodin}) and to check if it is isomorphic to an affine del Pezzo surface.  
\item
\emph{Elliptic equations.}
According to H. Sakai, ``the boss" is the equation $A_0^{(1)}$ at the top of his classification.
It is an \emph{elliptic equation} (\cf\ \cite{SakaiRat}, 7.4, \cite{KNY}, 5.5, \cite{JoshiNobutaka}
and the preprint by Nalini Joshi ``Elliptic-difference-type Painlev\'e equations'', Sidney University).
For this elliptic equation, the definition of the space of monodromy data is, as far as we know,
not known\footnote{In all the other cases it follows from Birkhoff work \cite{Birkhoff1}.}.
It seems possible to define it using some results of Igor Krichever \cite{Krichever}.
In \cite{CMR}, the authors conjecture that the space of monodromy data of the elliptic
Painlev\'e equation is an affine del Pezzo surface of degree $3$ (\cf\ Table $4$, first line).
Following their description, there is no line at infinity, therefore we expect $27$ lines on
this affine surface \footnote{Their choice of hyperplane at infinity seems mysterious.}.
\item
\emph{Conjecture}
We extend (boldly) the conjecture 7.1 of \cite{CMR} into the following conjecture for \emph{all}
elliptic and multiplicative/additive Painlev\'e equations. 

For these equations Sakai extended the Okamoto definition of spaces of initial conditions and
defined \emph{Okamoto pairs}. Such a pair $(S,\Delta)$ is the data of a generalized Halphen
surface\footnote{Namely a blow-up of $9$ points in $\Pr^2(\C)$ in non-generic position.} $S$
and an anticanonical divisor $\Delta$ satisfying some conditons. In particular, there exists
a rational two-form $\omega_S$ whose restriction to $S\setminus \Delta$ is a symplectic form.
Then, the Riemann-Hilbert correspondence\footnote{In some cases it is necessary to add
\emph{Stokes style elements} to the monodromy data and to condider a \emph{wild} version of
the correspondence.} assign to an Okamoto pair a \emph{Looijenga pair} $(Y,\tilde X)$, where:
\begin{itemize}
\item 
$Y$ is a (possibely singular) del Pezzo surface;
\item
  $\tilde X$ is ``the divisor at infinity" of $Y$, it is an \emph{anti-canonical cycle} in the
  sense of Looijenga  (\cf\ \cite{Loo}, Introduction);
\item
the affine surface\footnote{It is, in some sense, an affine Calabi-Yau surface.} 
$Y \setminus \tilde X $ is the space of monodromy data endowed with an algebraic quotient
structure\footnote{A good quotient for generic values of the parameters.};
\item
to the Looijenga pair is associated a rational two-form $\omega_Y$ whose restriction to 
$Y\setminus \tilde X$ is a symplectic form\footnote{Except at singular points if there are
such points.}, it is defined by \emph{generalized residues} (up to scaling by $\C^*$);
\item
the Riemann-Hilbert correspondence induces a surjective proper analytic map from 
$S\setminus \Delta$ to $Y \setminus \tilde X$ and the two-form $\omega_S$ is the pull-back
of the two-form $\omega_Y$ by this map (up to scaling).
\end{itemize}

For the classical Painlev\'e equations (differential case), the above conjecture is true
($\Delta$ is the polar divisor of the symplectic form defining the Hamiltonian structure).
In every case, the surface $Y$ is a cubic surface and the divisor at infinity $\tilde X$
is a triangle. In the $P_{VI}$ case (for all parameter values) the vertices of the triangle
are smooth points of $Y$, but when one travels to the right in the Ohyama-Okumura confluence
scheme \cite{OhyamaOkumura}, a smooth vertex can become singular or a singular vertex can become
``more singular''. This phenomena is the conductor of the confluence\footnote{Explicit computations
are easy.} (\cf\ the lecture by J.-P. Ramis ``Canonical dynamics on the Character Varieties
of the Painlev\'e Equations'' at the Strasbourg Conference on Painlev\'e Equations, in honor of
Yousuke Ohyama, october 2022).

In \cite{CMR}, in relation with Painlev\'e equations, there appears only del Pezzo surfaces of
degree $3$. (For the cases $A_0^{(1)}$, $A_0^{(1)*}$, the divisor at infinity is an irreducible
nodal curve, for $A_1^{(1)}$, $A_2^{(1)}$ it is a triangle). With $q$-$P_{VI}$ ($A_3^{(1)}$) there
appears surfaces of degree $4$. We do not know if others degrees can appear. The maximum number
of lines on the spaces of monodromy data of Painlev\'e equations (classical or discrete) seems
to be $27$. 

\end{enumerate}
%
%

%%%%%%%%%%%%%%%%%%%%%%%%%%%%%%%%%%%%%%%%%%%%%%%%%%%%%%%%%%%%%%%%%%%%%%%%%%%%%

\bibliographystyle{plain}

\bibliography{ORSGeometry}

\begin{thebibliography}{10}

\bibitem{ArinkinBorodin}
D.~{Arinkin} and A.~{Borodin}.
\newblock {Moduli spaces of \(d\)-connections and difference Painlev\'e
  equations.}
\newblock {\em {Duke Math. J.}}, 134(3):515--556, 2006.

\bibitem{Berge}
Claude Berge.
\newblock {\em The theory of graphs. {Translated} from the 1958 {French}
  edition by {Alison} {Doig}.}
\newblock Mineola, NY: Dover Publications, 2nd printing of the 1962 first
  {English} edition edition, 2001.

\bibitem{Birkhoff1}
George~D. Birkhoff.
\newblock The generalized {R}iemann problem for linear differential equations
  and the allied problems for linear difference and $q$-difference equations.
\newblock {\em Proc. Amer. Acad.}, 49:521--568, 1913.

\bibitem{Borel}
Armand Borel.
\newblock {\em Linear algebraic groups.}, volume 126 of {\em Grad. Texts Math.}
\newblock New York etc.: Springer-Verlag, 2nd enlarged ed. edition, 1991.

\bibitem{BAC89}
Nicolas Bourbaki.
\newblock {\em {\'E}l{\'e}ments de math{\'e}matique. {Alg{\`e}bre} commutative.
  {Chapitres} 8 et 9}.
\newblock Berlin: Springer, reprint of the 1983 original edition, 2006.

\bibitem{Brion}
Michel Brion.
\newblock Introduction to actions of algebraic groups, 2009.

\bibitem{CMR}
Leonid Chekhov, Marta Mazzocco, and Vladimir Rubtsov.
\newblock Quantised {Painlev{\'e}} monodromy manifolds, {Sklyanin} and
  {Calabi}-{Yau} algebras.
\newblock {\em Adv. Math.}, 376:53, 2021.
\newblock Id/No 107442.

\bibitem{Dolgachev}
Igor Dolgachev.
\newblock {\em Lectures on invariant theory}, volume 296 of {\em Lond. Math.
  Soc. Lect. Note Ser.}
\newblock Cambridge: Cambridge University Press, 2003.

\bibitem{DolgachevCAG}
Igor~V. Dolgachev.
\newblock {\em Classical algebraic geometry. {A} modern view}.
\newblock Cambridge: Cambridge University Press, 2012.

\bibitem{Drezet}
Jean-Marc Dr{\'e}zet.
\newblock Luna's slice theorem and applications.
\newblock In {\em Algebraic group actions and quotients. Notes of the 23rd
  autumn school in algebraic geometry, Wykno, Poland, September 3--10, 2000},
  pages 39--89. Cairo: Hindawi Publishing Corporation, 2004.

\bibitem{GriffithsHarris}
Phillip Griffiths and Joseph Harris.
\newblock {\em Principles of algebraic geometry.}
\newblock New York, NY: John Wiley \& Sons Ltd., 2nd ed. edition, 1994.

\bibitem{Humphreys}
James~E. Humphreys.
\newblock {\em Linear algebraic groups. {Corr}. 2nd printing}, volume~21 of
  {\em Grad. Texts Math.}
\newblock Springer, Cham, 1981.

\bibitem{JimboSakai}
Michio {Jimbo} and Hidetaka {Sakai}.
\newblock {A $q$-analog of the sixth Painlev\'e equation.}
\newblock {\em {Lett. Math. Phys.}}, 38(2):145--154, 1996.

\bibitem{JoshiNobutaka}
Nalini Joshi and Nakazono Nobutaka.
\newblock A review of elliptic difference painlev\'e equations (arxiv
  1902.07842v1), 2019.

\bibitem{JR2}
Nalini Joshi and Pieter Roffelsen.
\newblock On the {Riemann}-{Hilbert} problem for a {{\(q\)}}-difference
  {Painlev{\'e}} equation.
\newblock {\em Commun. Math. Phys.}, 384(1):549--585, 2021.

\bibitem{JR}
Nalini Joshi and Pieter Roffelsen.
\newblock On the monodromy manifold of $q$-painlevé vi and its riemann-hilbert
  problem (arxiv 2202.10597), 2022.

\bibitem{KNY}
K.~Kajiwara, M.~Noumi, and Y.~Yamada.
\newblock Geometric aspects of {P}ainlev\'e equations (arxiv 1509.08186v8),
  2017.

\bibitem{KKMSD}
G.~Kempf, F.~Knudsen, D.~Mumford, and Bernard Saint-Donat.
\newblock {\em Toroidal embeddings. {I}}, volume 339 of {\em Lect. Notes Math.}
\newblock Springer, Cham, 1973.

\bibitem{Klimes}
Martin Klimes.
\newblock The wild monodromy of the fifth painlevé equation and its action on
  wild character variety: an approach of confluence (arxiv 1609.05185), 2022.

\bibitem{Krichever}
I.~M. Krichever.
\newblock Analytic theory of difference equations with rational and elliptic
  coefficients and the {Riemann}-{Hilbert} problem.
\newblock {\em Russ. Math. Surv.}, 59(6):1117--1154, 2004.

\bibitem{Loo}
Eduard Looijenga.
\newblock Rational surfaces with an anti-canonical cycle.
\newblock {\em Ann. Math. (2)}, 114:267--322, 1981.

\bibitem{Luna}
Domingo Luna.
\newblock Slices {\'e}tal{\'e}s.
\newblock {\em Bull. Soc. Math. Fr., Suppl., M{\'e}m.}, 33:81--105, 1973.

\bibitem{Mochizuki}
Takuro Mochizuki.
\newblock Doubly periodic monopoles and q-difference modules (arxiv
  1902.03551), 2019.

\bibitem{MumfordCPV}
David Mumford.
\newblock {\em Algebraic {Geometry}. {I}: {Complex} projective varieties.}
\newblock Class. Math. Berlin: Springer-Verlag, reprint of the corr. 2nd print.
  1976 edition, 1995.

\bibitem{MumfordTataII}
David Mumford.
\newblock {\em Tata lectures on theta. {II}: {Jacobian} theta functions and
  differential equations. {With} the collaboration of {C}. {Musili}, {M}.
  {Nori}, {E}. {Previato}, {M}. {Stillman}, and {H}. {Umemura}}.
\newblock Mod. Birkh{\"a}user Classics. Basel: Birkh{\"a}user, reprint of the
  1984 edition edition, 2007.

\bibitem{Murata}
Mikio {Murata}.
\newblock {Lax forms of the \(q\)-Painlev\'e equations.}
\newblock {\em {J. Phys. A, Math. Theor.}}, 42(11):17, 2009.
\newblock Id/No 115201.

\bibitem{NewsteadTata}
P.~E. Newstead.
\newblock {\em Lectures on introduction to moduli problems and orbit spaces},
  volume~51 of {\em Lect. Math. Phys., Math., Tata Inst. Fundam. Res.}
\newblock Springer, Berlin; Tata Inst. of Fundamental Research, Bombay, 1978.

\bibitem{NewsteadHAL}
P.~E. Newstead.
\newblock Geometric invariant theory.
\newblock In {\em Moduli spaces and vector bundles. A tribute to Peter
  Newstead}, pages 99--127. Cambridge: Cambridge University Press, 2009.

\bibitem{OdaTorus}
Tadao Oda.
\newblock {\em Lectures on torus embeddings and applications. ({Based} on joint
  work with {Katsuya} {Miyake}.)}, volume~58 of {\em Lect. Math. Phys., Math.,
  Tata Inst. Fundam. Res.}
\newblock Springer, Berlin; Tata Inst. of Fundamental Research, Bombay, 1978.

\bibitem{OdaConvex}
Tadao Oda.
\newblock {\em Convex bodies and algebraic geometry. {An} introduction to the
  theory of toric varieties}, volume~15 of {\em Ergeb. Math. Grenzgeb., 3.
  Folge}.
\newblock Berlin etc.: Springer-Verlag, 1988.

\bibitem{OhyamaOkumura}
Yousuke Ohyama and Shoji Okumura.
\newblock A coalescent diagram of the {Painlev{\'e}} equations from the
  viewpoint of isomonodromic deformations.
\newblock {\em J. Phys. A, Math. Gen.}, 39(39):12129--12151, 2006.

\bibitem{ORS}
Yousuke Ohyama, Jean-Pierre Ramis, and Jacques Sauloy.
\newblock The space of monodromy data for the {Jimbo}-{Sakai} family of
  {{\(q\)}}-difference equations.
\newblock {\em Ann. Fac. Sci. Toulouse, Math. (6)}, 29(5):1119--1250, 2021.

\bibitem{PopovVinberg}
V.~L. Popov and Eh.~B. Vinberg.
\newblock Invariant theory.
\newblock Algebraic geometry. {IV}: {Linear} algebraic groups, invariant
  theory, {Encycl}. {Math}. {Sci}. 55, 123-278 (1994); translation from {Itogi}
  {Nauki} {Tekh}., {Ser}. {Sovrem}. {Probl}. {Mat}., {Fundam}. {Napravleniya}
  55, 137-309 (1989)., 1989.

\bibitem{RamisRuget}
Jean-Pierre Ramis and Gabriel Ruget.
\newblock R{\'e}sidus et dualit{\'e}.
\newblock {\em Invent. Math.}, 26:89--131, 1974.

\bibitem{SakaiRat}
Hidetaka Sakai.
\newblock Rational surfaces associated with affine root systems and geometry of
  the {Painlev{\'e}} equations.
\newblock {\em Commun. Math. Phys.}, 220(1):165--229, 2001.

\bibitem{JSAIF}
Jacques Sauloy.
\newblock Syst\`emes aux {$q$}-diff\'erences singuliers r\'eguliers:
  classification, matrice de connexion et monodromie.
\newblock {\em Ann. Inst. Fourier (Grenoble)}, 50(4):1021--1071, 2000.

\bibitem{JSENS}
Jacques Sauloy.
\newblock Galois theory of {Fuchsian} {{\(q\)}}-difference equations.
\newblock {\em Ann. Sci. {\'E}c. Norm. Sup{\'e}r. (4)}, 36(6):925--968, 2003.

\bibitem{SerreGAGA}
Jean-Pierre Serre.
\newblock Algebraic geometry and analytic geometry.
\newblock {\em Ann. Inst. Fourier}, 6:1--42, 1956.

\bibitem{SpringerGAIV}
T.~A. Springer.
\newblock Linear algebraic groups.
\newblock Algebraic geometry {IV}: {Linear} algebraic groups, invariant theory,
  {Encycl}. {Math}. {Sci}. 55, 1-121 (1994); translation from {Itogi} {Nauki}
  {Tekh}., {Ser}. {Sovrem}. {Probl}. {Mat}., {Fundam}. {Napravleniya} 55, 5-136
  (1989)., 1989.

\bibitem{SpringerLAG}
T.~A. Springer.
\newblock {\em Linear algebraic groups.}
\newblock Mod. Birkh{\"a}user Classics. Basel: Birkh{\"a}user, reprint of the
  1998 2nd ed. edition, 2009.

\end{thebibliography}

\printindex

\end{document}